\documentclass[12pt]{article}

\usepackage[a4paper,margin=20mm]{geometry}
\usepackage{amsmath,amssymb,amsfonts,amsthm,mathtools}
\usepackage{bm}
\usepackage{enumitem}
\usepackage{microtype}
\usepackage{graphicx}
\usepackage{float}
\usepackage{setspace}
\usepackage{comment}
\usepackage[hidelinks]{hyperref}

\numberwithin{equation}{section}
\setlist{itemsep=0.10em,topsep=0.3em}
\allowdisplaybreaks

\theoremstyle{plain}
\newtheorem{theorem}{Theorem}[section]
\newtheorem{proposition}[theorem]{Proposition}
\newtheorem{lemma}[theorem]{Lemma}
\newtheorem{corollary}[theorem]{Corollary}

\theoremstyle{definition}
\newtheorem{assumption}[theorem]{Assumption}

\newtheorem{remark}[theorem]{Remark}
\newtheorem{example}[theorem]{Example}

\newcommand{\R}{\mathbb R}

\newcommand{\Sph}{\mathbb S^{d-1}}
\newcommand{\Id}{\mathrm I}
\newcommand{\dd}{\mathrm d}
\newcommand{\ee}{\mathrm e}
\newcommand{\dist}{\operatorname{dist}}

\title{Chaotic Dynamics of Charged Particles in Magnetic Fields with Singular Forces}
\author{
Linhao Chen,
		\thanks{
			School of Mathematics, Guangxi University, Nanning 530004, P.R. China.
			}\quad
Yong Li,
		\thanks{
			College of Mathematics, Jilin University, Changchun 130012, PR China.\\
			School of Mathematics and Statistics, and Center for Mathematics and
Interdisciplinary Sciences, Northeast Normal University, Changchun, 130024, PR
China.			}\quad
Weiping Yan,\footnote{E-mail:\,{\tt chenlh@st.gxu.edu.cn (L. Chen)}, {\tt liyong@jlu.edu.cn (Y. Li)}, {\tt yanwp@gxu.edu.cn (W. Yan)}, Corresponding author: Weiping Yan.}
	\thanks{
		School of Mathematics, Center for Applied Mathematics of Guangxi (Guangxi University),
                 Guangxi University, Nanning 530004, P.R. China.}
}
\date{}

\begin{document}
\maketitle

\begin{abstract}
	We study charged-particle dynamics under a non-Lipschitz velocity-dependent force with magnetic and damping effects, focusing on continuation through zero velocity and chaos in the velocity direction.
For a time-independent directional force and fixed parameters, we introduce a $C^1$ regularization. When an unregularized trajectory reaches zero velocity in finite time, inner rescaling yields convergence to a reference orbit within the regularization ball, while the regularized trajectories converge locally uniformly outside it to a unique continuation.
For the directional dynamics, stereographic projection yields a parameterized Euclidean system. By choosing the skew-symmetric magnetic part, we obtain a nonautonomous parameterized center manifold reduction with a prescribed finite Taylor jet. Under an additional factorization assumption, the reduced dynamics admits Shilnikov homoclinic orbits and associated chaos.

\bigskip

\bigskip


\end{abstract}

\tableofcontents

\section{Introduction and Main results}
\subsection{Introduction}
The motion of charged particles in magnetic fields is a classical subject in physics and mathematical physics. The magnetic part of the Lorentz force can bend the trajectory of a moving charged particle, causing it to gyrate about magnetic field lines.
 Poincaré \cite{Poincare1896} investigated the motion of charged particles in an idealized magnetic monopole field and showed that their trajectories lie on a cone. Such motion is relevant to plasmas and the Earth's radiation belts \cite{Kruskal1958,Northrop1961}.
Chaotic bahovior of charged particles also plays an important role in some situations. Numerical calculations of relativistic electrons
interacting with whistler waves indicate that chaotic trajectories
can enhance diffusion through phase space. Under suitable conditions, this process may scatter
electron pitch angles in the Io torus and contribute to auroral
electron precipitation \cite{WykesChapmanRowlands2001}. Magnetic
reconnection changes magnetic field topology and converts magnetic
energy into kinetic energy. During collisionless reconnection, the
in-plane Hall electric fields generated within the process can
energize ions and drive their trajectories into a chaotic regime,
resulting in ion heating \cite{YoonBellan2021}. One can see \cite{AssenzaReberTerek2025,DairbekovPaternainStefanovUhlmann2007,HairerLubichShi2022,HairerLubichWang2020,Jollivet2007,YanTangYuan2024}
for more results on the related charged particle dynamics.

From the perspective of Lagrangian mechanics, Landau \& Lifshitz \cite{LandauLifshitz1971} formulate a more physically realistic model of motion.
Let
$x^{i}=(ct,\mathbf r)$ denote the spacetime position of a charged
particle and $(\phi,\mathbf A)$ the electromagnetic four-potential $(i=0,1,2,3)$.  For initial and final times $t_1$ and $t_2$, the action function for particle with mass $m$ and charge $e$, is given by
$$S=\int_{t_1}^{t_2}\Big[-mc^2\sqrt{1-\frac{v^2}{c^2}}+\frac{e}{c}\mathbf A\cdot\mathbf v-e\phi\Big]dt,$$
where $\mathbf r$ is the position vector, $c$ is the speed of light,
$\phi$ is the scalar potential, and $\mathbf A$ is the vector
potential depending on both position and time, the velocity
$\mathbf v=d\mathbf r/dt$ with magnitude $v=|\mathbf v|$. 
%The integrand defines the Lagrangian 
%$$L=-mc^2\sqrt{1-\frac{v^2}{c^2}}+\frac{e}{c}\mathbf A\cdot\mathbf v-e\phi.$$
In dimensionless case, the Euler--Lagrange equation for the particle is
%$$\frac{d}{dt}\left(\frac{\partial L}{\partial\mathbf v}\right)=\frac{\partial L}{\partial\mathbf r}.$$
%Substitution into the Euler--Lagrange equation therefore yields
$$\frac{d\mathbf p}{dt}=-\frac{e}{c}\frac{\partial\mathbf A}{\partial t}-e\nabla\phi+\frac{e}{c}\mathbf v\times(\nabla\times\mathbf A),$$
where $\mathbf p=m\mathbf v/\sqrt{1-v^2/c^2}$ is the momentum.
Defining the electric field $\mathbf E$ and magnetic field $\mathbf H$ by
$$\mathbf E=-\nabla\phi-\frac{1}{c}\frac{\partial\mathbf A}{\partial t},\qquad\mathbf H=\nabla\times\mathbf A,$$
then the equation can be rewritten as follows
$$\frac{d\mathbf p}{dt}=e\mathbf E+\frac{e}{c}\mathbf v\times\mathbf H.$$
In the nonrelativistic limit $v\ll c$, the momentum is approximately equal to
$m\mathbf v$. Thus a dimensionless parameter model is 
$$
\frac{d\mathbf v}{dt}=\mathbf E+\mathbf v\times\mathbf H.
$$

In the present paper, we study the dynamic behavior of charged particles in magnetic fields with non-Lipschitz force. Let $x,v\in\R^d$ be the position and the velocity, respectively. We consider the parameter-dependent nonautonomous ordinary differential equation
\begin{equation}\label{eq:general-system}
	\frac{\dd x}{\dd t}=v,\qquad\frac{\dd v}{\dd t}=\mathsf B(x,t)v-\lambda(t)v+|v|^\alpha f\Big(t,\frac{v}{|v|},\varepsilon\Big),
\end{equation}
with $0<\alpha<1$.
The parameter $\varepsilon\in\mathcal U_\varepsilon\subset\R^\ell$ is frozen along each trajectory. The magnetic field is represented by a real skew-symmetric matrix $\mathsf{B}(x,t)$. The function $\lambda(t)$ is a nonnegative damping coefficient. The function $f:\mathbb{R}\times\mathbb{S}^{d-1}\times\mathbb{R}^\ell\to\mathbb{R}^d$ is at least $C^1$. The loss of Lipschitz
regularity comes from the last term in \eqref{eq:general-system} at $v=0$ and leads to the nonuniqueness at this point. The equation describes the motion of charged particles in a special environment under the influence of a magnetic field and a damping term. 
In addition to magnetic and damping effects, we incorporate a physically motivated force that loses Lipschitz regularity at zero velocity \cite{FerreiraAraujo2024}. This loss of regularity may lead to nonunique continuations when a particle reaches rest, connecting our analysis to the theory of singular ordinary differential equations. Relevant developments include existence, uniqueness, and convergence results for regularized flows in a generalized-function framework \cite{KunzingerOberguggenbergerSteinbauerVickers2004}, continuation by local smoothing near isolated non-Lipschitz singularities \cite{DrivasMailybaev2021}, and statistical descriptions involving chaotic attractors and physical measures \cite{DrivasMailybaevRaibekas2024}. The target of this paper is to select a unique continuation by means of regularization and to investigate the complex dynamics in the particle motion.  

\subsection{Main results}

We fix the parameter and further assume that the directional
force is independent of time. Consider
\begin{equation}\label{eq:intro-regularization-system}
	\frac{\dd x}{\dd t}=v,\qquad\frac{\dd v}{\dd t}=\mathsf B(x,t)v-\lambda(t)v+|v|^\alpha f\left(\frac{v}{|v|}\right).
\end{equation}
Write
$$f_r(u)=f(u)\cdot u,\qquad f_s(u)=f(u)-f_r(u)u.$$
Set
$$q=w^{1-\alpha},\qquad\beta=\frac{2-\alpha}{1-\alpha}>2,\qquad\frac{\dd s}{\dd t}=\frac1q.$$
Then
\begin{align}
	\frac{\dd u}{\dd s}&=f_s(u)+q\mathsf B(x,t)u,\label{eq:intro-ren-u}\\
	\frac{\dd q}{\dd s}&=(1-\alpha)q\bigl[f_r(u)-\lambda(t)q\bigr],\label{eq:intro-ren-q}\\
	\frac{\dd t}{\dd s}&=q,\\
	\frac{\dd x}{\dd s}&=q^\beta u.\label{eq:intro-ren-tx}
\end{align}
Differentiating $q=w^{1-\alpha}$ in time $t$ gives
\begin{equation}\label{eq:q-physical}
	\frac{\dd q}{\dd t}=(1-\alpha)\bigl[f_r(u)-\lambda(t)q\bigr].
\end{equation}

We impose the following assumptions.

\begin{enumerate}
	
	\item[\textup{(R1)}]
	The force satisfies
	$f\in C^2(\mathbb S^{d-1};\R^d)$.
	There exist two locally exponentially stable equilibria $u_-$ and $u_+$ of
	$$\frac{\dd u}{\dd s}=f_s(u).$$
	There are neighborhoods $N_-$ and $N_+$ and constants $a_\pm>0$ such that
	\begin{align*}
	f_r(u)&\leq-a_-<0,\qquad&u\in N_-,\\
	f_r(u)&\geq a_+>0,\qquad&u\in N_+.
\end{align*}
We call $u_-$ the incoming direction because the speed contracts toward zero
near $u_-$, and we call $u_+$ the outgoing direction because the speed
tends to expand away from zero near $u_+$.
	
	\item[\textup{(R2)}]
	The functions $\mathsf B$ and $\lambda$ are $C^1$. They and their first
	derivatives are globally uniformly bounded. 
	
	\item[\textup{(R3)}]
The direction of velocity enters a sufficiently
small stable neighborhood of $u_-$ with sufficiently small positive $q$.
\end{enumerate}

For each $\nu>0$, let $F_\nu\in C^1(\R^d;\R^d)$ be a regularization of
the singular force
$|v|^\alpha f(v/|v|)$ of a suitable class. Let $(x^\nu,v^\nu)$ solve
\begin{equation}\label{eq:intro-regularized-equation}
	\frac{\dd x^\nu}{\dd t}=v^\nu,\qquad\frac{\dd v^\nu}{\dd t}=\mathsf B(x^\nu,t)v^\nu-\lambda(t)v^\nu+F_\nu(v^\nu).
\end{equation}
All singular and regularized systems are considered with the common
initial data $(x(t_0),v(t_0))=(x_0,v_0)$, with $v_0\ne0$. 
The velocity is called blowup if there exists a finite time $t_b$ such that $\lim\limits_{t\uparrow t_b}v(t)=0$ for some initial data $(x(0),v(0))=(x_0,v_0)$.
\begin{theorem}
	\label{thm:intro-regularization-selection}
	Assume \textup{(R1)}--\textup{(R3)}. Then the solution $(x(t),v(t))$ of \eqref{eq:intro-regularization-system} 
	reaches $v=0$ at a finite time $t_b$. There is a point $x_b\in\R^d$ such that
	$$\lim_{t\uparrow t_b}x(t)=x_b,\qquad\lim_{t\uparrow t_b}v(t)=0,\qquad\lim_{t\uparrow t_b}\frac{v(t)}{|v(t)|}=u_-.$$
	
	There are $\delta>0$ and a unique continuous curve
	$(x^*,v^*)\in C^0([t_0,t_b+\delta];\R^{2d})$
	such that 
	$(x^\nu,v^\nu)\longrightarrow(x^*,v^*)$
	in	$C^0([t_0,t_b+\delta];\R^{2d})$ as $\nu\to0$.
	For $t_b<t\leq t_b+\delta$, it solves the original system \eqref{eq:intro-regularization-system},
	$v^*(t)\neq0$, and
	$$\lim_{t\downarrow t_b}\frac{v^*(t)}{|v^*(t)|}=u_+.$$
\end{theorem}

We notice that the inner field $H$ with a
reference orbit $V_*$ should be fixed. This orbit solves
$$\frac{\dd V_*}{\dd\tau}=\mathcal H(V_*),\qquad V_*(0)=u_-.$$
For each $\nu>0$, let $t_{\rm ent}^\nu$ be the time at which the regularized
velocity first enters the regularization ball $|v|\leq\nu$. Introduce the inner time $\tau$ and the rescaled velocity $V^\nu$ by
$$t^\nu(\tau)=t_{\rm ent}^\nu+\nu^{1-\alpha}\tau,\qquad V^\nu(\tau)=\frac{v^\nu(t^\nu(\tau))}{\nu}.$$
Direct estimates show that $V^\nu\to V_*$ as $\nu\to0$ uniformly on every fixed interval. This convergence implies that, for small $\nu$, the regularized trajectories exit into the stable neighborhood of $u_+$. The Arzel\`a--Ascoli theorem gives a subsequential limit. The one-dimensional strong unstable manifold makes this limit solution unique. Hence all
subsequential limits coincide, and the full regularized family converges to
the same continuation.

We remark that the method of \cite{DrivasMailybaev2021} for the homogeneous autonomous case does not apply directly here: general magnetic and damping terms break the scaling invariance of \eqref{eq:general-system}, while the magnetic term couples the velocity direction and magnitude through the direction-dependent term $f_r$. When $u_+$ is an equilibrium in both settings, both results yield a unique regularization limit; ours is characterized by the strong unstable manifold of the extended system and generally lacks an explicit self-similar form. Moreover, adding a magnetic field may change the structure, stability, or basin of attraction of an attractor of $du/ds=f_s(s,u,\varepsilon)$, or destroy it altogether.

%Our second result focuses on the chaos of \ref{eq:general-system}. Chaos plays a significant role in particle transport and heating. Numerical studies of relativistic electrons interacting with whistler waves suggest that chaotic motion enhances phase-space diffusion and may contribute to pitch-angle scattering in the torus and auroral precipitation \cite{WykesChapmanRowlands2001}. During collisionless magnetic reconnection, in-plane Hall electric fields can likewise energize ions and induce chaotic trajectories, resulting in ion heating \cite{YoonBellan2021}.
 %The magnetic term also plays an important role in the generation of chaos. After reducing the magnetic term to a constant matrix multiplied by a variable coefficient, we can design the constant matrix to provide part of the linear structure that supports the onset of chaotic dynamics.
  
Our second result concerns chaos in (\ref{eq:general-system}), motivated by its role in particle transport and heating through enhanced phase-space diffusion \cite{WykesChapmanRowlands2001} and ion energization during magnetic reconnection \cite{YoonBellan2021}. Writing the magnetic term as a variable coefficient times a constant matrix, we choose the matrix to supply part of the linear structure supporting chaotic dynamics.  
  
We now study the velocity directional dynamics when the speed remains
positive. For
$$w=|v|>0,\qquad u=\frac{v}{|v|},$$
the direction satisfies
\begin{equation}\label{eq:intro-II-direction}
	\frac{\dd u}{\dd t}=\mathsf B(x,t)u+w^{\alpha-1}f_s(t,u,\varepsilon).
\end{equation}

We impose the following assumptions.

\begin{enumerate}
	
	\item[\textup{(H1)}]
The magnetic field has the separated form
$$\mathsf B(x,t)=k(t)\mathsf B_\star,$$
where $\mathsf B_\star\in\R^{d\times d}$ is a constant real skew-symmetric
matrix and
\begin{equation}\label{c1}
	k\in C^1(\R,(0,\infty)),\qquad0<\underline k\leq k(t)\leq\overline k<\infty.
\end{equation}
Moreover, there exists a common equilibrium direction $u_\star$ of \eqref{eq:intro-II-direction}:
$$\mathsf B(x,t)u_\star=0,\qquad f_s(t,u_\star,\varepsilon)=0$$
for all relevant $(x,t,\varepsilon)$. 

\item[\textup{(H2)}]  $f_r$ is independent of the spherical direction.
 Hence $w(t)$ is a known scalar function once the initial value $w_0>0$ and the parameter
are fixed. A suitable design of $\lambda(t)$ and $f_r(t,\varepsilon)$ yields solutions $w(t)$ with positive upper and lower bounds on $\mathbb{R}$. 
\end{enumerate}
For example, take $\lambda(t)=1/(1+t^2)$ and $f_r(t,u,\varepsilon)=t/(1+t^2)^2$. For initial value $w(0)=w_0>0$, $\dd w/\dd t=-\lambda(t)w+w^\alpha f_r(t,u,\varepsilon)$ has a solution 
$$w(t)=\Bigg\{e^{-(1-\alpha)\arctan t}\Bigg[w_0^{\,1-\alpha}+(1-\alpha)\int_0^te^{(1-\alpha)\arctan s}\frac{s}{(1+s^2)^2}\,\mathrm{d}s\Bigg]\Bigg\}^{\frac{1}{1-\alpha}},$$
which is positive and away from zero.

Let
$$\pi:\mathbb S^{d-1}\setminus\{-u_\star\}\longrightarrow\R^n,\qquad n=d-1,$$
be the stereographic projection centered at $u_\star$, satisfying
$\pi(u_\star)=0$.
Write
$X=\pi(u)$
and define the projected tangential force by
$$\mathcal N(t,X,\varepsilon)=D_u\pi(u)\,f_s(t,u,\varepsilon),\qquad u=\pi^{-1}(X).$$
Let $B\in\R^{n\times n}$ denote the matrix induced by the restriction of
$\mathsf B_\star$ to
$T_{u_\star}\mathbb S^{d-1}$. Then the velocity directional equation takes the local
form
\begin{equation}\label{eq:project}
	\frac{\dd X}{\dd t}=k(t)BX+w^{\alpha-1}(t,\varepsilon)\mathcal N(t,X,\varepsilon).
\end{equation}

\begin{enumerate}
	
	\item[\textup{(H3)}]
	Fix a finite integer $\rho\geq2$. On a fixed neighborhood of
	$\R\times\{(0,0)\}$,
	the projected tangential force
	$\mathcal N(t,X,\varepsilon)$
	is $C^1$ in $t$ and $C^{\rho+1}$ in $(X,\varepsilon)$. All state and parameter derivatives used
	up to $\rho+1$ orders are uniformly bounded for $t\in\R$. The common directional equilibrium point persists for every parameter.
	That is, in stereographic coordinates,
	\begin{equation}\label{persistent-projected-equilibrium}
	\mathcal N(t,0,\varepsilon)=0,\qquad(t,\varepsilon)\in\mathbb R\times\mathcal U_\varepsilon.
\end{equation}
	\item[\textup{(H4)}]
There exists a constant matrix
$A\in\R^{n\times n}$
such that for a preselected initial value $w(0)=w_0>0$
\begin{equation}\label{linear-force-factorization}
	D_X\mathcal N(t,0,0)=k(t)w^{1-\alpha}(t)A,\qquad t\in\mathbb R.
\end{equation}
\item[\textup{(H5)}]
Let
$$F_A=\frac{A+A^T}{2}$$
be the symmetric part of $A$. Let
$$\lambda_1\geq\lambda_2\geq\cdots\geq\lambda_n,\qquad n=d-1,$$
be the eigenvalues of $F_A$. Assume that
$\lambda_1\geq0$ and $\operatorname{tr}F_A<0$. Define
\begin{equation}\label{eq:intro-II-center-dimension}
	m=\max\left\{j\in\{1,\ldots,n\}:\sum_{i=1}^j\lambda_i\geq0\right\}.
\end{equation}
Then
$1\leq m<n$.
\end{enumerate}

We begin by defining the center manifold reduction as follows, which is adapted from \cite{NijholtPereiraQueirozTuraev2023}.

\textbf{Some Definitions.}
	Let
	$$H:\mathbb R\times\mathbb R^n\times\mathcal U_\varepsilon\to\mathbb R^n$$
	be a family of time-dependent vector fields, parameterized by
	$\varepsilon$ in an open neighborhood
	$\mathcal U_\varepsilon\subset\mathbb R^\ell$ of the origin. Assume that
	$H(t,0,\varepsilon)=0$ and $D_xH(t,0,0)=k(t)L$, where
	$k:\mathbb R\to(0,\infty)$ is continuous and the constant matrix $L$
	admits a fixed invariant splitting
	$\mathbb R^n=E^c\oplus E^s$.
	Assume that the variational equation
	$\dd x/\dd t=k(t)Lx$
	admits an exponential trichotomy on $\mathbb R$ with trivial unstable
	space.
	
	We join the parameter as
	a frozen variable satisfying $\dd\varepsilon/\dd t=0$. Thus the center
	space of the parameter-extended system is
	$E^c\times\mathbb R^\ell$. A local parameterized center manifold is a
	locally invariant graph of
	$x_s=h(t,x_c,\varepsilon)$,
	where $x_c\in E^c$, $x_s\in E^s$,
	$h(t,0,\varepsilon)=0$, and $D_{x_c}h(t,0,0)=0$.
	
	Let $P^c$ denote the projection onto $E^c$. The dynamics on the center
	manifold is governed by
	$$\frac{\dd x_c}{\dd t}=R(t,x_c,\varepsilon):=P^cH\bigl(t,x_c+h(t,x_c,\varepsilon),\varepsilon\bigr),\qquad\frac{\dd\varepsilon}{\dd t}=0.$$
	We call
	$R:\mathbb R\times E^c\times\mathcal U_\varepsilon\to E^c$
	the parameterized reduced vector field. It satisfies
	$R(t,0,\varepsilon)=0$
	and
	$D_{x_c}R(t,0,0)=k(t)L|_{E^c}$.
	The dynamical dimension of the parameterized center manifold is
	$\dim E^c$, while $\ell$ is the number of frozen parameters. $R(t,x_c,\varepsilon)$ is the $E^c$-valued
	vector field on each fixed-$\varepsilon$ fiber, while
	$\widetilde R(t,x_c,\varepsilon)=(R(t,x_c,\varepsilon),0)$ denotes
	its frozen-parameter extension.

	In the above definitions, the constant and linear parts of the parameterized
	reduced vector field $R(t,x_c,\varepsilon)$ are specified.  Fix an integer
	Taylor order $\rho\ge2$ and assume that $H(t,x,\varepsilon)$ is
	$C^{\rho+1}$ in $(x,\varepsilon)$.  For $2\le r\le\rho$, let
	$H_r(t,x,\varepsilon)$ denote the homogeneous term of degree $r$ in the Taylor
	expansion of $H(t,x,\varepsilon)$ at $(x,\varepsilon)=(0,0)$.  Thus
	$$H(t,x,\varepsilon)=H(t,0,0)+D_xH(t,0,0)x+D_\varepsilon H(t,0,0)\varepsilon+\sum_{r=2}^{\rho}H_r(t,x,\varepsilon)+\mathcal{O}(\|(x,\varepsilon)\|^{\rho+1}),$$
	where each $H_r(t,x,\varepsilon)$ is a homogeneous polynomial of degree $r$ in $(x,\varepsilon)$, whose coefficients are functions of $t$, thanks to the continuity of $H$ in $t$ and Lemma \ref{lem2}. For example, when $H$ is a mapping from $\mathbb{R}\times\mathbb{R}\times\mathcal U_\varepsilon$ to $\mathbb{R}$ with $\mathcal U_\varepsilon\subseteq\mathbb{R}$, the quadratic and cubic parts take the forms
	$$H_2(t,x,\varepsilon)=a_1(t)x^2+a_2(t)x\varepsilon+a_3(t)\varepsilon^2,$$
	and
	$$H_3(t,x,\varepsilon)=a_4(t)x^3+a_5(t)x^2\varepsilon+a_6(t)x\varepsilon^2+a_7(t)\varepsilon^3,$$
	for certain coefficient functions $a_i(t)$. It is worth noting that if $H(t,0,\varepsilon)=0$ for all $t,\varepsilon$, then $D_\varepsilon H(t,0,0)=0$ and the terms $\varepsilon^2$ and $\varepsilon^3$ are absent from $H_2$ and $H_3$ respectively.

\begin{theorem}
	\label{thm:structural-reduction}
	Assume \textup{(H1)}--\textup{(H5)}.
	Then there exists a constant skew-symmetric matrix
	$B\in\R^{n\times n}$
	induced by the magnetic field, such that the projected directional
	system \eqref{eq:project} has a local parameterized center manifold of dynamical dimension $m$.
	Its graph has uniformly bounded derivatives up to order
	$\rho$.
	
	The corresponding reduced vector field on the parameterized center manifold satisfies
	$$R(t,0,\varepsilon)=0,\qquad D_{x_c}R(t,0,0)=k(t)\Omega,$$
	where
	$\Omega:E^c\longrightarrow E^c$
	is nilpotent. Moreover, $B$ can be designed such that $\Omega$ is skew-symmetric if \textup{(H5)} is replaced by the condition that $F_A$ has zero eigenvalues with multiplicity $m$ and $\operatorname{tr}F_A<0$.
	
	Moreover, the Taylor expansion of degrees $2$ through $\rho$ of
	$\mathcal N$ can be chosen so that, for every $2\leq r\leq\rho$, the
	degree-$r$ Taylor expansion $R_r(t,x_c,\varepsilon)$ of the reduced vector field on the parameterized center manifold can be any $E^c$-valued homogeneous polynomial $\widehat a_r(t,x_c,\varepsilon)$ of degree $r$
	in $(x_c,\varepsilon)$ satisfying
	$\widehat a_r(t,0,\varepsilon)=0$.
	Its coefficient functions are continuously differentiable, and both the
	coefficient functions and their time derivatives are bounded on $\R$.
\end{theorem}

We impose the stronger factorization
\begin{equation}\label{eq:intro-II-original-factorization}
	f_s(t,u,\varepsilon)=k(t)w^{1-\alpha}(t,\varepsilon)g_s(u,\varepsilon),
\end{equation}
where
$g_s(u_\star,\varepsilon)=0$.
Then the directional equation has the form
$$\frac{\dd u}{\dd t}=k(t)\bigl(Bu+g_s(u,\varepsilon)\bigr),$$
where $B$ is extended by zero in the direction $u_\star$. The positive time
change
$$s(t)=\int_0^tk(\tau)\,\dd\tau$$
removes the scalar time dependence and gives an autonomous directional
equation on each fixed-parameter fiber.

On a three-dimensional parameterized center manifold, a rank-two nilpotent center block admits a generic three-parameter unfolding of a codimension-three nilpotent singularity, while a nonzero skew-symmetric center block yields the Hopf-zero case. Suitable nonlinearities produce Shilnikov homoclinic loops and hence chaos in the velocity direction system. Compared with the result of \cite{NijholtPereiraQueirozTuraev2023}, our result not only allows the Taylor coefficients of the reduced dynamics to be prescribed, subject to the persistence of the equilibrium, but also covers nonzero skew-symmetric center blocks and the resulting chaotic dynamics.

Section~2 studies deterministic selection through the non-Lipschitz point, covering regularized dynamics, subsequential limits, uniqueness of continuation, and a planar example.
Section~3 develops the stereographic formulation, spectral splitting, exponential trichotomy, and parameterized center-manifold reduction for the velocity direction dynamics, then establishes chaos in the nilpotent and skew-symmetric center cases.

\section{Deterministic selection through the non-Lipschitz point}

\subsection{Hypotheses and renormalized variables}
\label{sec:regime-regularization}
Throughout this section, we assume that \textup{(R1)}--\textup{(R3)} hold.
We state the assumptions in detail and introduce the associated notation.
Throughout this part we assume $f\in C^2(\Sph;\R^d)$.  Consequently, $f_s$,
$f_r$, and the local first- and second-order Taylor expansions used below are
bounded.

 The two points
$u_\pm\in\Sph$ are locally exponentially stable equilibrium points of
\begin{equation}\label{re1}
	\frac{\dd u}{\dd s}=f_s(u).
\end{equation}
More precisely, $f_s(u_\pm)=0$, and there exist
$\delta_\pm>0$, $M_\pm\ge1$, and $\omega_\pm>0$ such that whenever $\dist(u_0,u_\pm)<\delta_\pm$, the solution $U(s;u_0)$ of \eqref{re1} starting from $u_0$ exists for
every $s\ge0$ and satisfies
\begin{equation}\label{eq:pm-exp-stable}
	\dist(U(s;u_0),u_\pm)\le M_\pm\ee^{-\omega_\pm s}\dist(u_0,u_\pm),\qquad s\ge0.
\end{equation}
By this local exponential stability and $f_s\in C^1$, the converse Lyapunov
theorem for exponential stability
\cite[Theorem~4.5]{Khalil1992} permits us to take neighborhoods
$N_\pm\subset\Sph$ of $u_\pm$ and construct functions
$W_\pm\in C^1(N_\pm)$ and positive constants
$\mu_\pm,c_{1,\pm},c_{2,\pm},C_{W,\pm}$ such that
\begin{align}
	c_{1,\pm}\dist(u,u_\pm)^2&\le W_\pm(u)\le c_{2,\pm}\dist(u,u_\pm)^2,\label{eq:pm-W1}\\
	\nabla W_\pm(u)\cdot f_s(u)&\le-2\mu_\pm W_\pm(u),\label{eq:pm-W2}\\
	|\nabla W_\pm(u)|&\le C_{W,\pm}\sqrt{W_\pm(u)}.\label{eq:pm-W3}
\end{align}
The functions $W_\pm$
and these three estimates follow from local exponential stability.

Additionally, we assume that there are $a_\pm>0$ such that,
\begin{equation}\label{eq:pm-rad}
	\pm f_r(u)\ge a_\pm>0,\qquad u\in N_\pm.
\end{equation}
Continuity then gives the
displayed uniform inequalities, after reducing $N_\pm$, with
$a_+=f_r(u_+)/2$ and $a_-=-f_r(u_-)/2$.

After reducing $N_+$, choose $\rho_+>0$ such that the set
\begin{equation}\label{eq:plus-stable-block}
	\mathcal K_+:=\{u\in N_+:W_+(u)\le\rho_+^2\}
\end{equation}
is compact and contained in the interior of $N_+$.
It is merely a sufficiently small sublevel set of the constructed function
$W_+$ and is not a new dynamical assumption.

The exponential estimate
and the resulting Lyapunov estimates control convergence of the direction to
$u_\pm$.  The additional radial condition controls convergence of the speed
to zero at $u_-$ and its outward growth from zero at $u_+$.

The magnetic field $\mathsf B(x,t)$ is $C^1$ in $(x,t)$, and the damping
coefficient $\lambda(t)$ is $C^1$ in $t$.  Both functions and all their first
derivatives are bounded on their full domains.  Hence there are constants $M_B,\Lambda<\infty$ such that
\begin{align}
	\|\mathsf B(x,t)\|&\le M_B,\label{eq:B-bound}\\
	0\le\lambda(t)&\le\Lambda.\label{eq:lambda-bound}
\end{align}

Fix a solution of
\eqref{eq:intro-regularization-system} with initial data
$(x(t_0),v(t_0))=(x_0,v_0)$, where $v_0\neq0$. On the maximal forward interval on which $v(t)\neq0$, define the
renormalized time by
$$s(t)=\int_{t_0}^{t}\frac{\dd\tau}{q(\tau)},\qquad q(t)=|v(t)|^{1-\alpha}.$$
Thus
$s(t_0)=0$.

\subsection{Incoming dynamics and finite-time blowup}
In this section, we give a description that the particle velocity reaches zero in finite time. Following \cite{DrivasMailybaev2021}, this scenario is called finite-time blow-up. There is an $r_->0$ such that
$$K_-:=\{u\in N_-:W_-(u)\le r_-^2\}$$
is compact and is contained in the interior of $N_-$.  Choose a compact
entrance set $N_-^{\mathrm{in}}$ that is contained in the interior of $K_-$.
\begin{theorem}
\label{thm:pre}
Assume that there is a renormalized time $s_0$ such that $u(s_0)\in N_-^{\mathrm{in}}$ and $q(s_0)>0$ is sufficiently small, that is,
\begin{equation}\label{eq:boundary-smallness}
	C_0q(s_0)<\mu_-r_-,\qquad C_0:=\frac{C_{W,-}M_B}{2}.
\end{equation}
Then the renormalized solution exists for every $s\ge s_0$ and
$$u(s)\in K_-,\qquad0<q(s)\le q(s_0)\ee^{-c_-(s-s_0)},\qquad c_-=(1-\alpha)a_->0.$$
Moreover,
$$\lim_{s\to\infty}u(s)=u_-,\qquad\lim_{s\to\infty}q(s)=0.$$
\end{theorem}

\begin{proof}
Define the first time at which the direction leaves the interior of $K_-$ by
$$s_*:=\inf\{s>s_0:u(s)\notin\operatorname{int}K_-\},$$
with $s_*=\infty$ if no such time exists.  Since
$u(s_0)\in\operatorname{int}K_-$ and $u$ is continuous, $s_*>s_0$.

Fix $s<s_*$.  For every $\sigma\in[s_0,s]$,
$$u(\sigma)\in K_-\subset N_-.$$
Consequently, on $[s_0,s]$, equations \eqref{eq:intro-ren-q} and
\eqref{eq:pm-rad}, together with $\lambda\ge0$, give
$$\frac{\dd q}{\dd s}=(1-\alpha)q[f_r(u)-\lambda q]\le-(1-\alpha)a_-q.$$
Divide by $q>0$ and integrate from $s_0$ to $s$:
$$\log q(s)-\log q(s_0)\le-(1-\alpha)a_-(s-s_0).$$
Exponentiation yields
\begin{equation}\label{eq:q-decay}
	q(s)\le q(s_0)\ee^{-c_-(s-s_0)},\qquad c_-=(1-\alpha)a_->0.
\end{equation}
At this stage, \eqref{eq:q-decay} has been proved only for
$s_0\le s<s_*$.

Let $R(s)=\sqrt{W_-(u(s))}$ and $R(s_0)=\sqrt{W_-(u(s_0))}<r_-$. By \eqref{eq:intro-ren-u},
\begin{align*}
	\frac{\dd W_-(u(s))}{\dd s}&=\nabla W_-(u)\cdot f_s(u)+q\nabla W_-(u)\cdot\mathsf B(x,t)u\\
	&\le-2\mu_-W_-(u)+q|\nabla W_-(u)|\,\|\mathsf B(x,t)\|\\
	&\le-2\mu_-R^2+C_{W,-}M_BqR.
\end{align*}
Then
\begin{equation}\label{eq:R-diff}
	\frac{\dd R}{\dd s}\le-\mu_-R+C_0q,\qquad C_0=\frac{C_{W,-}M_B}{2}.
\end{equation}
For any $s<s_*$, \eqref{eq:q-decay} gives $q(s)\le q(s_0)$.  Applying an
integrating factor to \eqref{eq:R-diff} and integrating from $s_0$ to $s$
gives
$$R(s)\le\ee^{-\mu_-(s-s_0)}R(s_0)+\frac{C_0q(s_0)}{\mu_-}\bigl(1-\ee^{-\mu_-(s-s_0)}\bigr).$$
The right-hand side is a convex combination of $R(s_0)$ and
$C_0q(s_0)/\mu_-$.  By $R(s_0)<r_-$ and
\eqref{eq:boundary-smallness}, $R(s)<r_-$. In fact, define the common fixed margin from the boundary
by
\begin{equation}\label{eq:R-fixed-margin}
	\varepsilon_-:=\min\left\{r_--R(s_0),\quad r_--\frac{C_0q(s_0)}{\mu_-}\right\}>0.
\end{equation}
This number depends only on the entrance data and not on $s$.  Put
$$\theta(s):=\ee^{-\mu_-(s-s_0)}\in(0,1],\qquad s\ge s_0.$$
Then $R(s_0)\le r_--\varepsilon_-$ and
$C_0q(s_0)/\mu_-\le r_--\varepsilon_-$.  Thus, for every
$s_0\le s<s_*$,
\begin{align*}
	R(s)&\le\theta(s)R(s_0)+(1-\theta(s))\frac{C_0q(s_0)}{\mu_-}\\
	&\le\theta(s)(r_--\varepsilon_-)+(1-\theta(s))(r_--\varepsilon_-)\\
	&=r_--\varepsilon_-<r_-.
\end{align*}

Suppose that $s_*<\infty$.  By the definition of the first exit time, the
continuity of $u$, and the fact that the compact set $K_-$ is contained in the
interior of $N_-$, one must have $u(s_*)\in\partial K_-$ and hence
$R(s_*)=r_-$.
But $R$ is continuous and the number $\varepsilon_->0$ in
\eqref{eq:R-fixed-margin} is independent of $s$, so letting $s\uparrow s_*$
gives
$$R(s_*)=\lim_{s\uparrow s_*}R(s)\le r_--\varepsilon_-<r_-,$$
which is a contradiction.  Therefore $s_*=\infty$.  This proves for $s\ge s_0$, $u(s)\in K_-\subset N_-$
and the validity of \eqref{eq:q-decay} for every $s\ge s_0$.

Apply the integrating factor $\ee^{\mu_-s}$ to \eqref{eq:R-diff}:
$$\frac{\dd\bigl(\ee^{\mu_-s}R(s)\bigr)}{\dd s}\le C_0\ee^{\mu_-s}q(s).$$
Integration from $s_0$ to $s$ gives
\begin{equation}\label{eq:R-integral}
	R(s)\le\ee^{-\mu_-(s-s_0)}R(s_0)+C_0\int_{s_0}^{s}\ee^{-\mu_-(s-\sigma)}q(\sigma)\,\dd\sigma.
\end{equation}
Substitute \eqref{eq:q-decay} into \eqref{eq:R-integral}.  If $\mu_-\ne c_-$, then
$$\int_{s_0}^{s}\ee^{-\mu_-(s-\sigma)}\ee^{-c_-(\sigma-s_0)}\,\dd\sigma=\frac{\ee^{-c_-(s-s_0)}-\ee^{-\mu_-(s-s_0)}}{\mu_--c_-}.$$
If $\mu_-=c_-$, the integral equals
$(s-s_0)\ee^{-\mu_-(s-s_0)}$.
Both expressions tend to zero as $s\to\infty$.  Hence $R(s)\to0$ and, by
\eqref{eq:pm-W1}, $\dist(u(s),u_-)\to0$ as $s\to\infty$.
\end{proof}

\begin{proposition}\label{prop:pre-reconstruction}
Under the hypotheses and conclusions of Theorem~\ref{thm:pre}, there are a
finite $t_b$ and a point $x_b\in\R^d$ such that
$$\lim_{s\to\infty}t(s)=t_b<\infty,\qquad\lim_{s\to\infty}x(s)=x_b,\qquad\lim_{s\to\infty}v(s)=0,\qquad s_b:=\lim_{t\uparrow t_b}s(t)=\infty.$$
\end{proposition}

\begin{proof}
Set
$z(s):=\log w(s)=\frac{1}{1-\alpha}\log q(s)$.
As long as $q(s)>0$, equation \eqref{eq:intro-ren-q} gives
\begin{equation}\label{eq:pre-log-w}
	\frac{\dd z}{\dd s}=\frac{1}{1-\alpha}\frac1q\frac{\dd q}{\dd s}=f_r(u(s))-\lambda(t(s))q(s).
\end{equation}
Theorem~\ref{thm:pre} gives $q(s)\to0$ as $s\to\infty$, and therefore
\begin{equation}\label{eq:pre-z-minus-infinity}
	\lim_{s\to\infty}w(s)=0,\qquad\lim_{s\to\infty}z(s)=-\infty.
\end{equation}
Theorem~\ref{thm:pre} also gives $u(s)\in K_-$ and
$0<q(s)\le q(s_0)$.  Since $K_-$ is compact, $f_r$ is bounded there.  Hence
\begin{equation}\label{eq:pre-log-w-bound}
	\left|\frac{\dd z}{\dd s}\right|\le\max_{u\in K_-}|f_r(u)|+\Lambda q(s_0)=:M_z<\infty.
\end{equation}
Integrating over any finite interval from $s_0$ gives
$$z(s)\ge z(s_0)-M_z(s-s_0)>-\infty.$$
Thus $\log w$ cannot tend to $-\infty$ at a finite $s$, and $w$ cannot reach
zero at a finite $s$.  Together with \eqref{eq:pre-z-minus-infinity}, this
shows that the renormalized time corresponding to the vanishing speed is
$s_b=\infty$.

Since $\dd t/\dd s=q(s)$, for every finite $s\ge s_0$,
\begin{equation}\label{eq:pre-t-inverse}
	t(s)=t(s_0)+\int_{s_0}^{s}q(\sigma)\,\dd\sigma.
\end{equation}
The improper integral converges by \eqref{eq:q-decay}, so
\begin{equation}\label{eq:pre-tb-integral}
	t_b=t(s_0)+\int_{s_0}^{\infty}q(\sigma)\,\dd\sigma.
\end{equation}
Moreover,
\begin{equation}\label{eq:time-tail}
	t_b-t(s_0)=\int_{s_0}^{\infty}q(s)\,\dd s\le q(s_0)\int_0^\infty\ee^{-c_-\sigma}\,\dd\sigma=\frac{q(s_0)}{c_-}<\infty,
\end{equation}
and, for every finite $s\ge s_0$,
\begin{equation}\label{eq:time-tail-from-s}
	0<t_b-t(s)=\int_s^\infty q(\sigma)\,\dd\sigma\le\frac{q(s_0)}{c_-}\ee^{-c_-(s-s_0)}.
\end{equation}
Consequently,
$$\lim_{s\to\infty}\bigl(t_b-t(s)\bigr)=0.$$
Thus $t_b<\infty$, $t(s)<t_b$ for every finite $s$, and
$t(s)\to t_b$ as $s\to\infty$.

Since $\dd x/\dd s=q^\beta u$ and $|u|=1$,
\begin{equation}
	\int_{s_0}^{\infty}\left|\frac{\dd x}{\dd s}\right|\,\dd s=\int_{s_0}^{\infty}q(s)^\beta\,\dd s\le q(s_0)^\beta\int_0^\infty\ee^{-\beta c_-\sigma}\,\dd\sigma=\frac{q(s_0)^\beta}{\beta c_-}<\infty.
\end{equation}
Thus $x(s)$ satisfies the Cauchy
criterion as $s\to\infty$.  Completeness of $\R^d$ gives
$x(s)\to x_b$ as $s\to\infty$.  

Finally, $w=q^{1/(1-\alpha)}\to0$ and $v=wu\to0$ as $s\to\infty$.
\end{proof}

\subsection{The regularization condition}
\label{sec:reg-layer}

Following Section 4.1 and Equation (33) of
\cite{DrivasMailybaev2021}, we first define the patched vector field.  The
construction keeps the original homogeneous vector field unchanged
outside the ball of radius $\nu$ and replaces it inside the ball by a
scaled one.  Choose
$$H\in C^1(\overline{B_1};\R^d),\qquad\overline{B_1}:=\{V\in\R^d:|V|\le1\}$$
so that the patched force 
\begin{equation}\label{eq:Fnu}
	F_\nu(v)=\begin{cases}
	|v|^\alpha f(v/|v|),&|v|>\nu,\\
	\nu^\alpha H(v/\nu),&|v|\le\nu
\end{cases}
\end{equation}
 is continuously
 differentiable everywhere. For \eqref{eq:Fnu}, this means that the values and
 first derivatives of the inner and outer expressions must agree on
$|v|=\nu$.

  For each $\nu>0$, let
$(x^\nu,v^\nu)$ solve
\begin{equation}\label{eq:regularized-system}
	\frac{\dd x^\nu}{\dd t}=v^\nu,\qquad\frac{\dd v^\nu}{\dd t}=\mathsf B(x^\nu,t)v^\nu-\lambda(t)v^\nu+F_\nu(v^\nu),
\end{equation}
with initial data $(x^\nu(t_0),v^\nu(t_0))=(x_0,v_0)$. This initial value problem has a unique global solution for every
fixed $\nu>0$.

Define
\begin{equation}\label{eq:H-piece}
	\mathcal H(V)=\begin{cases}
	H(V),&|V|\le1,\\
	|V|^\alpha f(V/|V|),&|V|\ge1.
\end{cases}
\end{equation}
Then
$$F_\nu(\nu V)=\nu^\alpha\mathcal H(V),$$
and $\mathcal H$ is a $C^1$ vector field across
$|V|=1$. We now state the detailed assumptions on the function 
$H$ used in the regularization.

\begin{assumption}\label{ass:expel}
	
 Consider the equation
\begin{equation}\label{eq:inner-limit}
	\frac{\dd V}{\dd\tau}=\mathcal H(V),\qquad V(0)=u_-.
\end{equation}
Denote its solution by $V_*(\tau)$.  We impose the following requirements.
\begin{enumerate}
 \item The entrance is strictly inward: there is a $\kappa_->0$ such that
 \begin{equation}\label{eq:inner-entry-transverse}
	u_-\cdot H(u_-)=f_r(u_-)\le-2\kappa_-<0.
\end{equation}
 \item There is a finite $\tau_*>0$ such that
 \begin{equation}\label{eq:reference-first-exit}
	|V_*(\tau)|<1,\qquad0<\tau<\tau_*,\qquad|V_*(\tau_*)|=1.
\end{equation}
 \item Put $y_e=V_*(\tau_*)$.  There is a compact exit set $K_e\subset\mathbb S^{d-1}$ such that $y_e$ lies in the interior
 of $K_e$, and $K_e$ is contained in the interior of $\mathcal K_+$.
 Since $\mathcal{K}_+$ is a subset of $N_+$, there exists $\kappa_e$ such that
 \begin{equation}\label{eq:reference-exit-transverse}
	y_e\cdot\mathcal H(y_e)=y_e\cdot f(y_e)=f_r(y_e)\ge2\kappa_e>0.
\end{equation}
\end{enumerate}
\end{assumption}

\subsection{Comparison with the reference solution}

For $|v|>\nu$, equation \eqref{eq:Fnu} agrees with the singular force in
\eqref{eq:intro-regularization-system}.  Hence, for sufficiently small $\nu$, the unregularized and regularized trajectories from the same initial data coincide up
to their entrance into the regularization ball.  Let $t_{\rm ent}^\nu$ be the time of entry into the regularization region,  at
which
$$|v(t_{\rm ent}^\nu)|=|v^\nu(t_{\rm ent}^\nu)|=\nu.$$
Define the common entrance data by
$$x_{\rm ent}^\nu:=x(t_{\rm ent}^\nu)=x^\nu(t_{\rm ent}^\nu),\qquad v_{\rm ent}^\nu:=v(t_{\rm ent}^\nu)=v^\nu(t_{\rm ent}^\nu),\qquad u_{\rm ent}^\nu:=\frac{v_{\rm ent}^\nu}{\nu}.$$
Then $q(t_{\rm ent}^\nu)=\nu^{1-\alpha}$. By \eqref{eq:q-physical} and $u(t)\to u_-$ as $t\to t_b$, there are
$0<c<C$  such that, 
$$-C\le\frac{\dd q(t)}{\dd t}\le-c<0$$
near $t=t_b$.
Integrating from $t_{\rm ent}^\nu$ to $t_b$, and using $q(t_b)=0$, gives
\begin{equation}\label{eq:entry-time-scale}
	\frac{\nu^{1-\alpha}}{C}\le t_b-t_{\rm ent}^\nu\le\frac{\nu^{1-\alpha}}{c}.
\end{equation}
On this interval $|v|\le\nu$, and hence
\begin{equation}\label{eq:entry-x-scale}
	|x_b-x_{\rm ent}^\nu|\le\int_{t_{\rm ent}^\nu}^{t_b}|v(t)|\,\dd t\le\nu(t_b-t_{\rm ent}^\nu)\le\frac{\nu^{2-\alpha}}{c}.
\end{equation}
The corresponding renormalized entrance times tend to infinity as
$\nu\to0$, so
\begin{equation}\label{eq:entry-direction}
	\lim_{\nu\to0}u_{\rm ent}^\nu=u_-.
\end{equation}

Starting at the entrance time, introduce the scaling
\begin{equation}\label{eq:inner-scale}
	t^\nu(\tau):=t_{\rm ent}^\nu+\nu^{1-\alpha}\tau,\qquad V^\nu(\tau):=\frac{v^\nu(t^\nu(\tau))}{\nu},\qquad X^\nu(\tau):=\frac{x^\nu(t^\nu(\tau))-x_{\rm ent}^\nu}{\nu^{2-\alpha}}.
\end{equation}
Then
\begin{align}
	\frac{\dd X^\nu}{\dd\tau}&=V^\nu,\\
	\frac{\dd V^\nu}{\dd\tau}&=\mathcal H(V^\nu)+\nu^{1-\alpha}\bigl[\mathsf B(x^\nu(t^\nu(\tau)),t^\nu(\tau))-\lambda(t^\nu(\tau))\Id\bigr]V^\nu,\label{eq:inner-V}
\end{align}
with initial data $X^\nu(0)=0$ and
$V^\nu(0)=u_{\rm ent}^\nu$, where $\mathcal H$ is the unit-scale field obtained by matching the fixed
inner function $H$ to the outer field in \eqref{eq:H-piece}.

Recall that $V_*$ is the reference solution fixed in
Assumption~\ref{ass:expel}; equivalently, it is the unique solution of
$$\frac{\dd V_*}{\dd\tau}=\mathcal H(V_*),\qquad V_*(0)=u_-.$$
By \eqref{eq:reference-first-exit}, its exit data are denoted by
$$\tau_e:=\tau_*,\qquad y_e:=V_*(\tau_e),\qquad|y_e|=1.$$
For every $\nu>0$, define the extended first exit time
$$\tau_e^\nu:=\inf\{\tau>0:|V^\nu(\tau)|=1\},$$
and $\tau_e^\nu=\infty$ if $|V^\nu(\tau)|<1$ for all $\tau$. When $\tau_e^\nu<\infty$, define the corresponding exit point by
$y_e^\nu:=V^\nu(\tau_e^\nu)$.

\begin{lemma}\label{lem:inner}\label{lem:inner-exit}
If
$u_{\rm ent}^\nu\to u_-$ as $\nu\to0$, then there is a fixed
$T>\tau_*$, independent of $\nu$, such that
$$\lim_{\nu\to0}\sup_{0\le\tau\le T}|V^\nu(\tau)-V_*(\tau)|=0.$$
For all sufficiently small $\nu$, $\tau_e^\nu$ is finite and
\begin{align}
	\lim_{\nu\to0}\tau_e^\nu&=\tau_e=\tau_*,\label{eq:actual-exit-time-reference}\\
	\lim_{\nu\to0}y_e^\nu&=y_e,\label{eq:actual-exit-point-reference}\\
	y_e^\nu&\in K_e.
\end{align}
\end{lemma}

\begin{proof}
Set $\Theta(V):=|V|^2-1$.
Along this reference solution,
$$\frac{\dd\Theta(V_*(\tau))}{\dd\tau}=2V_*(\tau)\cdot\mathcal H(V_*(\tau)).$$
By \eqref{eq:inner-entry-transverse} and
\eqref{eq:reference-exit-transverse},
$$\left.\frac{\dd\Theta(V_*(\tau))}{\dd\tau}\right|_{\tau=0}\le-4\kappa_-,\qquad\left.\frac{\dd\Theta(V_*(\tau))}{\dd\tau}\right|_{\tau=\tau_e}\ge4\kappa_e.$$
Then there exist $h_->0$ and $h_+>0$, with
$h_-<\tau_e-h_+$, satisfying
\begin{align}
	\frac{\dd\Theta(V_*(\tau))}{\dd\tau}&\le-3\kappa_-,\qquad0\le\tau\le h_-,\label{eq:theta-entry-negative}\\
	\frac{\dd\Theta(V_*(\tau))}{\dd\tau}&\ge3\kappa_e,\qquad\tau_e-h_+\le\tau\le\tau_e+h_+.\label{eq:theta-reference-exit-neighborhood}
\end{align}
Set
$$T:=\tau_e+h_+,\qquad I_{\rm mid}:=[h_-,\tau_e-h_+].$$
By \eqref{eq:reference-first-exit}, the reference solution $V_*$ lies strictly
inside the unit ball on $I_{\rm mid}$, so
$$m:=-\max_{\tau\in I_{\rm mid}}\Theta(V_*(\tau))>0.$$
Integrating \eqref{eq:theta-entry-negative} and \eqref{eq:theta-reference-exit-neighborhood} gives
\begin{align}
	\Theta(V_*(\tau_e-h_+))&\le-3\kappa_eh_+,\label{eq:theta-reference-before-exit}\\
	\Theta(V_*(\tau_e+h_+))&\ge3\kappa_eh_+.\label{eq:theta-reference-after-exit}
\end{align}

Since $\mathcal H$ is globally Lipschitz on $\R^d$, there exists a global Lipschitz constant $L$ such that for $P,Q\in\R^d$,
$$|\mathcal H(P)-\mathcal H(Q)|\le L|P-Q|.$$
Finally, set
$$M_*:=\max_{0\le\tau\le T}|V_*(\tau)|<\infty,\qquad M_H:=\max_{0\le\tau\le T}|\mathcal H(V_*(\tau))|<\infty.$$

Write
\begin{equation}\label{eq:delta}
	\delta_\nu:=\nu^{1-\alpha}(M_B+\Lambda),\qquad\lim_{\nu\to0}\delta_\nu=0.
\end{equation}
Indeed, for any fixed $T>0$, global boundedness gives 
\begin{align*}
	\sup_{0\le\tau\le T}\nu^{1-\alpha}\bigl(\|\mathsf B(x^\nu(t^\nu(\tau)),t^\nu(\tau))\|+|\lambda(t^\nu(\tau))|\bigr)&\le\delta_\nu,\\
	\lim_{\nu\to0}\sup_{0\le\tau\le T}\nu^{1-\alpha}\bigl(\|\mathsf B(x^\nu(t^\nu(\tau)),t^\nu(\tau))\|+|\lambda(t^\nu(\tau))|\bigr)&=0,
\end{align*}
along the trajectory.
Set
$R_\nu(\tau):=\nu^{1-\alpha}\bigl[\mathsf B(x^\nu(t^\nu(\tau)),t^\nu(\tau))-\lambda(t^\nu(\tau))\Id\bigr]$.
Hence
\begin{equation}\label{eq:Rnu-bound}
	\sup_{0\le\tau\le T}\|R_\nu(\tau)\|\le\nu^{1-\alpha}(M_B+\Lambda)=\delta_\nu.
\end{equation}
Let $D(\tau):=V^\nu(\tau)-V_*(\tau)$.
Integrating equations \eqref{eq:inner-V} and \eqref{eq:inner-limit} from $0$ to any $0\le\tau\le T$ and subtracting
gives
\begin{align*}
	D(\tau)=u_{\rm ent}^\nu-u_-+\int_0^\tau\bigl[\mathcal H(V^\nu(\sigma))-\mathcal H(V_*(\sigma))\bigr]\,\dd\sigma+\int_0^\tau R_\nu(\sigma)V^\nu(\sigma)\,\dd\sigma.
\end{align*}
Use Lipschitz continuity of $\mathcal{H}$ and
$$|V^\nu(\sigma)|\le|V_*(\sigma)|+|D(\sigma)|\le M_*+|D(\sigma)|$$
to obtain
\begin{align*}
	|D(\tau)|&\le|u_{\rm ent}^\nu-u_-|+M_*T\delta_\nu+(L+\delta_\nu)\int_0^\tau|D(\sigma)|\,\dd\sigma.
\end{align*}
Set
$a_\nu:=|u_{\rm ent}^\nu-u_-|+M_*T\delta_\nu$.
By hypothesis $a_\nu\to0$ as $\nu\to0$.  For sufficiently small $\nu$,
$\delta_\nu\le1$, and hence
$$|D(\tau)|\le a_\nu+(L+1)\int_0^\tau|D(\sigma)|\,\dd\sigma.$$
Gronwall's inequality gives
\begin{align*}
	\sup_{0\le\tau\le T}|V^\nu(\tau)-V_*(\tau)|&\le a_\nu\ee^{(L+1)T},\\
	\lim_{\nu\to0}\sup_{0\le\tau\le T}|V^\nu(\tau)-V_*(\tau)|&=0.
\end{align*}
This proves the first limit. For sufficiently small $\nu$, $V_*(\tau)$ and $V^\nu(\tau)$ have uniform upper bound $M_0$.

Direct calculation gives
\begin{equation}\label{eq:theta-nu-derivative}
	\frac{\dd\Theta(V^\nu(\tau))}{\dd\tau}=2V^\nu(\tau)\cdot\mathcal H(V^\nu(\tau))+2V^\nu(\tau)\cdot R_\nu(\tau)V^\nu(\tau).
\end{equation}
The function $G(V):=2V\cdot\mathcal H(V)$ is Lipschitz on the closed ball
$|V|\le M_0$, and the absolute value of the second term on the right-hand side of \eqref{eq:theta-nu-derivative} is at most
$2M_0^2\delta_\nu$.  Hence there is a constant $C_0$, independent of $\nu$,
such that
\begin{align}
	\sup_{0\le\tau\le T}\left|\frac{\dd\Theta(V^\nu(\tau))}{\dd\tau}-\frac{\dd\Theta(V_*(\tau))}{\dd\tau}\right|&\le C_0(a_\nu\ee^{(L+1)T}+\delta_\nu),\label{eq:theta-derivative-difference}\\
	\lim_{\nu\to0}\sup_{0\le\tau\le T}\left|\frac{\dd\Theta(V^\nu(\tau))}{\dd\tau}-\frac{\dd\Theta(V_*(\tau))}{\dd\tau}\right|&=0.\nonumber
\end{align}
By \eqref{eq:theta-entry-negative} and
\eqref{eq:theta-derivative-difference}, for sufficiently small $\nu$,
$$\frac{\dd\Theta(V^\nu(\tau))}{\dd\tau}\le-2\kappa_-<0,\qquad0\le\tau\le h_-.$$
Since $|u_{\rm ent}^\nu|=1$, one has $\Theta(V^\nu(0))=0$.  Integrating from
$0$ to $0<\tau\le h_-$ gives
\begin{equation}\label{eq:actual-entry-inside}
	\Theta(V^\nu(\tau))\le-2\kappa_-\tau<0.
\end{equation}
Thus $V^\nu$ lies inside the unit ball strictly for $0<\tau\le h_-$.

Take $\nu$ sufficiently small such that $V^\nu(\tau)$ is close enough to $V_*(\tau)$. Hence,
$$|\Theta(V^\nu(\tau))-\Theta(V_*(\tau))|=|(V^\nu-V_*)\cdot(V^\nu+V_*)|<\min\left\{\frac m2,\,\kappa_eh_+\right\}.$$
Then
$$\Theta(V^\nu(\tau))\le-\frac m2<0,\qquad\tau\in I_{\rm mid}.$$
Also \eqref{eq:theta-reference-before-exit} and
\eqref{eq:theta-reference-after-exit} give
\begin{align*}
	\Theta(V^\nu(\tau_e-h_+))&\le-2\kappa_eh_+<0,\\
	\Theta(V^\nu(\tau_e+h_+))&\ge2\kappa_eh_+>0.
\end{align*}
The intermediate value theorem gives a
$\tau_e^\nu\in(\tau_e-h_+,\tau_e+h_+)$
such that $\Theta(V^\nu(\tau_e^\nu))=0$.

Again using \eqref{eq:theta-derivative-difference} and
\eqref{eq:theta-reference-exit-neighborhood}, for all sufficiently small
$\nu$,
$\dd\Theta(V^\nu)/\dd\tau\ge2\kappa_e$
throughout $[\tau_e-h_+,\tau_e+h_+]$.  Thus
$\Theta(V^\nu(\tau))$ is strictly increasing there and its zero is unique.
Together with \eqref{eq:actual-entry-inside} and the strict negative sign on
the middle interval, this shows that $\tau_e^\nu$ is the first exit after the
entrance. 

We now establish the remaining limit relations.
Both $\tau_e$ and $\tau_e^\nu$ lie in
$[\tau_e-h_+,\tau_e+h_+]$, where
$\dd\Theta(V^\nu)/\dd\tau\ge2\kappa_e$.  Hence
\begin{align*}
	2\kappa_e|\tau_e^\nu-\tau_e|&\le|\Theta(V^\nu(\tau_e))-\Theta(V^\nu(\tau_e^\nu))|\\
	&=|\Theta(V^\nu(\tau_e))-\Theta(V_*(\tau_e))|\\
	&\le2M_0|V^\nu(\tau_e)-V_*(\tau_e)|,
\end{align*}
where we used
$\Theta(V^\nu(\tau_e^\nu))=0=\Theta(V_*(\tau_e))$.  Therefore
\begin{equation}\label{eq:inner-exit-time-convergence}
	\lim_{\nu\to0}|\tau_e^\nu-\tau_e|=0.
\end{equation}

Finally, since $y_e^\nu=V^\nu(\tau_e^\nu)$ and
$y_e=V_*(\tau_e)$,
$$|y_e^\nu-y_e|\le|V^\nu(\tau_e^\nu)-V_*(\tau_e^\nu)|+|V_*(\tau_e^\nu)-V_*(\tau_e)|.$$
Both terms on the right tend to zero as $\nu\to0$. Hence $\lim_{\nu\to0}y_e^\nu=y_e$. Since $y_e$ lies in the interior of $K_e$ and
$y_e^\nu\to y_e$ as $\nu\to0$, one has $y_e^\nu\in K_e$ for all sufficiently small
$\nu$.
\end{proof}

\subsection{Subsequential limits}

Define the exit time, position, and velocity by
$$t_e^\nu:=t_{\rm ent}^\nu+\nu^{1-\alpha}\tau_e^\nu,\qquad x_e^\nu:=x^\nu(t_e^\nu),\qquad v_e^\nu:=v^\nu(t_e^\nu).$$
Equation \eqref{eq:inner-scale} and the uniform boundedness of the exit times
give
\begin{equation}\label{eq:inner-physical}
	t_e^\nu-t_{\rm ent}^\nu=O(\nu^{1-\alpha}),\qquad x_e^\nu-x_{\rm ent}^\nu=O(\nu^{2-\alpha}),\qquad|v_e^\nu|=\nu.
\end{equation}
Together with \eqref{eq:entry-time-scale}--\eqref{eq:entry-x-scale}, this gives
$\lim_{\nu\to0}(x_e^\nu,t_e^\nu,v_e^\nu)=(x_b,t_b,0)$.

For each $\varepsilon>0$, let $\mathcal S_+(\varepsilon)$ be the set of all solutions of \eqref{eq:intro-regularization-system} for which:
\begin{enumerate}
	\item $(x,v)\in C^0([t_b,t_b+\varepsilon];\mathbb R^{2d})\cap C^1((t_b,t_b+\varepsilon);\mathbb R^{2d})$;
	\item $x(t_b)=x_b$ and $v(t_b)=0$;
	\item $v(t)\ne0$ for every $t\in(t_b,t_b+\varepsilon]$;
	\item $\displaystyle\lim_{t\downarrow t_b}\frac{v(t)}{|v(t)|}=u_+$.
\end{enumerate}
On every exit interval on which $v^\nu(t)\ne0$, use the variables
\begin{align*}
	u^\nu(t)&:=\frac{v^\nu(t)}{|v^\nu(t)|},\\
	q^\nu(t)&:=|v^\nu(t)|^{1-\alpha},\\
	q_e^\nu&:=q^\nu(t_e^\nu)=\nu^{1-\alpha},\\
	\ell^\nu(t)&:=\log\frac{q^\nu(t)}{q_e^\nu}.
\end{align*}
For every $\delta>0$, write
$I_\delta:=[t_b-\delta,t_b+\delta]$.

\begin{proposition}
\label{prop:subsequence-limit}
Suppose the family of 
regularized trajector satisfies the entrance hypothesis
of Lemma~\ref{lem:inner}.  Then there is a $\delta>0$, independent of $\nu$,
such that every sequence $\nu_n\to0$ has a subsequence
$\nu_{n_k}\to0$ and
$(x_{\lim},v_{\lim})\in C^0(I_\delta;\R^{2d})$ satisfying
$$\lim_{k\to\infty}(x^{\nu_{n_k}},v^{\nu_{n_k}})=(x_{\lim},v_{\lim})\quad\text{in }C^0(I_\delta;\R^{2d}).$$
On the right of $t_b$, every such limit satisfies
$$(x_{\lim},v_{\lim})|_{[t_b,t_b+\delta]}\in\mathcal S_+(\delta).$$
\end{proposition}

\begin{proof}

First fix a sufficiently small $q_\rho>0$ such that
$\Lambda q_\rho\le a_+/2$.
Whenever $u\in N_+$ and $0<q\le q_\rho$, equation
\eqref{eq:q-physical} gives
\begin{equation}\label{eq:q-two-sided}
	c_+\le\frac{\dd q}{\dd t}\le C_+,\qquad c_+=\frac{(1-\alpha)a_+}{2}>0,
\end{equation}
where one may take
$$C_+=(1-\alpha)\left(\sup_{u\in N_+}|f_r(u)|+\Lambda q_\rho\right).$$

Use
$\ell^\nu=\log(q^\nu/q_e^\nu)$
as the independent variable. By \eqref{eq:q-two-sided}, $\ell^\nu$ increases monotonically with $t$ when $u\in N_+$ and $0<q\le q_\rho$.  Equation \eqref{eq:intro-ren-u} gives
\begin{equation}\label{eq:ell-u}
	\frac{\dd u^\nu}{\dd\ell^\nu}=\frac{f_s(u^\nu)+q^\nu\mathsf B(x^\nu,t)u^\nu}{(1-\alpha)[f_r(u^\nu)-\lambda(t)q^\nu]}.
\end{equation}
Because $f_r(u)\ge a_+$ and
$\lambda q\le\Lambda q_\rho\le a_+/2$, the denominator has a uniform positive
lower bound.  Since $f_r$ is bounded on $\mathcal K_+$, it also has a uniform
upper bound.  Thus there are $0<d_0\le d_1<\infty$ such that
$$d_0\le(1-\alpha)[f_r(u^\nu)-\lambda(t)q^\nu]\le d_1\qquad u^\nu\in\mathcal K_+,\ 0<q^\nu\le q_\rho.$$
Let $R_+(u)=\sqrt{W_+(u)}$.  When $R_+(u)>0$, equations
\eqref{eq:pm-W2}--\eqref{eq:pm-W3} imply
\begin{align}
	\frac{\dd R_+(u^\nu)}{\dd\ell^\nu}&=\frac{\nabla W_+(u^\nu)\cdot[f_s(u^\nu)+q^\nu\mathsf B(x^\nu,t)u^\nu]}{2R_+(u^\nu)(1-\alpha)[f_r(u^\nu)-\lambda(t)q^\nu]}\nonumber\\
	&\le-\frac{\mu_+}{d_1}R_+(u^\nu)+\frac{C_{W,+}M_B}{2d_0}q^\nu=-c_WR_+(u^\nu)+C_Bq^\nu,\label{eq:ell-Wplus}
\end{align}
where
$$c_W:=\frac{\mu_+}{d_1}>0,\qquad C_B:=\frac{C_{W,+}M_B}{2d_0}>0.$$
Reduce $q_\rho$ further so that
$C_Bq_\rho<c_W\rho_+/2$.
Then $\dd R_+(u^\nu)/\dd\ell^\nu$ is  strictly negative when $R_+\ge\rho_+/2$.  Starting from the interior of
$\mathcal K_+$, a direction with $q\le q_\rho$ cannot cross out of
$\mathcal K_+$, so \eqref{eq:q-two-sided} remains available.

We next choose $\delta$ for every sufficiently small
$\nu$.  
It holds that $q(t)\to0$ as $t\uparrow t_b$, so first
choose $\delta>0$ such that
$$q(t)<q_\rho,\qquad t_b-\delta\le t<t_b.$$
The constants $C_+$ and $q_\rho$ have already been fixed. Reduce $\delta$ so
that
$2C_+\delta<q_\rho/2$.
This $\delta$ depends only on the fixed constants,
not on $\nu$.
Since $t_e^\nu\to t_b$ and $q_e^\nu=\nu^{1-\alpha}\to0$ as $\nu\to0$, there is a
$\nu_0>0$ such that, for $0<\nu<\nu_0$,
$$t_b-\frac\delta2<t_e^\nu<t_b+\frac\delta2,\qquad q_e^\nu<\frac{q_\rho}{2}.$$

At the exit point, the direction $y_e^\nu$ belongs to the compact set $K_e$ contained
in the interior of $\mathcal K_+$, and the exit estimate gives
$q_e^\nu<q_\rho$. We claim that, for every sufficiently small $\nu$,
$q^\nu(t)<q_\rho$ and $u^\nu(t)\in\mathcal K_+$ throughout
$[t_e^\nu,t_b+\delta]$. Suppose the claim fails and define the first failure
time as
$$\tau_\nu:=\inf\bigl\{t\in[t_e^\nu,t_b+\delta]:q^\nu(t)=q_\rho\ \text{or}\
u^\nu(t)\in\partial\mathcal K_+\bigr\}.$$
By definition, on $[t_e^\nu,\tau_\nu)$ one has
$q^\nu<q_\rho$ and $u^\nu\in\mathcal K_+$, so the inward boundary estimate
for the direction and \eqref{eq:q-two-sided} both apply.

By \eqref{eq:ell-Wplus}, $\dd R_+(u^\nu)/\dd\ell^\nu$ is strictly negative on the directional boundary,
so the trajectory cannot first reach $\partial\mathcal K_+$ from the inside.
Also, integrating $\dd q^\nu/\dd t\le C_+$ from $t_e^\nu$ to
$\tau_\nu$ gives
$$q^\nu(\tau_\nu)\le q_e^\nu+C_+(\tau_\nu-t_e^\nu).$$
Because $\tau_\nu\le t_b+\delta$ and
$t_e^\nu>t_b-\delta/2$, we have
$\tau_\nu-t_e^\nu<3\delta/2<2\delta$.
The choices of $\delta$ and $\nu_0$ then give
$q^\nu(\tau_\nu)<q_\rho/2+2C_+\delta<q_\rho$.
Thus such $\tau_\nu$ cannot exist.  For every
$t\in[t_e^\nu,t_b+\delta]$,
$$q^\nu(t)<q_\rho,\qquad u^\nu(t)\in\mathcal K_+.$$
Hence, for every $0<\nu<\nu_0$ one has
$q^\nu(t)<q_\rho$ on the same fixed interval
$I_\delta:=[t_b-\delta,t_b+\delta]$.
Equivalently, $|v^\nu(t)|\le\rho$, where
$\rho:=q_\rho^{1/(1-\alpha)}$.

Since $H$ is continuous on $\overline{B_1}$, set
$$C_H:=\sup_{|V|\le1}|H(V)|<\infty,\qquad C_{F,0}:=C_H+\|f\|_\infty,\qquad C_{F,1}:=\max\{C_H,\|f\|_\infty\}.$$
On $I_\delta$, $|v^\nu|\le\rho$.  The position equation in \eqref{eq:regularized-system} gives
$$\left|\frac{\dd x^\nu}{\dd t}\right|=|v^\nu|\le\rho.$$
The velocity equation in \eqref{eq:regularized-system} and uniform boundedness of
$F_\nu$ on $\{|v|\le\rho\}$ give
$$\left|\frac{\dd v^\nu}{\dd t}\right|\le M_B\rho+\Lambda\rho+C_{F,1}\rho^\alpha.$$
Thus the family $\{(x^\nu,v^\nu)\}$ is uniformly bounded and equicontinuous on
$I_\delta$.  The derivative bounds show equicontinuity directly.  The Arzel\`a--Ascoli theorem therefore
applies.  For every sequence $\nu_n\to0$ there
is a strictly increasing sequence $n_1<n_2<\cdots$ such that
\begin{equation}\label{eq:subsequence-conv}
	\lim_{k\to\infty}(x^{\nu_{n_k}},v^{\nu_{n_k}})=(x_{\lim},v_{\lim})\quad\text{in }C^0(I_\delta;\R^{2d}).
\end{equation}

It remains to verify that $(x_{\lim},v_{\lim})$ solves the original singular equation. Define
$$F_0(v)=\begin{cases}
	|v|^\alpha f(v/|v|),&v\ne0,\\
	0,&v=0.
\end{cases}$$
Since $\alpha>0$, $F_0$ is continuous at $v=0$.
For $|v|>\nu$, $F_\nu(v)=F_0(v)$. For $|v|\le\nu$,
$$|F_\nu(v)-F_0(v)|\le C_{F,0}\nu^\alpha.$$
Consequently,
\begin{equation}\label{eq:F-uniform}
	\sup_{|v|\le\rho}|F_\nu(v)-F_0(v)|\le C_{F,0}\nu^\alpha,\qquad\lim_{\nu\to0}\sup_{|v|\le\rho}|F_\nu(v)-F_0(v)|=0.
\end{equation}
The integral equations for the regularized solutions are
\begin{align*}
	x^{\nu_{n_k}}(t)&=x^{\nu_{n_k}}(t_b-\delta)+\int_{t_b-\delta}^{t}v^{\nu_{n_k}}(\tau)\,\dd\tau,\\
	v^{\nu_{n_k}}(t)&=v^{\nu_{n_k}}(t_b-\delta)+\int_{t_b-\delta}^{t}\bigl[\mathsf B(x^{\nu_{n_k}}(\tau),\tau)v^{\nu_{n_k}}(\tau)-\lambda(\tau)v^{\nu_{n_k}}(\tau)+F_{\nu_{n_k}}(v^{\nu_{n_k}}(\tau))\bigr]\,\dd\tau.
\end{align*}
By \eqref{eq:subsequence-conv}, continuity of $\mathsf B$, and
\eqref{eq:F-uniform}, the integrands converge uniformly as $k\to\infty$.
Passing to the limit term by term gives
\begin{align*}
	x_{\lim}(t)&=x_{\lim}(t_b-\delta)+\int_{t_b-\delta}^{t}v_{\lim}(\tau)\,\dd\tau,\\
	v_{\lim}(t)&=v_{\lim}(t_b-\delta)+\int_{t_b-\delta}^{t}\bigl[\mathsf B(x_{\lim}(\tau),\tau)v_{\lim}(\tau)-\lambda(\tau)v_{\lim}(\tau)+F_0(v_{\lim}(\tau))\bigr]\,\dd\tau.
\end{align*}
Hence $(x_{\lim},v_{\lim})$ solves the original singular equation
\eqref{eq:intro-regularization-system} on $I_\delta$.

Now we prove that every subsequential limit belongs to
$\mathcal S_+(\delta)$.
By \eqref{eq:inner-physical},
$$\lim_{k\to\infty}\bigl(x^{\nu_{n_k}}(t_e^{\nu_{n_k}}),t_e^{\nu_{n_k}},v^{\nu_{n_k}}(t_e^{\nu_{n_k}})\bigr)=(x_b,t_b,0).$$
Together with uniform convergence, this gives
$$(x_{\lim}(t_b),v_{\lim}(t_b))=(x_b,0).$$

Fix any $t\in(t_b,t_b+\delta]$.  For sufficiently large $k$,
$t>t_e^{\nu_{n_k}}$.  Integrating \eqref{eq:q-two-sided} from
$t_e^{\nu_{n_k}}$ to $t$ gives
$$q_e^{\nu_{n_k}}+c_+(t-t_e^{\nu_{n_k}})\le q^{\nu_{n_k}}(t)\le q_e^{\nu_{n_k}}+C_+(t-t_e^{\nu_{n_k}}).$$
At this fixed time, $q_e^{\nu_{n_k}}\to0$,
$t_e^{\nu_{n_k}}\to t_b$, and $v^{\nu_{n_k}}(t)\to v_{\lim}(t)$ as
$k\to\infty$.  Passing to the limit
gives
\begin{equation}\label{eq:limit-q-linear}
	c_+(t-t_b)\le q_{\lim}(t)\le C_+(t-t_b),\qquad q_{\lim}(t):=|v_{\lim}(t)|^{1-\alpha}.
\end{equation}
In particular $q_{\lim}(t)>0$, so $v_{\lim}(t)\ne0$ to the right of $t_b$.

Finally, we prove that the direction of the limiting
solution tends to $u_+$. For
$t_b<t\leq t_b+\delta$,
define
$$u_{\lim}(t):=\frac{v_{\lim}(t)}{|v_{\lim}(t)|}.$$
Recall that
$R_+(u)=\sqrt{W_+(u)}$
and
$\ell^\nu(t)=\log(q^\nu(t)/q_e^\nu)$.
At the exit time, $\ell^\nu(t_e^\nu)=0$ and
$u^\nu(t_e^\nu)=y_e^\nu\in K_e$. Hence
$R_+(y_e^\nu)\leq C_e$ and $C_e:=\max_{u\in K_e}R_+(u)<\infty$.
Integrating \eqref{eq:ell-Wplus} from $\ell=0$ to
$\ell=\ell^\nu(t)$ gives
\begin{equation}\label{eq:one-orbit-direction}
	R_+(u^\nu(t))\leq C_e\left(\frac{q_e^\nu}{q^\nu(t)}\right)^{c_W}+\frac{C_B}{1+c_W}q^\nu(t),\qquad t_e^\nu\leq t\leq t_b+\delta.
\end{equation}
Fix $t\in(t_b,t_b+\delta]$. For all sufficiently large $k$,
$t>t_e^{\nu_{n_k}}$, so \eqref{eq:one-orbit-direction} applies with
$\nu=\nu_{n_k}$. Moreover, as $k\to\infty$,
$$q_e^{\nu_{n_k}}\to0,\qquad q^{\nu_{n_k}}(t)\to q_{\lim}(t)>0,$$
and uniform convergence of the velocities implies
$u^{\nu_{n_k}}(t)\to u_{\lim}(t)$.
Passing to the limit in \eqref{eq:one-orbit-direction} gives
$$R_+(u_{\lim}(t))\leq\frac{C_B}{1+c_W}q_{\lim}(t).$$
As $t\downarrow t_b$, equation \eqref{eq:limit-q-linear} gives
$q_{\lim}(t)\to0$. Therefore
$W_+(u_{\lim}(t))\to0$.
By \eqref{eq:pm-W1},
$\lim_{t\downarrow t_b}u_{\lim}(t)=u_+$.

Every subsequential limit therefore satisfies
$$(x_{\lim}(t_b),v_{\lim}(t_b))=(x_b,0),\qquad\lim_{t\downarrow t_b}u_{\lim}(t)=u_+.$$
and $v_{\lim}(t)\ne0$ for $t>t_b$. This shows that
$(x_{\lim},v_{\lim})|_{[t_b,t_b+\delta]}\in\mathcal S_+(\delta)$.
\end{proof}

\subsection{Uniqueness of the continuation}
\label{sec:reg-outgoing}
In this subsection, we prove that there is only one solution in $\mathcal S_+(\delta)$.

The renormalized time used below is not a continuation of the
pre-blowup renormalized time. We introduce a new renormalized time
$s_{\rm out}$ also satisfying
$$\frac{\dd t}{\dd s_{\rm out}}=q,\qquad\frac{\dd s_{\rm out}}{\dd t}=\frac1q.$$
Choose a fixed time $t_{\rm ref}>t_b$ in the domain of the outgoing
solution and define
\begin{equation}\label{eq:outgoing-s-normalization}
	s_{\rm out}(t):=\int_{t_{\rm ref}}^t\frac{\dd\tau}{q(\tau)}.
\end{equation}
Then $s_{\rm out}(t_{\rm ref})=0$. We will prove below that
$\lim_{t\downarrow t_b}s_{\rm out}(t)=-\infty$. For simplicity, throughout this subsection we still write $s=s_{\rm out}$.

\begin{proposition}\label{prop:Splus-unique}\label{prop:outgoing-unique}
Let $0<\alpha<1$ and $f\in C^2(\Sph;\R^d)$. There is an
$\varepsilon>0$ such that the set $\mathcal S_+(\varepsilon)$ contains exactly one solution, denoted by
$(x^+,v^+)$.  
\end{proposition}

\begin{proof}
The proof has five steps.  Steps~1--2 construct a $C^1$ extended system and
identify its unique strictly unstable direction.  Step~3 shows that any local solution in $\mathcal S_+(\varepsilon)$, if
one exists, has a renormalized trajectory converging
exponentially to the equilibrium in backward time.  Step~4 constructs
such a solution from the strong unstable manifold, and Step~5 uses its
characterization property to prove uniqueness.

\medskip
\noindent\textbf{Step 1: construct a $C^1$ extended system.}

Equations \eqref{eq:intro-ren-u}--\eqref{eq:intro-ren-tx} give, in the region $v\ne0$,
\begin{align}
	\frac{\dd u}{\dd s}&=f_s(u)+q\mathsf B(x,t)u,\nonumber\\
	\frac{\dd q}{\dd s}&=(1-\alpha)q[f_r(u)-\lambda(t)q],\nonumber\\
	\frac{\dd x}{\dd s}&=q^\beta u,\nonumber\\
	\frac{\dd t}{\dd s}&=q,\label{eq:outgoing-ren-system}
\end{align}
where $\beta=(2-\alpha)/(1-\alpha)>2$.
Choose a smooth local coordinate block near $u_+$ whose closure is contained
in $N_+$, and write
$$u=\Phi(\xi),\qquad\Phi(0)=u_+.$$
Here $D\Phi(\xi)$ maps coordinate vectors to tangent vectors on the sphere,
and $(D\Phi(\xi))^{-1}$ below is its coordinate inverse on that tangent space.
Set
$$X=x-x_b,\qquad\theta=t-t_b.$$
The chain rule gives
$$\frac{\dd u}{\dd s}=D\Phi(\xi)\frac{\dd\xi}{\dd s}.$$
Define
\begin{align*}
	F(\xi)&:=(D\Phi(\xi))^{-1}f_s(\Phi(\xi)),\\
	M(\xi,X,\theta)&:=(D\Phi(\xi))^{-1}\mathsf B(x_b+X,t_b+\theta)\Phi(\xi),\\
	r(\xi)&:=f_r(\Phi(\xi)).
\end{align*}
Substituting this into the directional equation in
\eqref{eq:outgoing-ren-system} yields
\begin{align}
	\frac{\dd\xi}{\dd s}&=F(\xi)+qM(\xi,X,\theta),\label{eq:unstable-exact-xi}\\
	\frac{\dd q}{\dd s}&=(1-\alpha)q\bigl[r(\xi)-\lambda(t_b+\theta)q\bigr],\label{eq:unstable-exact-q}\\
	\frac{\dd X}{\dd s}&=q^\beta\Phi(\xi),\label{eq:unstable-exact-X}\\
	\frac{\dd\theta}{\dd s}&=q.\label{eq:unstable-exact-theta}
\end{align}
Since $f_s(u_+)=0$, $F(0)=0$. Set
\begin{equation}\label{eq:unstable-linear-coefficients}
	A:=DF(0),\qquad b:=M(0,0,0),\qquad r_+:=r(0)=f_r(u_+),\qquad a_0:=(1-\alpha)r_+>0.
\end{equation}
Because $f\in C^2$ and
$\mathsf B,\lambda\in C^1$, Taylor's formula yields
a constant $C$, independent of $(\xi,q,X,\theta)$, such that
\begin{align}
	|F(\xi)-A\xi|&\le C|\xi|^2,\label{eq:unstable-Taylor-F}\\
	|M(\xi,X,\theta)-b|&\le C(|\xi|+|X|+|\theta|),\label{eq:unstable-Taylor-M}\\
	|r(\xi)-r_+|&\le C|\xi|.\label{eq:unstable-Taylor-r}
\end{align}
Define the remainders
\begin{align*}
	\mathcal N_\xi(\xi,q,X,\theta)&:=[F(\xi)-A\xi]+q[M(\xi,X,\theta)-b],\\
	\mathcal N_q(\xi,q,X,\theta)&:=(1-\alpha)q[r(\xi)-r_+]-(1-\alpha)\lambda(t_b+\theta)q^2.
\end{align*}
Then \eqref{eq:unstable-exact-xi}--
\eqref{eq:unstable-exact-theta} are equivalent to
\begin{align}
	\frac{\dd\xi}{\dd s}&=A\xi+bq+\mathcal N_\xi(\xi,q,X,\theta),\label{eq:unstable-local-xi}\\
	\frac{\dd q}{\dd s}&=a_0q+\mathcal N_q(\xi,q,X,\theta),\label{eq:unstable-local-q}\\
	\frac{\dd X}{\dd s}&=q^\beta\Phi(\xi),\label{eq:unstable-local-X}\\
	\frac{\dd\theta}{\dd s}&=q.\label{eq:unstable-local-theta}
\end{align}
Equations \eqref{eq:unstable-Taylor-F}--\eqref{eq:unstable-Taylor-r} give
\begin{align}
	|\mathcal N_\xi|&\le C\bigl[|\xi|^2+|q|(|\xi|+|X|+|\theta|)\bigr],\label{eq:unstable-Nxi-bound}\\
	|\mathcal N_q|&\le C(|q||\xi|+q^2).\label{eq:unstable-Nq-bound}
\end{align}

The standard invariant manifold
theorem is formulated in an open neighborhood of the equilibrium point.  So we extend
$q^\beta$ in \eqref{eq:unstable-local-X} to negative $q$ by
\begin{equation}\label{eq:unstable-q-extension}
	\widehat p_\beta(q):=q|q|^{\beta-1}.
\end{equation}
Since $\beta>2$, this function is at least $C^2$ at $q=0$, and it equals
$q^\beta$ for $q\ge0$.  The remaining terms extend to small negative $q$ by
the same formulas.  

 The specified point $(x_b,t_b)$
corresponds to
$E_b:=(0,0,0,0)$,
which is an equilibrium point of system \eqref{eq:unstable-exact-xi}--\eqref{eq:unstable-exact-theta}.
We have obtained a $C^1$ autonomous system on an open neighborhood of $E_b$,
whose positive-$q$ part is the original system.  The extension to
$q<0$ is used only to apply the standard invariant manifold theorem.

\medskip
\noindent\textbf{Step 2: $E_b$ has exactly one strictly unstable direction.}

With the variables ordered as $(\xi,q,X,\theta)$, set
$Z:=(\xi,q,X,\theta)^{\mathsf T}$.  The linearized system at $E_b$ is
\begin{equation}\label{eq:unstable-full-linearization}
	\frac{\dd Z}{\dd s}=LZ,\qquad L:=\begin{pmatrix}
	A&b&0&0\\
	0&a_0&0&0\\
	0&0&0&0\\
	0&1&0&0
\end{pmatrix}.
\end{equation}
The third equation of the linearized system has zero linear part because $\beta>2$. By the local exponential stability of $u_+$, every
eigenvalue of the coordinate linearization $A=DF(0)$ has strictly negative
real part.  

The spectrum of the coefficient matrix $L$ consists of the eigenvalues of $A$ with strictly negative real
parts, one positive eigenvalue $a_0$, and zero eigenvalues in the
$(X,\theta)$ directions.  Thus there is exactly one positive eigenvalue.

\medskip
\noindent\textbf{Step 3: every possible local outgoing solution converges
exponentially backward to $E_b$.}

Let $(x,v)$ be any local solution satisfying the definition
of $\mathcal S_+$ on a right neighborhood of $t_b$.  Since $v(t)\to0$ and
$v(t)/|v(t)|\to u_+$ as $t\downarrow t_b$, we may write
$u(t)=\Phi(\xi(t))$ near $t_b$, with $\xi(t)\to0$.  Equation
\eqref{eq:q-physical}, boundedness of $\lambda$, and
$f_r(u(t))\to f_r(u_+)$ give
\begin{equation}\label{eq:unstable-any-dqdt-limit}
	\lim_{t\downarrow t_b}\frac{\dd q(t)}{\dd t}=a_0.
\end{equation}
Choose $t_1>t_b$ so that
$$\frac{a_0}{2}\le\frac{\dd q(t)}{\dd t}\le2a_0,\qquad t_b<t\le t_1.$$
Since $q(t_b)=0$, integration yields
\begin{equation}\label{eq:unstable-any-q-t-bounds}
	\frac{a_0}{2}(t-t_b)\le q(t)\le2a_0(t-t_b),\qquad t_b<t\le t_1.
\end{equation}
Apply \eqref{eq:outgoing-s-normalization} to this solution with
$t_{\rm ref}=t_1$.  Thus
\begin{equation}\label{eq:unstable-any-s-integral}
	s(t)=-\int_t^{t_1}\frac{\dd\tau}{q(\tau)}.
\end{equation}
The upper bound in \eqref{eq:unstable-any-q-t-bounds} shows that the
integral \eqref{eq:unstable-any-s-integral} diverges to $+\infty$ as $t\downarrow t_b$.  Hence
\begin{equation}\label{eq:unstable-any-s-minus-infty}
	\lim_{t\downarrow t_b}s(t)=-\infty.
\end{equation}
Thus the solution $(x,v)\in\mathcal S_+$ gives an extended trajectory with
$(\xi,q,X,\theta)\to E_b$ as $s\to-\infty$.

We now strengthen convergence to exponential convergence, beginning
with $q$.  Equation \eqref{eq:unstable-exact-q} gives
$$\frac1q\frac{\dd q}{\dd s}=(1-\alpha)[r(\xi)-\lambda(t_b+\theta)q],\qquad\lim_{s\to-\infty}\frac1q\frac{\dd q}{\dd s}=a_0.$$
Thus there is  an $s_1$ such that, for all $s\le s_1$,
\begin{equation}\label{eq:unstable-any-logq-bounds}
	\frac{a_0}{2}\le\frac1q\frac{\dd q}{\dd s}\le2a_0.
\end{equation}
For $\rho\le s\le s_1$, integrate the left inequality from $\rho$ to $s$:
$$\log q(s)-\log q(\rho)\ge\frac{a_0}{2}(s-\rho).$$
Therefore
\begin{equation}\label{eq:unstable-any-q-exp}
	q(\rho)\le q(s)\ee^{-a_0(s-\rho)/2},\qquad\rho\le s\le s_1.
\end{equation}

We first show that $\xi=O(q)$.  Write \eqref{eq:unstable-local-xi} as
$$\frac{\dd\xi}{\dd s}=A\xi+G(s),\qquad G(s):=bq(s)+\mathcal N_\xi(\xi,q,X,\theta)(s).$$
For any $r<s\le s_1$, variation of constants gives
$$\xi(s)=\ee^{A(s-r)}\xi(r)+\int_r^s\ee^{A(s-\sigma)}G(\sigma)\,\dd\sigma.$$
Since $(\xi,q,X,\theta)\to0$ as $s\to-\infty$, for any small $\epsilon>0$ we may take $s_1$
more negative so that
$$|\xi(s)|+|q(s)|+|X(s)|+|\theta(s)|\le\epsilon,\qquad s\le s_1.$$
After increasing $C$ if necessary, \eqref{eq:unstable-Nxi-bound} yields
\begin{equation}\label{eq:unstable-any-G-bound}
	|G(s)|\le C\epsilon|\xi(s)|+Cq(s),\qquad s\le s_1.
\end{equation}
Since every eigenvalue of $A$ has strictly negative real part, there are
$C_A\ge1$ and $\gamma>0$ such that
$$\|\ee^{A(s-r)}\|\le C_A\ee^{-\gamma(s-r)}.$$
As $r\to-\infty$, $\xi(r)\to0$, and therefore
$$|\ee^{A(s-r)}\xi(r)|\le C_A\ee^{-\gamma(s-r)}|\xi(r)|,\qquad\lim_{r\to-\infty}|\ee^{A(s-r)}\xi(r)|=0.$$
Letting $r\to-\infty$ gives
\begin{equation}\label{eq:unstable-any-xi-integral}
	\xi(s)=\int_{-\infty}^s\ee^{A(s-\sigma)}G(\sigma)\,\dd\sigma.
\end{equation}
Define
$Y(s):=\sup_{\rho\le s}|\xi(\rho)|$.
For $\tau\le s\le s_1$, equations
\eqref{eq:unstable-any-xi-integral},
\eqref{eq:unstable-any-G-bound}, and \eqref{eq:unstable-any-q-exp} give
\begin{align*}
	|\xi(\tau)|&\le C\epsilon Y(s)\int_{-\infty}^{\tau}\ee^{-\gamma(\tau-\sigma)}\,\dd\sigma+Cq(\tau)\int_{-\infty}^{\tau}\ee^{-(\gamma+a_0/2)(\tau-\sigma)}\,\dd\sigma\\
	&\le\frac{C\epsilon}{\gamma}Y(s)+\frac{C}{\gamma+a_0/2}q(s).
\end{align*}
Taking the supremum over $\tau\le s$ yields
$$Y(s)\le\frac{C\epsilon}{\gamma}Y(s)+\frac{C}{\gamma+a_0/2}q(s).$$
Choose $\epsilon$ small enough so that $C\epsilon/\gamma\le1/2$.  Then
\begin{equation}\label{eq:unstable-any-xi-Oq}
	|\xi(s)|\le Y(s)\le\frac{2C}{\gamma+a_0/2}q(s),\qquad s\le s_1.
\end{equation}

The lower bound in \eqref{eq:unstable-any-q-t-bounds} gives
$0\le\theta=t-t_b\le2q/a_0$.  Since $X(s)\to0$ as $s\to-\infty$ and
$\dd X/\dd s=q^\beta\Phi(\xi)$,
$$X(s)=\int_{-\infty}^sq(\rho)^\beta\Phi(\xi(\rho))\,\dd\rho.$$
Because $|\Phi(\xi)|=1$, equation \eqref{eq:unstable-any-q-exp} gives
\begin{equation}\label{eq:unstable-any-X-Oqbeta}
	|X(s)|\le q(s)^\beta\int_{-\infty}^s\ee^{-\beta a_0(s-\rho)/2}\,\dd\rho=\frac{2}{\beta a_0}q(s)^\beta.
\end{equation}

After $s_1$ has been fixed, \eqref{eq:unstable-any-logq-bounds} gives a
constant $C_1>0$ such that
$q(s)\le C_1\ee^{a_0s/2}$ for $s\le s_1$.  Combining this with
\eqref{eq:unstable-any-xi-Oq}, $\theta=O(q)$, and
\eqref{eq:unstable-any-X-Oqbeta} yields a constant $C_2>0$ such that
\begin{equation}\label{eq:unstable-any-full-exp}
	|\xi(s)|+|q(s)|+|X(s)|+|\theta(s)|\le C_2\ee^{a_0s/2},\qquad s\le s_1.
\end{equation}

\medskip
\noindent\textbf{Step 4: construct a candidate solution from the strong
unstable manifold.}

  If the linearization of a $C^1$ autonomous
system at an equilibrium has a stable, center, and unstable spectral splitting
with the unstable spectrum separated positively from the imaginary axis, the strong unstable manifold theorem~\cite{CapinskiWasieczko2015,HirschPughShub1977} gives a unique local strong unstable manifold
$W^{uu}_{\rm loc}$ with the following properties:
\begin{enumerate}
 \item it passes through the equilibrium and is tangent there to the unstable
 subspace;
 \item it is locally invariant under the flow, and its points converge
 exponentially to the equilibrium in backward time;
 \item conversely, every trajectory that converges exponentially to the equilibrium in backward time
 belongs to this local strong unstable manifold.
\end{enumerate}
The theorem gives a unique one-dimensional local manifold
$W^{uu}_{\rm loc}(E_b)$.  On a connected half-branch of a one-dimensional
manifold where the vector field never vanishes, ordinary local ODE uniqueness
makes the local trajectories through different points join into the whole
half-branch.  Thus the positive-$q$ half branch is one trajectory.

Choose the half branch of $W^{uu}_{\rm loc}(E_b)$ on which $q>0$.
Let $(\xi(s),q(s),X(s),\theta(s))$ denote its backward orbit.  By the
characterization of the strong unstable manifold,
$$(\xi(s),q(s),X(s),\theta(s))\to E_b\qquad\text{as }s\to-\infty.$$
Along this half-branch,
$$\frac{\dd t}{\dd s}=q(s)>0,\qquad t(s)=t_b+\theta(s),$$
so $t(s)$ is strictly increasing and maps $(-\infty,s_1]$ onto
$(t_b,t_1]$ for a suitable $s_1$.  Let $s=s(t)$ denote its inverse and set
$$x^+(t):=x_b+X(s(t)),\qquad v^+(t):=q(s(t))^{1/(1-\alpha)}\Phi(\xi(s(t))).$$
  Hence $(x^+,v^+)$ solves system \eqref{eq:intro-regularization-system} on $(t_b,t_1]$.  Moreover, as $t\downarrow t_b$,
$$x^+(t)\to x_b,\qquad v^+(t)\to0,\qquad\frac{v^+(t)}{|v^+(t)|}\to u_+.$$
Thus $(x^+,v^+)\in\mathcal S_+(t_1-t_b)$. Therefore $W^{uu}_{\rm loc}(E_b)$ determines an element of
$\mathcal S_+(\varepsilon)$ for some $\varepsilon>0$.

\medskip

\noindent\textbf{Step 5: uniqueness from the characterization property of
the strong unstable manifold.}

By the third property of the strong unstable manifold theorem and
\eqref{eq:unstable-any-full-exp}, the local trajectory associated with every
element of $\mathcal S_+(\varepsilon)$ must lie in
$W^{uu}_{\rm loc}(E_b)$.  We already proved that the positive-$q$ half branch
of this manifold is a single trajectory.  Thus all elements of
$\mathcal S_+(\varepsilon)$ correspond to the same extended trajectory, up to
a translation of renormalized time $s$.

\end{proof}

Under the hypotheses of Propositions~\ref{prop:subsequence-limit} and
\ref{prop:Splus-unique}, reduce the number $\delta$ from
Proposition~\ref{prop:subsequence-limit} so that it does not exceed the
$\varepsilon$ supplied by Proposition~\ref{prop:Splus-unique}. Define
$$(x^-(t),v^-(t)):=\begin{cases}
	(x(t),v(t)),&t_0\le t<t_b,\\
	(x_b,0),&t=t_b.
\end{cases}$$
Then $(x^-,v^-)\in C^0([t_0,t_b];\mathbb R^{2d})$. Denote the
unique element of $\mathcal S_+(\delta)$ by $(x^+,v^+)$ and define the
candidate full limit on $I_\delta=[t_b-\delta,t_b+\delta]$ by
$$(x^*(t),v^*(t)):=\begin{cases}
	(x^-(t),v^-(t)),&t_b-\delta\le t\le t_b,\\
	(x^+(t),v^+(t)),&t_b\le t\le t_b+\delta.
\end{cases}$$
The two branches have the common value $(x_b,0)$ at $t_b$, so this definition
gives a continuous curve on $I_\delta$.

\begin{proposition}\label{prop:local-full-convergence}
With the notation fixed above, the full regularization family satisfies
\begin{equation}\label{eq:local-full-conv}
	\lim_{\nu\to0}(x^\nu,v^\nu)=(x^*,v^*)\quad\text{in }C^0(I_\delta;\R^{2d}).
\end{equation}
\end{proposition}

\begin{proof}
Suppose, on the contrary, that the full regularization family does not
converge to $(x^*,v^*)$ on $I_\delta$.  Then there exist an $\varepsilon_0>0$
and a sequence $\nu_n\to0$ such that
$$\sup_{t\in I_\delta}\bigl(|x^{\nu_n}(t)-x^*(t)|+|v^{\nu_n}(t)-v^*(t)|\bigr)\ge\varepsilon_0.$$
Proposition~\ref{prop:subsequence-limit} supplies a subsequence of $\nu_n$ whose
regularized solutions converge uniformly on $I_\delta$, with a limit lying in
$\mathcal S_+(\delta)$ to the right of $t_b$ and equal to $(x^-,v^-)$ to the
left of $t_b$, contradicting the positive lower bound above. 
\end{proof}

\subsection{Proof of the main result}

\begin{proof}[Proof of Theorem~\ref{thm:intro-regularization-selection}]
	Theorem~\ref{thm:pre} and
	Proposition~\ref{prop:pre-reconstruction} give a finite blowup time
	$t_b$ and a point $x_b\in\R^d$ such that 
	$$\lim_{t\uparrow t_b}x(t)=x_b,\qquad\lim_{t\uparrow t_b}v(t)=0,\qquad\lim_{t\uparrow t_b}\frac{v(t)}{|v(t)|}=u_-.$$
	Proposition~\ref{prop:Splus-unique} gives a unique local solution
	$(x^+,v^+)\in\mathcal S_+(\delta)$. Hence
	$$x^+(t_b)=x_b,\qquad v^+(t_b)=0,$$
	and for $t_b<t\leq t_b+\delta$,
	$$v^+(t)\neq0,\qquad\lim_{t\downarrow t_b}\frac{v^+(t)}{|v^+(t)|}=u_+.$$
	Extend the curve defined on
	$I_\delta=[t_b-\delta,t_b+\delta]$ to
	$[t_0,t_b+\delta]$ by setting
	$$(x^*(t),v^*(t)):=\begin{cases}
	(x^-(t),v^-(t)),&t_0\leq t\leq t_b,\\
	(x^+(t),v^+(t)),&t_b\leq t\leq t_b+\delta.
\end{cases}$$
	The two branches have the common value $(x_b,0)$ at $t_b$, so
	$$(x^*,v^*)\in C^0([t_0,t_b+\delta];\R^{2d}).$$
	Proposition~\ref{prop:local-full-convergence} gives as $\nu\to0$
	$$(x^\nu,v^\nu)\longrightarrow(x^*,v^*)\quad\text{in } C^0([t_b-\delta,t_b+\delta];\R^{2d}).$$
	
It remains only to combine the local convergence with the incoming
trajectory away from the collision. By \eqref{eq:entry-time-scale},
$t_{\rm ent}^\nu\to t_b$ as $\nu\to0$.
Hence, for all sufficiently small $\nu$,
$t_{\rm ent}^\nu>t_b-\frac{\delta}{2}$.
Since the original and regularized trajectories coincide up to
$t_{\rm ent}^\nu$, one has
$$(x^\nu,v^\nu)=(x^-,v^-)\qquad\text{on }[t_0,t_b-\delta/2].$$
Combining this identity with the convergence on
$I_\delta=[t_b-\delta,t_b+\delta]$ gives
$$(x^\nu,v^\nu)\longrightarrow(x^*,v^*)\quad\text{in } C^0([t_0,t_b+\delta];\R^{2d}).$$
\end{proof}

\subsection{Planar example}

\begin{example}
	\label{ex:planar-magnetic}
	
	Let $d=2$, fix $0<\alpha<1$ and $b\neq0$, and set
	$$J=\begin{pmatrix}
	0&-1\\
	1&0
\end{pmatrix},\qquad\mathsf B(x,t)\equiv bJ,\qquad\lambda(t)\equiv0.$$
	For $u=(u_1,u_2)\in\mathbb S^1$, define
	$$f(u):=u_1u-2u_1u_2Ju.$$
	Then
	$$f_r(u)=u_1,\qquad f_s(u)=-2u_1u_2Ju.$$
Hence
	$u_\pm=(\pm1,0)$
	are locally exponentially stable equilibria of equation $\dd u/\dd s=f_s(u)$, while
	$f_r(u_\pm)=\pm1$.
	On sufficiently small neighborhoods of $u_\pm$, 
	$$W_+(u)=1-u_1,\qquad W_-(u)=1+u_1$$
	are the Lyapunov functions of $u_\pm$.
	
	We next give an explicit $C^1$ regularization.  Let
	$$\chi(r):=3r^2-2r^3,\qquad0\le r\le1,$$
	and, for $V\ne0$, set
	$$\mathcal G(V):=|V|^\alpha f\!\left(\frac{V}{|V|}\right).$$
	Define
	$$H(V):=\chi(|V|)\mathcal G(V)+\bigl[1-\chi(|V|)\bigr]\,2e_1,\qquad0<|V|\le1,\qquad H(0)=2e_1,$$
	where $e_1=(1,0)^T$.
	Since $\chi(r)=O(r^2)$ as $r\downarrow0$, the first term is $C^1$ at
	$V=0$.  Since $\chi(1)=1$ and $\chi'(1)=0$, the values and first
	derivatives of $H$ and $\mathcal G$ agree on $|V|=1$.  Thus the patched
	field $F_\nu$ defined by \eqref{eq:Fnu} is $C^1$.
	
	The reference orbit entering at $u_-=-e_1$ remains on the
	$e_1$ axis:
	$$V_*(\tau)=\zeta(\tau)e_1,\qquad\zeta(0)=-1,$$
	where
	$$\dot\zeta=\chi(|\zeta|)|\zeta|^\alpha+2\bigl[1-\chi(|\zeta|)\bigr]>0.$$
	
  The
	inner equation for the regularized problem is
	$$\frac{\dd V^\nu}{\dd\tau}=H(V^\nu)+\nu^{1-\alpha}bJV^\nu,$$
	so its magnetic contribution vanishes in the unit-scale limit.   In particular, the trajectory satisfying the entrance condition has a unique
	selected continuation through $v=0$, and the selected outgoing direction
	is $u_+=e_1$.

	For the numerical illustrations below, we set $\alpha=1/3$ and $b=4$.
	The plotted vector field is precisely the velocity equation
	\begin{equation}\label{eq:planar-plot-vector-field}
	\dot v=4Jv+F_\nu(v),
\end{equation}
	where $F_\nu$ is the patched field in \eqref{eq:Fnu}, with the present
	choices of $f$ and $H$.  Equivalently, if $u=v/|v|$ for $|v|>\nu$ and
	$V=v/\nu$ for $|v|\le\nu$, then
	$$F_\nu(v)=\begin{cases}
	|v|^{1/3}\bigl(u_1u-2u_1u_2Ju\bigr),&|v|>\nu,\\
	\nu^{1/3}\!\left[\chi(|V|)|V|^{1/3}f\!\left(\dfrac{V}{|V|}\right)+\bigl(1-\chi(|V|)\bigr)2e_1\right],&0<|v|\le\nu,\\
	\nu^{1/3}2e_1,&v=0.
\end{cases}$$
	where $f(u)=u_1u-2u_1u_2Ju$. Thus the displayed curves are velocity trajectories only; the $x$-equation
	does not enter the computation because the magnetic matrix is constant and
	the force is independent of $x$.  They share the initial value
	$$v(0)=0.25\bigl(\cos(\pi-0.15),\sin(\pi-0.15)\bigr)^\mathsf T$$
	and are integrated on $0\le t\le1$.

\begin{figure}[H]
	\centering
	\includegraphics[width=0.74\textwidth]{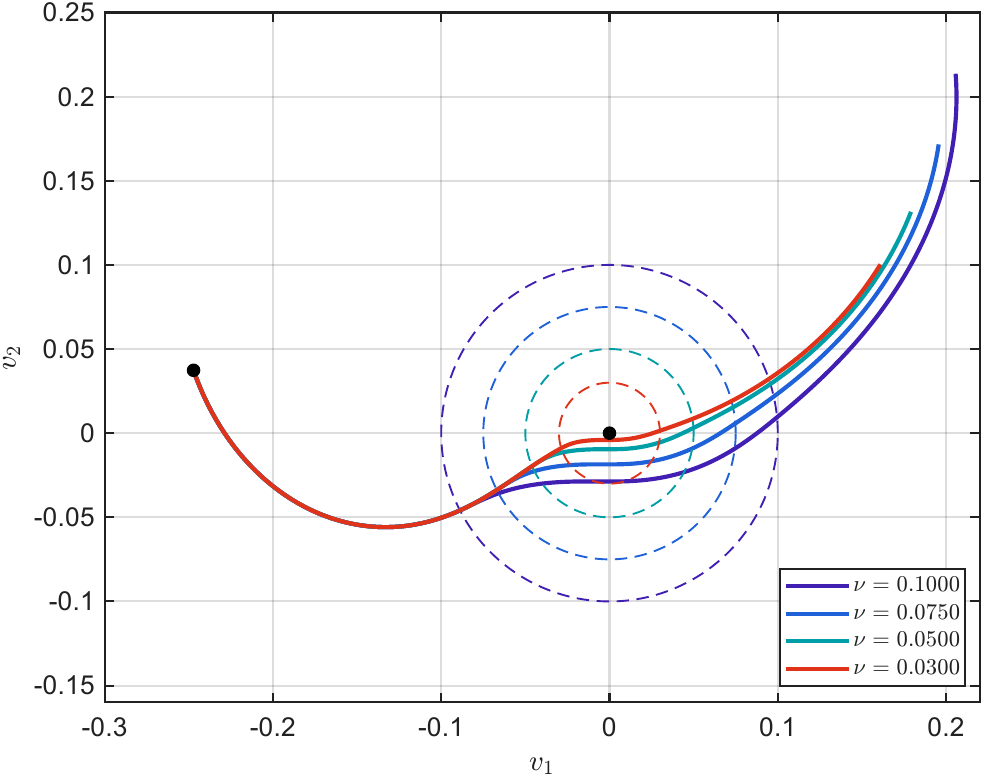}
	\caption{Velocity trajectories of \eqref{eq:planar-plot-vector-field} in the $v$-plane for $\nu=0.10$ (purple), $\nu=0.075$ (blue),
$\nu=0.05$ (cyan), and $\nu=0.03$ (red).  The dashed circles have radii $\nu$ and mark the
	corresponding regularization boundaries.}
	\label{fig:planar-magnetic-full-plane}
	\end{figure}

\begin{figure}[H]
	\centering
	\includegraphics[width=0.74\textwidth]{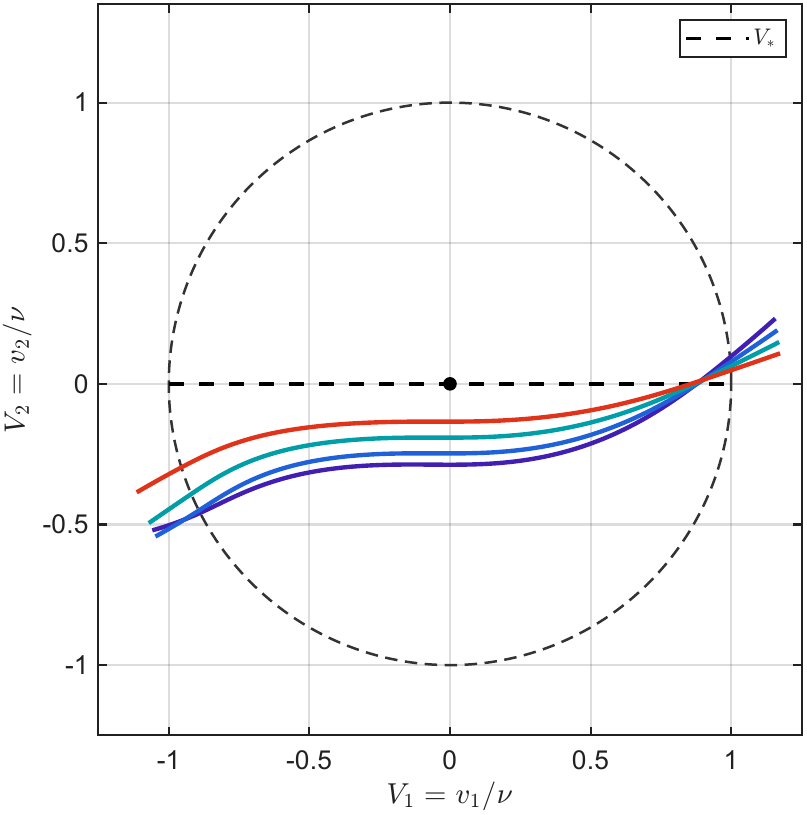}
	\caption{The same solutions in the rescaled variable $V=v/\nu$ near the
	regularization layer.  The dashed circle is $|V|=1$; each curve is shown
	up to $|V|\le1.18$, so that both its entrance and exit across the layer are
	visible.  The black dashed segment is the reference orbit $V_*$.}
	\label{fig:planar-magnetic-inner-layer}
	\end{figure}
\end{example}

\section{Parameterized dynamics and chaos}

Unless explicitly stated otherwise, the variables and parameters introduced
in this section are independent of those used in the preceding sections. In this section, we assume that  \textup{(H1)}--\textup{(H4)} hold.

\subsection{Stereographic projection}
\label{sec:chaos-stereographic}

We assume first that the magnetic field and the tangential force have a fixed
common equilibrium direction $u_\star\in\mathbb S^{d-1}$, namely
\begin{equation}\label{common-directional-equilibrium}
	\mathsf B(x,t)u_\star=0,\qquad f_s(t,u_\star,\varepsilon)=0
\end{equation}
for all relevant $(x,t,\varepsilon)$.   A fixed
orthogonal change of coordinates preserves skew-symmetry and tangency and
sends this common equilibrium to the South Pole.  Hence, without loss of
generality, we take
$u_\star=(-1,0,\ldots,0)$.
We then use the stereographic chart from the North Pole, so that the chosen
common equilibrium is mapped to the coordinate origin.  Writing
$u=(u^{(0)},\widehat u)$ where $u^{(0)}\in\mathbb{R}$ and $\widehat{u}\in\mathbb{R}^{d-1}$, set
\begin{equation}\label{eq:stereographic}
	X=\frac{\widehat u}{1-u^{(0)}}\in\R^{n},\qquad n=d-1,
\end{equation}
with inverse
$$u^{(0)}(X)=\frac{|X|^2-1}{1+|X|^2},\qquad\widehat u(X)=\frac{2X}{1+|X|^2}.$$
In particular,
$$\pi(u_\star)=0,\qquad u(0)=u_\star.$$
For a tangent vector $z=(z^{(0)},\widehat z)\in T_u\Sph$, differentiation of
$\pi(u)=\widehat u/(1-u^{(0)})$ gives
\begin{equation}\label{eq:stereographic-derivative}
	D_u\pi(u)z=\frac{\widehat z+Xz^{(0)}}{1-u^{(0)}}=\frac{1+|X|^2}{2}\bigl(\widehat z+Xz^{(0)}\bigr).
\end{equation}
By the first condition in \eqref{common-directional-equilibrium}, the
restriction of $\mathsf B(x,t)$ to the tangent coordinates is a real
skew-symmetric matrix $\mathsf B_0(x,t)$.  Indeed, skew-symmetry together
with $\mathsf B(x,t)u_\star=0$ gives the block form
$$\mathsf B(x,t)=\begin{pmatrix}
	0&0\\
	0&\mathsf B_0(x,t)
\end{pmatrix},$$
and hence
$$D_u\pi(u(X))\mathsf B(x,t)u(X)=\mathsf B_0(x,t)X.$$
Define the projected tangential force by
$$\mathcal N(t,X,\varepsilon)=D_u\pi(u(X))f_s(t,u(X),\varepsilon).$$
Equation
\eqref{eq:intro-II-direction} becomes
\begin{equation}\label{eq:projected-system}
	\frac{\dd X}{\dd t}=\mathsf B_0(x,t)X+w^{\alpha-1}(t)\mathcal N(t,X,\varepsilon).
\end{equation}
Writing $f_s=((f_s)^{(0)},\widehat f_s)$, formula
\eqref{eq:stereographic-derivative} yields the explicit expression
\begin{equation}\label{eq:projected-N-explicit}
	\mathcal N(t,X,\varepsilon)=\frac{1+|X|^2}{2}\left[\widehat f_s(t,u(X),\varepsilon)+X(f_s)^{(0)}(t,u(X),\varepsilon)\right].
\end{equation}
The second condition in \eqref{common-directional-equilibrium} gives
$$\mathcal N(t,0,\varepsilon)=D_u\pi(u_\star)f_s(t,u_\star,\varepsilon)=0,$$
so the common equilibrium $u_\star$ becomes the equilibrium $X=0$ of the
projected system.  

The stereographic projection of the charged-particle system gives a projected equation of the form
\begin{equation}\label{projected-system}
	\frac{\dd X}{\dd t}=\mathsf B_0(x,t)X+w^{\alpha-1}(t)\mathcal{N}(t,X,\varepsilon),
\end{equation}
where $X\in\mathbb R^n$, the parameter $\varepsilon$ takes values in an open
neighbourhood $\mathcal U_\varepsilon\subseteq\mathbb R^\ell$ of the origin,
and $\dd\varepsilon/\dd t=0$.

Set
$$\mathcal G(t,X,\varepsilon):=\mathcal N(t,X,\varepsilon)-k(t)w^{1-\alpha}(t)AX.$$
Then
$$\mathcal G(t,0,\varepsilon)=0,\qquad D_{(X,\varepsilon)}\mathcal G(t,0,0)=0,$$
and the projected equation has the form
\begin{equation}\label{nonautonomous-projected-system}
	\frac{\dd X}{\dd t}=k(t)(B+A)X+w^{\alpha-1}(t)\mathcal G(t,X,\varepsilon).
\end{equation}

	\begin{lemma}\label{lem:realizability-projected-jets}
		Let $n=d-1$, and use the stereographic chart near the South Pole
		$$u_0(X)=\frac{|X|^2-1}{1+|X|^2},\qquad u_i(X)=\frac{2X_i}{1+|X|^2},\quad i=1,\ldots,n.$$
		Let $G(t,X,\varepsilon)$ be any prescribed local $C^r$ vector field on
		$\mathbb R^n$.
		Then there exists a local $C^r$ tangential force
		$f_s(t,u,\varepsilon)\in T_u\mathbb S^{d-1}$, defined near the South Pole, such that
		the nonlinear term induced by $f_s$ through the stereographic projection is exactly
		$$\mathcal N(t,X,\varepsilon)=G(t,X,\varepsilon).$$
	\end{lemma}
	
	\begin{proof}
		Define $f_s$ in the stereographic chart by
		$$f_s(t,u(X),\varepsilon)=D_Xu(X)\,G(t,X,\varepsilon).$$
		Equivalently, writing $G=(G_1,\ldots,G_n)^T$, we have
		$$(f_s)_0=\frac{4X\cdot G}{(1+|X|^2)^2},$$
		and, for $i=1,\ldots,n$,
		$$(f_s)_i=\frac{2G_i}{1+|X|^2}-\frac{4X_i(X\cdot G)}{(1+|X|^2)^2}.$$
		Since $f_s$ is obtained by differentiating the sphere parametrization $u(X)$ in the
		direction $G$, it is tangent to the sphere:
		$$f_s(t,u,\varepsilon)\cdot u=0.$$
		Thus this construction is compatible with the original decomposition of the force into
		radial and tangential components.
		
		We now compute the projected nonlinear term. From the stereographic projection formula
		derived above,
		$$\mathcal N_i=\frac{1+|X|^2}{2}\bigl(X_i(f_s)_0+(f_s)_i\bigr),\qquad i=1,\ldots,n.$$
		Substituting the preceding expressions for $(f_s)_0$ and $(f_s)_i$ gives
		$$\begin{aligned}
	\mathcal N_i&=\frac{1+|X|^2}{2}\left[X_i\frac{4X\cdot G}{(1+|X|^2)^2}+\frac{2G_i}{1+|X|^2}-\frac{4X_i(X\cdot G)}{(1+|X|^2)^2}\right]\\
	&=G_i.
\end{aligned}$$
		Therefore the induced projected term satisfies
		$\mathcal N(t,X,\varepsilon)=G(t,X,\varepsilon)$.
	\end{proof}

\subsection{Statement of chaotic behavior}

Impose in addition the stronger assumption
\begin{equation}\label{force-factorization}
	\mathcal{N}(t,X,\varepsilon)=k(t)w^{1-\alpha}(t)\bigl(AX+g_0(X,\varepsilon)\bigr),
\end{equation}
where $g_0(0,\varepsilon)=0$ and $D_Xg_0(0,0)=0$. Under this additional assumption, \eqref{projected-system} becomes
\begin{equation}\label{factorized-projected-system}
	\frac{\dd X}{\dd t}=k(t)\bigl((B+A)X+g_0(X,\varepsilon)\bigr).
\end{equation}
All explicit time dependence enters the positive scalar factor $k(t)$.  The time change
$$s(t):=\int_0^tk(\tau)\,\dd\tau$$
is a strictly increasing $C^1$ diffeomorphism from $\mathbb R$ onto
$\mathbb R$.  Indeed, \eqref{c1} implies linear growth of $s(t)$,
so both $s(t)\to+\infty$ as $t\to+\infty$ and $s(t)\to-\infty$ as
$t\to-\infty$. After this time change, equation \eqref{factorized-projected-system} becomes 
\begin{equation}
	\frac{dX}{ds}=(B+A)X+g_0(X,\varepsilon).
\end{equation}
We note that chaotic behavior arises through Shilnikov homoclinic orbits, which can be preserved by this time reparametrization.

\begin{corollary}\label{cor:three-dimensional-shilnikov}
	Assume that the conditions of
	Theorem~\ref{thm:structural-reduction} hold with $m=3$ and
	$\rho=2$, and impose the additional factorization
	\eqref{force-factorization}. Suppose also that
	$\operatorname{rank}\Omega=2$,
	so that $\Omega$ is similar to a single nilpotent Jordan block
	of size three. Then the nonlinear terms can be chosen so that
	the reduced dynamics realizes a generic three-parameter
	unfolding in which chaos emerges through the formation of a
	Shilnikov homoclinic loop on the parameterized center manifold.
\end{corollary}

\begin{corollary}
	\label{cor:persistent-skew-shilnikov}
 Assume that the
	conditions of Theorem~\ref{thm:structural-reduction} hold with
	$\rho=\infty$, and impose the factorization
	\eqref{force-factorization}. Suppose also that
	$\dim\ker F_A=3<n$,
	and choose $B$ so that $\Omega$ is nonzero and skew-symmetric. Then the 
	nonlinear terms and a $C^\infty$-flat remainder can
	be chosen so that chaos emerges through the formation of a Shilnikov
	homoclinic loop on the parameterized center manifold for
	a sequence of parameter values $(\varepsilon_{1,j},\varepsilon_{2,j})$ tending to $(0,0)$ as $j\to\infty$.
\end{corollary}

\subsection{Spectral construction of center--stable splittings}

Recall that any square matrix can be written as the sum of a symmetric matrix and a skew-symmetric matrix, that is, $$A=\frac{A+A^T}{2}+\frac{A-A^T}{2}.$$ Since the skew-symmetric part of $A$ can be absorbed into the freely chosen skew-symmetric matrix $B$, the spectral construction only depends on the symmetric part $F_A=(A+A^T)/2$.
Recall that the eigenvalues of a real symmetric matrix are all real. Moreover, the geometric multiplicity of each eigenvalue coincides with its algebraic multiplicity.

The Schur--Horn theorem gives a complete characterization of the possible diagonal entries of a real symmetric matrix with prescribed eigenvalues. In the present paper, we use its real symmetric version: for a fixed real symmetric matrix $F$, the theorem characterizes all vectors that can occur as the diagonal of $Q^TFQ$ under orthogonal changes of coordinates, see \cite{Horn1954,HornJohnson2013}.
For completeness, we recall the real symmetric version of the Schur--Horn theorem.

For two vectors $x=(x_1,\ldots,x_n)$ and $y=(y_1,\ldots,y_n)$ in $\mathbb R^n$, arranged in nonincreasing order,
$$x_1\ge\cdots\ge x_n,\qquad y_1\ge\cdots\ge y_n,$$
we say that $x$ is majorized by $y$, denoted by $x\prec y$, if $\sum_{i=1}^nx_i=\sum_{i=1}^ny_i$ and
$$\sum_{i=1}^kx_i\le\sum_{i=1}^ky_i,\qquad k=1,\ldots,n-1.$$

Denote by $F_{ij}$ the $(i,j)$ entry of matrix $F\in\mathbb{R}^{n\times n}$. Write $\operatorname{diag}F$ as the vector of diagonal entries of $F$, that is,  $\operatorname{diag}F:=(F_{11},\cdots,F_{nn})$. 

\begin{lemma}[Schur--Horn theorem]
	\label{lem:schur-horn-real}
	Let $\lambda=(\lambda_1,\ldots,\lambda_n),\,d=(d_1,\ldots,d_n)\in\mathbb R^n$ with both vectors arranged in nonincreasing order. There exists a real symmetric matrix $F\in\mathbb{R}^{n\times n}$ such that the eigenvalues of $F$ are $\lambda_1,\ldots,\lambda_n$ and $\operatorname{diag}F=d$ 
	if and only if $d$ is majorized by $\lambda$. 
\end{lemma}

Let $F\in\mathbb R^{n\times n}$ be a real symmetric matrix, and let
$\lambda_1\ge\lambda_2\ge\cdots\ge\lambda_n$
be the eigenvalues of $F$. Let $d=(d_1,\ldots,d_n)\in\mathbb R^n$, and assume that its components are arranged in nonincreasing order,
$d_1\ge d_2\ge\cdots\ge d_n$.
Then the Schur--Horn theorem guarantees that if $d$ is majorized by $\lambda=(\lambda_1,\ldots,\lambda_n)$, there exists a symmetric real matrix $F_0$ such that $F_0$ has the same eigenvalues with $F$ and prescribed diagonal entries $d_1,\cdots,d_n$. That is, there exists an orthogonal matrix $Q\in\mathbb R^{n\times n}$ such that 
$$(Q^TFQ)_{ii}=d_i,\qquad i=1,\ldots,n,$$
provided that real symmetric matrices with the same eigenvalues are orthogonally similar. 

\begin{lemma}
	\label{lem:constant-diagonalization}
	Let $F\in\mathbb R^{n\times n}$ be a real symmetric matrix. Then there exists an orthogonal matrix $Q\in\mathbb R^{n\times n}$ such that
	$$\bigl(Q^TFQ\bigr)_{11}=\bigl(Q^TFQ\bigr)_{22}=\cdots=\bigl(Q^TFQ\bigr)_{nn}=\frac{\operatorname{tr}F}{n}.$$
\end{lemma}
\begin{proof} Let the eigenvalues of $F$ be ordered as $\lambda_1\ge\lambda_2\ge\cdots\ge\lambda_n$. Set $$\tau=\frac{\operatorname{tr}F}{n}=\frac{1}{n}\sum_{i=1}^n\lambda_i$$ and take $d=(\tau,\ldots,\tau)\in\mathbb R^n$.  For each $k=1,\ldots,n-1$, since $\lambda_1,\ldots,\lambda_k$ are the largest $k$ eigenvalues, their average is not smaller than the average of all eigenvalues. Hence $d$ is majorized by $\lambda=(\lambda_1,\ldots,\lambda_n)$. By the Schur--Horn theorem, there exists an orthogonal matrix $Q\in\mathbb R^{n\times n}$ such that  $$\bigl(Q^TFQ\bigr)_{ii}=\tau=\frac{\operatorname{tr}F}{n},\qquad i=1,\ldots,n.$$  
\end{proof}

Recall that a matrix is called Hurwitz if all of its eigenvalues have negative real parts.
\begin{proposition}
	\label{prop:hurwitz-skew-perturbation}
	Let $F\in\mathbb R^{n\times n}$ be a real symmetric matrix. If
	$\operatorname{tr}F<0$,
	then there exists a real skew-symmetric matrix $B\in\mathbb R^{n\times n}$ such that $F+B$ is Hurwitz.
\end{proposition}

\begin{proof}
	Let
	$$\tau=\frac{\operatorname{tr}F}{n}<0.$$
	By Lemma~\ref{lem:constant-diagonalization}, there exists an orthogonal matrix $Q\in\mathbb R^{n\times n}$ such that
	$H:=Q^TFQ$
	satisfies
	$H_{11}=H_{22}=\cdots=H_{nn}=\tau$.
	
	Now construct a real matrix $K$ by
	$$K_{ij}=\begin{cases}
	-H_{ij},&i<j,\\
	H_{ij},&i>j,\\
	0,&i=j.
\end{cases}$$
	Since $H$ is symmetric, namely $H_{ij}=H_{ji}$, it follows directly that
	$K^T=-K$.
	Thus $K$ is skew-symmetric.
	
	Moreover, for $i<j$, we have
	$(H+K)_{ij}=H_{ij}+K_{ij}=H_{ij}-H_{ij}=0$.
	Therefore $H+K$ is a lower triangular matrix. Its diagonal entries are
	$(H+K)_{ii}=H_{ii}+K_{ii}=\tau$.
	Hence all eigenvalues of $H+K$ are equal to $\tau$. Since $\tau<0$, the matrix $H+K$ is Hurwitz.
	
	Finally, set
	$B=QKQ^T$.
	Then
	$$B^T=(QKQ^T)^T=QK^TQ^T=-QKQ^T=-B,$$
	so $B$ is skew-symmetric. Furthermore,
	$$Q^T(F+B)Q=Q^TFQ+Q^TBQ=H+K.$$
	Thus $F+B$ is orthogonally similar to $H+K$. Since $H+K$ is Hurwitz, $F+B$ is also Hurwitz. 
\end{proof}

\begin{proposition}
	\label{prop:center-stable-decomposition}
	Let $F\in\mathbb R^{n\times n}$ be a real symmetric matrix with real eigenvalues $\lambda_1\ge\lambda_2\ge\cdots\ge\lambda_n$ arranged in nonincreasing order.
	Assume that $\operatorname{tr}F=\sum_{i=1}^{n}\lambda_i<0$ and for some integer $m$ with $1\le m<n$,
	$$\sum_{i=1}^m\lambda_i\ge0,\qquad\sum_{i=1}^{m+1}\lambda_i<0.$$
	Then there exist a real skew-symmetric matrix $B\in\mathbb R^{n\times n}$ and an orthogonal matrix $P\in\mathbb R^{n\times n}$ such that
	$$P^T(F+B)P=\begin{pmatrix}
	C&0\\
	E&D
\end{pmatrix},$$
	where $C\in\mathbb{R}^{m\times m}$ is lower triangular with vanishing diagonal entries	and $D\in\mathbb{R}^{(n-m)\times(n-m)}$ is lower triangular with negative diagonal entries.
\end{proposition}

\begin{proof}
	Let
	$$\tau=\frac{\operatorname{tr}F}{n-m}.$$
	Since $\operatorname{tr}F<0$, we have
	$\tau<0$.
	Consider the vector
	$$d=(\underbrace{0,\ldots,0}_{m\ {\rm times}},\underbrace{\tau,\ldots,\tau}_{n-m\ {\rm times}}).$$
	We first verify that $d$ is majorized by
	$\lambda=(\lambda_1,\ldots,\lambda_n)$.
	
	For $1\le k\le m$, since
	$\lambda_1\ge\lambda_2\ge\cdots\ge\lambda_n$
	and
	$\sum_{i=1}^m\lambda_i\ge0$,
	the average of the first $k$ eigenvalues is not smaller than the average of the first $m$ eigenvalues. Hence
	$$\frac{1}{k}\sum_{i=1}^k\lambda_i\ge\frac{1}{m}\sum_{i=1}^m\lambda_i\ge0.$$
	Therefore
	$$\sum_{i=1}^k\lambda_i\ge0=\sum_{i=1}^kd_i,\qquad k=1,\ldots,m.$$
	
	Next let $m<k<n$ and write
	$r=k-m$.
	Then $1\le r<n-m$. Since
	$\lambda_{m+1}\ge\lambda_{m+2}\ge\cdots\ge\lambda_n$,
	the average of the first $r$ terms among
	$\lambda_{m+1},\ldots,\lambda_n$ is not smaller than the average of all these
	$n-m$ terms. Thus
	$$\sum_{i=m+1}^{m+r}\lambda_i\ge\frac{r}{n-m}\sum_{i=m+1}^{n}\lambda_i.$$
	Hence
	\begin{eqnarray*}
		\sum_{i=1}^k\lambda_i
	=
	\sum_{i=1}^m\lambda_i
	+
	\sum_{i=m+1}^{m+r}\lambda_i
	&\ge&
	\sum_{i=1}^m\lambda_i
	+
	\frac{r}{n-m}\sum_{i=m+1}^{n}\lambda_i\\
&=&	\sum_{i=1}^m\lambda_i
	+
	\frac{r}{n-m}
	\left(
	\operatorname{tr}F-\sum_{i=1}^m\lambda_i
	\right)\\
		&=&
	\frac{r}{n-m}\operatorname{tr}F
	+
	\left(1-\frac{r}{n-m}\right)
	\sum_{i=1}^m\lambda_i\\
	\end{eqnarray*}
	Since
	$\sum_{i=1}^m\lambda_i\ge0$
	and
	$1-r/(n-m)>0$,
	it follows that for $m<k<n$
	$$\sum_{i=1}^k\lambda_i\ge\frac{r}{n-m}\operatorname{tr}F=r\tau=\sum_{i=1}^{m+r}d_i=\sum_{i=1}^kd_i.$$
	
	Finally,
	$$\sum_{i=1}^nd_i=(n-m)\tau=\operatorname{tr}F=\sum_{i=1}^n\lambda_i.$$
	Therefore $d$ is majorized by $\lambda$.
	
	By the real symmetric version of the Schur--Horn theorem, there exists an orthogonal matrix $P\in\mathbb R^{n\times n}$ such that, for
	$H:=P^TFP$,
	one has
	$\operatorname{diag}(H)=d$.
	That is,
	$$H_{11}=\cdots=H_{mm}=0,\qquad H_{m+1,m+1}=\cdots=H_{nn}=\tau.$$
	
	By Proposition \ref{prop:hurwitz-skew-perturbation}, construct a real matrix $K$ by
	$$K_{ij}=\begin{cases}
	-H_{ij},&i<j,\\
	H_{ij},&i>j,\\
	0,&i=j.
\end{cases}$$
	Since $H$ is symmetric, we have
	$K^T=-K$.
	Moreover, $H+K$ is a lower triangular matrix.
	
	Set
	$B=PKP^T$.
	Then
	$B^T=(PKP^T)^T=PK^TP^T=-PKP^T=-B$,
	so $B$ is skew-symmetric. Moreover,
	$P^T(F+B)P=P^TFP+P^TBP=H+K$.
	
	Since $H+K$ is lower triangular and its diagonal entries are
	$$\underbrace{0,\ldots,0}_{m\ {\rm times}},\underbrace{\tau,\ldots,\tau}_{n-m\ {\rm times}},$$
	we can write
	$$P^T(F+B)P=\begin{pmatrix}
	C&0\\
	E&D
\end{pmatrix},$$
	where $C\in\mathbb R^{m\times m}$ is lower triangular with all diagonal entries equal to zero, and $D\in\mathbb R^{(n-m)\times(n-m)}$ is lower triangular with all diagonal entries equal to $\tau<0$.
\end{proof}

The block lower triangular form obtained in Proposition \ref{prop:center-stable-decomposition}
implies that the spectrum of $B+F$ consists of the spectrum of $C$ and the
spectrum of $D$. Since $C$ is lower triangular with zero diagonal entries, all
eigenvalues of $C$ are zero. While $D$ is lower triangular with negative diagonal entries, $D$ is Hurwitz. Hence $B+F$ has an $m$-dimensional generalized center spectral subspace associated with the eigenvalue $0$ and an $(n-m)$-dimensional stable spectral subspace associated with eigenvalues of negative real part.
These spectral subspaces give an invariant splitting
$\mathbb R^n=E^c\oplus E^s$.
With respect to this splitting, the restriction of $B+F$ to $E^c$ is nilpotent, while the restriction to $E^s$ is Hurwitz.

We remark that the necessary condition is that $F$ cannot be negative definite: a zero eigenvalue of $F$ may also produce a center block, as in Corollary \ref{cor:skew-hurwitz-decomposition-zero-eigenvalue}. Indeed, let $F\in\mathbb{R}^{n\times n}$ be symmetric negative definite and let $B\in\mathbb{R}^{n\times n}$ be any real skew-symmetric matrix, so $B^T=-B$. We claim that $M=F+B$ has all eigenvalues with strictly negative real parts. To see this, take any eigenvalue $\lambda\in\mathbb{C}$ of $M$ and let $x\in\mathbb{C}^n$ be a corresponding unit eigenvector, so $Mx=\lambda x$. Multiplying on the left by $x^*$, the conjugate transpose of $x$, gives $\lambda=x^*Mx=x^*Fx+x^*Bx$. Because $F$ is real symmetric and negative definite, the quantity $x^*Fx$ is a negative real number. On the other hand, since $B$ is real skew-symmetric, we have $\overline{x^*Bx}=x^*B^Tx=-x^*Bx$, which means that $x^*Bx$ is purely imaginary and hence has zero real part. Taking real parts in the expression for $\lambda$, we obtain $\operatorname{Re}(\lambda)=\operatorname{Re}(x^*Fx)+\operatorname{Re}(x^*Bx)=x^*Fx<0$. Since $\lambda$ was arbitrary, every eigenvalue of $M$ has negative real part, so $M$ is Hurwitz. The underlying reason is that the skew-symmetric part contributes only an imaginary term to the Rayleigh quotient and therefore cannot affect the sign of the real part, which is governed entirely by the negative definite symmetric part.

\begin{remark}
	The center block $C$ obtained in Proposition \ref{prop:center-stable-decomposition} is nilpotent, because it is lower triangular and all its diagonal entries are zero. If, in addition, the first subdiagonal entries satisfy
	$$C_{21}C_{32}\cdots C_{m,m-1}\neq0,$$
	then $\operatorname{rank}C=m-1$, and therefore $\dim\ker C=1$. Thus the zero eigenvalue of $C$ has geometric multiplicity one. Equivalently, $C$ is similar to a single nilpotent Jordan block of size $m$. In what follows, whenever a single nilpotent Jordan block is needed, we impose this nondegeneracy condition.
\end{remark}

In the case where 
$F$ has zero eigenvalues, one can choose 
$B$ such that $F+B$ contains a skew-symmetric submatrix.

\begin{corollary}
	\label{cor:skew-hurwitz-decomposition-zero-eigenvalue}
	Let $F\in\mathbb R^{n\times n}$ be a real symmetric matrix with
	$\operatorname{tr}F<0$.
	If $0$ is an eigenvalue of $F$ with multiplicity $m$, where $1\le m<n$, then there exist a real skew-symmetric matrix $B\in\mathbb R^{n\times n}$ and an orthogonal matrix $P\in\mathbb R^{n\times n}$ such that
	$$P^T(F+B)P=\begin{pmatrix}
	C&0\\
	0&D
\end{pmatrix},$$
	where $C$ is skew-symmetric and $D$ is Hurwitz. 
\end{corollary}

\begin{proof}
	Assume that $0$ is an eigenvalue of $F$ with multiplicity $m$. Since $F$
	is real symmetric, there exists an orthogonal matrix
	$P\in\mathbb R^{n\times n}$ such that
	$$P^TFP=\begin{pmatrix}
	0_m&0\\
	0&F_1
\end{pmatrix},$$
	where $F_1\in\mathbb R^{(n-m)\times(n-m)}$ is real symmetric. Moreover,
	$\operatorname{tr}F_1=\operatorname{tr}F<0$.
	By Proposition \ref{prop:hurwitz-skew-perturbation}, there exists a real
	skew-symmetric matrix $K_1\in\mathbb R^{(n-m)\times(n-m)}$ such that
	$F_1+K_1$ is Hurwitz.
	
	Let $C\in\mathbb R^{m\times m}$ be any real skew-symmetric matrix and set
	$$K=\begin{pmatrix}
	C&0\\
	0&K_1
\end{pmatrix}.$$
	Then $K$ is skew-symmetric. Define
	$B=PKP^T$.
	It follows that $B^T=-B$, and
	$$P^T(F+B)P=P^TFP+K=\begin{pmatrix}
	C&0\\
	0&F_1+K_1
\end{pmatrix}.$$
	Setting $D=F_1+K_1$, we obtain the desired block decomposition, with $C$
	skew-symmetric and $D$ Hurwitz.
\end{proof}

\begin{example}
	Consider
	$$F=\begin{pmatrix}
	F_{11}&0&0\\
	0&F_{22}&0\\
	0&0&F_{33}
\end{pmatrix},$$
	where $F_{11}<0$, $F_{22}<0$, $F_{33}>0$, and $\operatorname{tr}F=F_{11}+F_{22}+F_{33}<0$. Thus $F$ has exactly one positive eigenvalue and two negative eigenvalues, with negative trace.
	
	Take the skew-symmetric matrix
	$$B=\begin{pmatrix}
	0&b&b\\
	-b&0&b\\
	-b&-b&0
\end{pmatrix},$$
	where $b\in\mathbb{R}$ is a parameter to be chosen. Then
	$$M=F+B=\begin{pmatrix}
	F_{11}&b&b\\
	-b&F_{22}&b\\
	-b&-b&F_{33}
\end{pmatrix}.$$
	
	A direct computation gives the determinant
	$$\det(M)=F_{11}F_{22}F_{33}+(F_{11}+F_{22}+F_{33})b^2.$$
	Setting $\det(M)=0$ yields
	$$b^2=-\frac{F_{11}F_{22}F_{33}}{F_{11}+F_{22}+F_{33}}.$$
	Since $F_{11}<0$, $F_{22}<0$, $F_{33}>0$, we have $F_{11}F_{22}F_{33}>0$. Since $\operatorname{tr}F=F_{11}+F_{22}+F_{33}<0$, the right-hand side is strictly positive. Hence we may choose a real value of $b$, for instance the positive square root, such that $\det(M)=0$. Consequently, $M$ has a zero eigenvalue.

\end{example}

\subsection{Exponential trichotomy}
As discussed in Proposition \ref{prop:center-stable-decomposition} and Corollary \ref{cor:skew-hurwitz-decomposition-zero-eigenvalue}, we have found some conditions under which the restriction of the matrix $B+F$ to an invariant subspace is skew-symmetric or nilpotent, and is Hurwitz on the complementary subspace. Let us now assume that $B$ and $F$ are such that $$B+F:\mathbb{R}^n\to\mathbb{R}^n$$ has an invariant center subspace $E^c$ of dimension $m$ and an invariant stable subspace $E^s$. 

Denote by $\Phi(t,s)$ the evolution operator of linear system $\dd x/\dd t=k(t)(B+F)x$. The linear system is said to admit an exponential trichotomy if there exist projections
$P^c,P^s,P^u$,
with
$$P^c+P^s+P^u=I,$$
such that these projections are invariant under $\Phi$, that is,
$$\Phi(t,s)P^\ast=P^\ast\Phi(t,s),\qquad\ast\in\{c,s,u\},$$
and the following estimates hold: there exists $\eta,C_\eta>0$ such that
$$\|\Phi(t,s)P^c\|\le C_\eta e^{\eta|t-s|},\qquad t,s\in\mathbb{R}.$$
There exist $p,K>0$ such that
$$\|\Phi(t,s)P^s\|\le Ke^{-p(t-s)},\qquad t\ge s,$$
and
$$\|\Phi(t,s)P^u\|\le Ke^{p(t-s)},\qquad t\le s.$$
The matrix norm $\|\cdot\|$ is the operator norm induced by the Euclidean norm on $\mathbb{R}^n$.

\begin{lemma}\label{lem1}
	Let $S\in\mathbb{R}^{n\times n}$ be a Hurwitz matrix. Then there exist constants $M\ge1$ and $\omega>0$ such that for all $\xi\ge0$,
	$$\|e^{\xi S}\|\le Me^{-\omega\xi}.$$
\end{lemma}

\begin{proof}
	Since $S$ is Hurwitz, we begin by taking its Jordan normal form. There exists an invertible matrix $P$ such that $S=PJP^{-1}$, where $J=\operatorname{diag}(J_1,J_2,\dots,J_k)$ and each Jordan block is given by $J_i=\lambda_iI_{m_i}+N_i$, where $N_i$ is nilpotent with $N_i^{m_i}=0$, and $m_i$ is the size of the block.
	
	It holds that $e^{\xi S}=Pe^{\xi J}P^{-1}$, and $e^{\xi J}=\operatorname{diag}\bigl(e^{\xi J_1},e^{\xi J_2},\dots,e^{\xi J_k}\bigr)$.
	
	For a single Jordan block, nilpotency of $N_i$ yields
	$$e^{\xi J_i}=e^{\lambda_i\xi}e^{\xi N_i}=e^{\lambda_i\xi}\sum_{j=0}^{m_i-1}\frac{\xi^j}{j!}N_i^{\,j}.$$
	Taking norms gives
	$$\|e^{\xi J_i}\|\le e^{\operatorname{Re}(\lambda_i)\xi}\sum_{j=0}^{m_i-1}\frac{\xi^j}{j!}\|N_i\|^{\,j}.$$
	The sum on the right is a polynomial in $\xi$ of degree at most $m_i-1$.
	
	Let $m=\max\{m_1,\dots,m_k\}$ and define
	$p=\min_i\bigl(-\operatorname{Re}(\lambda_i)\bigr)>0$.
	Because $e^{\operatorname{Re}(\lambda_i)\xi}\le e^{-p\xi}$ for every $i$, there exists a constant $C_1>0$ such that for all $\xi\ge0$,
	\begin{equation}\label{c2}
	\|e^{\xi J}\|=\max_i\|e^{\xi J_i}\|\le C_1(1+\xi^{m-1})e^{-p\xi}.
\end{equation}
	
	Now choose any $\omega$ satisfying $0<\omega<p$. Rewrite the estimate \eqref{c2} as
	$$\|e^{\xi J}\|\le C_1\bigl[(1+\xi^{m-1})e^{-(p-\omega)\xi}\bigr]e^{-\omega\xi}.$$
	The function $\phi(\xi)=(1+\xi^{m-1})e^{-(p-\omega)\xi}$ is continuous on $[0,\infty)$ and tends to $0$ as $\xi\to\infty$ (the exponential decay dominates any polynomial growth). Hence it attains a finite maximum. Set
	$C_2=\sup_{\xi\ge0}\phi(\xi)<\infty$.
	Then we obtain
	$\|e^{\xi J}\|\le C_1C_2e^{-\omega\xi}$.
	
	Returning to the original matrix $S$,  $e^{\xi S}=Pe^{\xi J}P^{-1}$ yields
	$$\|e^{\xi S}\|\le\|P\|\|P^{-1}\|\cdot\|e^{\xi J}\|\le\|P\|\|P^{-1}\|\,C_1C_2e^{-\omega\xi}.$$
	Let $M=\|P\|\|P^{-1}\|C_1C_2$. Then for every $\xi\ge0$,
	$\|e^{\xi S}\|\le Me^{-\omega\xi}$.
\end{proof}

\begin{lemma}\label{lem4}
	Let $\Omega\in\mathbb{R}^{n\times n}$ be nilpotent. Then there exists a constant $C>0$ depending only on $\Omega$ such that for all $t\in\mathbb{R}$,
	$\|e^{t\Omega}\|\le C(1+|t|)^{n-1}$.
	In particular, the norm of $e^{t\Omega}$ grows at most polynomially in $|t|$.
\end{lemma}

\begin{proof}
	Since $\Omega$ is nilpotent, there exists a positive integer $k\le n$ such that $\Omega^k=0$. The matrix exponential is defined by the power series
	$$e^{t\Omega}=\sum_{j=0}^{\infty}\frac{t^j}{j!}\Omega^j=\sum_{j=0}^{k-1}\frac{t^j}{j!}\Omega^j.$$
 Applying the norm to the finite sum and using the triangle inequality, we obtain
	$$\|e^{t\Omega}\|=\left\|\sum_{j=0}^{k-1}\frac{t^j}{j!}\Omega^j\right\|\le\sum_{j=0}^{k-1}\frac{|t|^j}{j!}\|\Omega^j\|.$$
	By submultiplicativity, $\|\Omega^j\|\le\|\Omega\|^j$ for each $j$. Hence
	$$\|e^{t\Omega}\|\le\sum_{j=0}^{k-1}\frac{|t|^j}{j!}\|\Omega\|^j.$$
	Define $M=\max\{1,\|\Omega\|\}$. Then for every $j=0,1,\dots,k-1$, we have $\|\Omega\|^j\le M^j\le M^{k-1}$. Therefore,
	$$\|e^{t\Omega}\|\le M^{k-1}\sum_{j=0}^{k-1}\frac{|t|^j}{j!}\le C(1+|t|)^{n-1}.$$
\end{proof}
	
	\begin{lemma}\label{lem5}
		Let $\Omega\in\mathbb{R}^{n\times n}$ be skew-symmetric  and let $t\in\mathbb{R}$. Then the matrix exponential $e^{t\Omega}$ is orthogonal and its norm satisfies $\|e^{t\Omega}\|=1$.
	\end{lemma}
	
	\begin{proof}
		Because the matrix exponential is defined by the convergent power series $$e^{t\Omega}=\sum_{k=0}^{\infty}\frac{t^{k}\Omega^{k}}{k!},$$ taking the transpose and using the skew-symmetry $\Omega^T=-\Omega$ yields $$(e^{t\Omega})^T=\sum_{k=0}^{\infty}\frac{t^{k}(\Omega^T)^{k}}{k!}=\sum_{k=0}^{\infty}\frac{t^{k}(-\Omega)^{k}}{k!}=e^{-t\Omega}.$$ Since $e^{t\Omega}e^{-t\Omega}=e^{t\Omega+(-t\Omega)}=e^{0}=I$, we have $e^{-t\Omega}=(e^{t\Omega})^{-1}$ and therefore $(e^{t\Omega})^T=(e^{t\Omega})^{-1}$, which means $e^{t\Omega}$ is an orthogonal matrix. An orthogonal matrix preserves the Euclidean norm of every vector, i.e. $\|e^{t\Omega}x\|=\|x\|$ for all $x\in\mathbb{R}^{n}$. Consequently, the norm  satisfies $$\|e^{t\Omega}\|=\max_{\|x\|=1}\|e^{t\Omega}x\|=\max_{\|x\|=1}\|x\|=1.$$
	\end{proof}

	Under the assumption \eqref{c1}, the linear equation
	$$\frac{\dd x}{\dd t}=k(t)(B+F)x$$
	admits an exponential trichotomy on $\mathbb{R}$ with center space $E^c$, stable space $E^s$, and trivial unstable space
	$E^u=\{0\}$. Indeed,	since the coefficient matrix is of the form $k(t)(B+F)$, where $B+F$ is constant, the matrices at different times commute. Therefore the evolution operator is given by
	$$\Phi(t,s)=\exp\!\left(\int_s^tk(\tau)\,\dd\tau\,(B+F)\right).$$
	Set
	$$K(t,s):=\int_s^tk(\tau)\,\dd\tau.$$
	Then
	$$\Phi(t,s)=e^{K(t,s)(B+F)}.$$
	
	Because $E^c$ and $E^s$ are invariant under $B+F$, they are also invariant under $\Phi(t,s)$ for all $t,s\in\mathbb R$. Let
	$$P^c,P^s:\mathbb R^n\to\mathbb R^n$$
	be the projections onto $E^c$ and $E^s$, respectively, associated with the decomposition
	$\mathbb R^n=E^c\oplus E^s$. Set $P^u:=0$.
	Then $P^c+P^s+P^u=I$,
	and, by invariance of the splitting,
	$$\Phi(t,s)P^\ast=P^\ast\Phi(t,s),\qquad\ast\in\{c,s,u\}.$$
	
	Let
	$$\Omega:=(B+F)|_{E^c}:E^c\to E^c,\qquad S:=(B+F)|_{E^s}:E^s\to E^s$$
	be the restrictions of $B+F$ to the two invariant subspaces. Since $B+F$ is constant and the splitting is invariant, the restrictions of the evolution operator to $E^c$ and $E^s$ are given by
	$$\Phi(t,s)|_{E^c}=e^{K(t,s)\Omega},\qquad\Phi(t,s)|_{E^s}=e^{K(t,s)S},$$
	respectively. Hence, for every $x\in\mathbb R^n$,
	$$\Phi(t,s)P^cx=e^{K(t,s)\Omega}P^cx,\qquad\Phi(t,s)P^sx=e^{K(t,s)S}P^sx,$$
	given that $P^cx\in E^c$ and $P^sx\in E^s$. Therefore,
	$$\|\Phi(t,s)P^c\|\le\|e^{K(t,s)\Omega}\|\,\|P^c\|,\qquad\|\Phi(t,s)P^s\|\le\|e^{K(t,s)S}\|\,\|P^s\|.$$

	We first estimate the center part. In the case where $\Omega$ is skew-symmetric,
	by Lemma \ref{lem5}, $e^{K(t,s)\Omega}$ is orthogonal with norm 1.
	Consequently,
	$$\|e^{K(t,s)\Omega}\|\le C_\eta e^{\eta|t-s|},\qquad t,s\in\mathbb{R}.$$
	This yields
	$$\|\Phi(t,s)P^c\|\le\|P^c\|C_\eta e^{\eta|t-s|},\qquad t,s\in\mathbb{R}.$$
	When $\Omega$ is nilpotent, by Lemma \ref{lem4} 
	$$\|e^{K(t,s)\Omega}\|\le C(1+|K(t,s)|)^{n-1}\le C(1+\overline{k}|t-s|)^{n-1}.$$
	Since any polynomial can be dominated by an exponential function, there exist $\eta$ and $C_\eta$ such that 	$$\|e^{K(t,s)\Omega}\|\le C_\eta e^{\eta|t-s|},\qquad t,s\in\mathbb{R}.$$  
	Here $\eta$ can be chosen arbitrarily small in both cases.
	
	We now estimate the stable part.
	Since $S$ is Hurwitz, Lemma \ref{lem1} yields there exist constants $M\ge1$ and $\omega>0$ such that 
	$$\|e^{\xi S}\|\le Me^{-\omega\xi},\qquad\xi\ge0.$$
	For $t\ge s$, Assumption \eqref{c1} gives
	$$K(t,s)=\int_s^tk(\tau)\,\dd\tau\ge\underline{k}(t-s).$$
	Therefore,
	$$\|\Phi(t,s)P^s\|\le\|P^s\|Me^{-\omega K(t,s)}\le\|P^s\|Me^{-\omega\underline{k}(t-s)},\qquad t\ge s.$$

	Finally, the unstable projection is zero by construction, so the unstable subspace is trivial. Hence the equation admits an exponential trichotomy with center space $E^c$, stable space $E^s$, and $E^u=\{0\}$.

\subsection{Parameterized center-manifold reduction}

	Consider now the projected system after the invariant splitting
	$\mathbb R^n=E^c\oplus E^s$.
	Write $x_c\in E^c$, $\varepsilon\in\mathcal U_\varepsilon$, and $v\in E^s$. After a change of coordinates, together with the frozen-parameter equation $\dd\varepsilon/\dd t=0$, the system reads
	\begin{eqnarray*}
		\frac{\dd x_c}{\dd t} &=& k(t)\Omega x_c+P(t,x_c,\varepsilon,v),\\
		\frac{\dd\varepsilon}{\dd t}&=&0,\\
		\frac{\dd v}{\dd t} &= &k(t)Sv+Q(t,x_c,\varepsilon,v).
	\end{eqnarray*}
	where $\Omega:E^c\to E^c$ is the dynamical center matrix, $S:E^s\to E^s$ is Hurwitz, $P$ takes values in $E^c$, and $Q$ takes values in $E^s$. The higher order terms $P,Q$ satisfy that for all $t\in\mathbb R$
	$$P(t,0,\varepsilon,0)=Q(t,0,\varepsilon,0)=0,\qquad D_{(x_c,\varepsilon,v)}P(t,0,0,0)=D_{(x_c,\varepsilon,v)}Q(t,0,0,0)=0.$$ Moreover, $P,Q$ are $C^1$ in $t\in\mathbb{R}$ and $C^{\rho+1}$ in $(x_c,\varepsilon,v)\in E^c\times\mathcal U_\varepsilon\times E^s$, with the derivatives used below uniformly bounded on a fixed neighbourhood. 

	We apply the smooth nonautonomous center-manifold theorem to the splitting
	$$E^s\oplus\bigl(E^c\times\mathbb R^\ell\bigr),$$
	where $E^s$ is the exponentially stable block and
	$E^c\times\mathbb R^\ell$ is the extended center block of the
	frozen-parameter extension.  Since $E^u=\{0\}$, the resulting graph
	over the extended center block is precisely the nonautonomous center
	manifold required here.  The preceding trichotomy and
	$(\rho+1)\eta<\omega\underline k$ give the required order-$\rho$ gap for
	the center-manifold construction; see \cite{Aulbach1982}.
	Consequently, there exist $r_0>0$ and a local time-dependent center
	manifold of the graph form
	$$\mathcal W^c_{\mathrm{loc}}=\{(t,x_c,\varepsilon,v)\in\mathbb{R}\times E^c\times\mathcal U_\varepsilon\times E^s:\ v=h(t,x_c,\varepsilon),\ |(x_c,\varepsilon)|<r_0\},$$
	where the graph function
	$$h:\mathbb R\times\{(x_c,\varepsilon)\in E^c\times\mathcal U_\varepsilon:|(x_c,\varepsilon)|<r_0\}\to E^s$$
	is $C^\rho$ in $(x_c,\varepsilon)$, is of class $C^1$ in $t$, and has
	uniformly bounded derivatives.
	Moreover, $h(t,x_c,\varepsilon)$ satisfies
	$$h(t,0,\varepsilon)=0,\qquad D_{x_c}h(t,0,0)=0$$
	for all $t\in\mathbb{R}$.
	The dynamics on $\mathcal W^c_{\mathrm{loc}}$ is governed by the parameterized reduced equation
	$$\frac{\dd x_c}{\dd t}=R(t,x_c,\varepsilon):=k(t)\Omega x_c+P(t,x_c,\varepsilon,h(t,x_c,\varepsilon)),\qquad\frac{\dd\varepsilon}{\dd t}=0.$$
	The function $h$ satisfies the invariance equation
	\begin{equation}\label{c3}
	\frac{\partial h}{\partial t}(t,x_c,\varepsilon)+D_{x_c}h(t,x_c,\varepsilon)R(t,x_c,\varepsilon)=k(t)S\,h(t,x_c,\varepsilon)+Q(t,x_c,\varepsilon,h(t,x_c,\varepsilon)).
\end{equation}

	In particular, the local dynamics of the full system near the origin is determined by
	the reduced vector field $R$ on the center manifold. Furthermore, nearby solutions admit
	an asymptotic phase relative to $\mathcal W^c_{\mathrm{loc}}$ and approach appropriate
	solutions on the center manifold at an exponential rate. If the origin of the reduced
	equation is stable, then the full system and the reduced system have the same local
	stability type.
	
	We now construct the polynomial coefficients recursively.  Set
	$\sigma_1=0$.  Suppose that $h_2,\ldots,h_{r-1}$ have already been
	constructed, and define
	$$\sigma_{r-1}(t,u):=\sum_{j=2}^{r-1}h_j(t,u),\qquad u=(x_c,\varepsilon).$$
	Let $H_r(t,u)$ be the degree-$r$ homogeneous part of
	$$Q\bigl(t,u,\sigma_{r-1}(t,u)\bigr)-D_{x_c}\sigma_{r-1}(t,u)P\bigl(t,u,\sigma_{r-1}(t,u)\bigr).$$
	Matching the degree-$r$ terms in the formal invariance equation gives
	\begin{equation}\label{c4}
	\frac{\partial h_r}{\partial t}(t,x_c,\varepsilon)+k(t)D_{x_c}h_r(t,x_c,\varepsilon)\Omega x_c-k(t)Sh_r(t,x_c,\varepsilon)=H_r(t,x_c,\varepsilon).
\end{equation}
	 To show that the differential equation \eqref{c4} has a solution, we start by introducing a useful lemma, aiming to separate variables.

\begin{lemma}\label{lem2}
	Let $k$ be a positive integer and $f:\mathbb{R}^n\times\mathbb{R}^m\to\mathbb{R}$ satisfy
$f(\lambda x,y)=\lambda^kf(x,y)$ for all $\lambda>0$, $x\in\mathbb{R}^n$, $y\in\mathbb{R}^m$. For each fixed $y$, $f(x,y)$ is of class $C^{k+1}(\mathbb{R}^n)$ as a function of $x$. For each fixed $x$, the map $f(x,y)$ is of class $C^1(\mathbb{R}^m)$ with respect to $y$.
	Then $f$ admits a unique representation
	\begin{equation}\label{c5}
	f(x,y)=\sum_{|p|=k}c_p(y)\,x^p,
\end{equation}
	where $p=(p_1,\dots,p_n)$ is a multi-index, $x^p=x_1^{p_1}\cdots x_n^{p_n}$, and each coefficient
	\begin{equation}\label{c6}
	c_p(y)=\frac{1}{p!}\frac{\partial^{|p|}f}{\partial x^p}(0,y)
\end{equation}
	is a $C^1$ function of $y\in\mathbb{R}^m$.
\end{lemma}

\begin{proof}
	Fix an arbitrary $y\in\mathbb{R}^m$ first and set $g(x):=f(x,y)$. By  assumption, we have $g\in C^{k+1}(\mathbb{R}^n)$, and that
	$$g(\lambda x)=\lambda^kg(x)\qquad\forall\,\lambda>0,\;\forall x\in\mathbb{R}^n.$$
	For any multi-index $p$ with $|p|\le k+1$, differentiate both sides with respect to $x$ using the chain rule to obtain
	$$(\partial^pg)(\lambda x)\cdot\lambda^{|p|}=\lambda^k\partial^pg(x),$$
	or equivalently $(\partial^pg)(\lambda x)=\lambda^{k-|p|}\,\partial^pg(x)$ for all $\lambda>0$ and $x$.  Now let $\lambda\to0^+$.  Because $g$ is of class $C^{k+1}$, each derivative $\partial^pg$ is continuous and the left-hand side tends to $\partial^pg(0)$.  The limit of the right-hand side is determined by the exponent $k-|p|$ (recall $k$ is a positive integer).  If $|p|<k$, then $k-|p|\ge1$ and $\lambda^{k-|p|}\to0$, forcing $\partial^pg(0)=0$.  If $|p|=k$, then $k-|p|=0$ and the relation becomes $(\partial^pg)(\lambda x)=\partial^pg(x)$. Letting $\lambda\to0^+$ shows that $\partial^pg$ is constant on $\mathbb{R}^n$, equal to $\partial^pg(0)$.  If $|p|=k+1$, then $k-|p|=-1$ and we obtain $(\partial^pg)(\lambda x)=\lambda^{-1}\partial^pg(x)$.  For a fixed $x\neq0$, the left side approaches the finite value $\partial^pg(0)$ as $\lambda\to0^+$, so the right side can remain bounded only if $\partial^pg(x)=0$.  Hence every derivative of order $k+1$ vanishes identically on $\mathbb{R}^n$.
	
	By the multivariate Taylor theorem with integral remainder,
	$$g(x)=\sum_{|p|\le k}\frac{\partial^pg(0)}{p!}x^p+R_{k+1}(x),$$
	where
	$$R_{k+1}(x)=(k+1)\sum_{|p|=k+1}\frac{x^p}{p!}\int_0^1(1-\tau)^k\partial^pg(\tau x)\,\mathrm d\tau.$$
	As shown above, every derivative $\partial^pg$ with
	$|p|=k+1$ vanishes identically. Hence
	$R_{k+1}(x)=0$. Moreover,
	$\partial^pg(0)=0$ whenever $|p|<k$. Therefore
	$$g(x)=\sum_{|p|=k}\frac{\partial^pg(0)}{p!}x^p.$$
	
	It remains to show that each $c_p(y)$ is a $C^1$ function of $y$.  The monomials $\{x^p\}_{|p|=k}$ are linearly independent, so we may select finitely many points $x^{(1)},\dots,x^{(N)}\in\mathbb{R}^n$  such that the matrix $\bigl[(x^{(i)})^p\bigr]_{i,p}$ is invertible.  Evaluating \eqref{c5} at these points gives
	$$f(x^{(i)},y)=\sum_{|p|=k}c_p(y)\,(x^{(i)})^p,\qquad i=1,\dots,N.$$
	By invertibility, each $c_p(y)$ is a constant linear combination of the functions $f(x^{(i)},y)$.  Because every $f(x^{(i)},\cdot)$ is $C^1$ by the stated $C^1$ dependence on $y$, and the class $C^1$ is closed under linear combinations, we conclude $c_p\in C^1(\mathbb{R}^m)$.
\end{proof}

For finite-dimensional real vector spaces $X$ and $Y$, and for $r\in\mathbb N$, we denote by
$\mathcal P_r(X,Y)$
the vector space of $Y$-valued homogeneous polynomials of degree $r$ on $X$. Equivalently, set
$n=\dim X$, and fix a basis $\{x_1,\cdots,x_n\}$ of $X$. Then
$h\in\mathcal P_r(X,Y)$ if and only if 
$h$ can be written uniquely in the form
$$h(x)=\sum_{|p|=r}x^pa_p,\qquad a_p\in Y,$$
where $p=(p_1,\dots,p_n)$ is a multi-index, $x^p=x_1^{p_1}\cdots x_n^{p_n}$.
In particular,
$$\dim\mathcal P_r(X,Y)=\dim(Y)\binom{\dim(X)+r-1}{r}.$$
For the coefficient calculations below, set
$$\mathcal U:=E^c\times\mathbb R^\ell,\qquad u=(x_c,\varepsilon)\in\mathcal U,$$
and define the two linear maps
$$\Omega_{\rm ext}(x_c,\varepsilon):=(\Omega x_c,0)\in\mathcal U,\qquad\Omega_c(x_c,\varepsilon):=\Omega x_c\in E^c.$$
Thus $\Omega_{\rm ext}$ is the frozen-parameter extension used when a
polynomial on $\mathcal U$ is differentiated, whereas $\Omega_c$ is the
$E^c$-valued linear part of the reduced vector field.  Thanks to
Lemma~\ref{lem2}, the terms $h_r(t,u)$ and $H_r(t,u)$ are homogeneous
polynomials.  We first reduce the homological equation to a finite-dimensional
coefficient equation.

\begin{lemma}\label{lem3}
	For fixed $r$, equation \eqref{c4} for a degree-$r$ homogeneous polynomial
	$h_r(t,u)$ is equivalent to the coefficient equation \eqref{c8} below.
\end{lemma}
\begin{proof}
	
	Let $q=\dim E^c+\ell$, the dimension of the extended variable $u=(x_c,\varepsilon)$, and let
$$\nu_r:=\dim\mathcal P_r(E^c\times\mathbb R^\ell,\mathbb R)=\binom{q+r-1}{r}.$$
Choose the monomial basis
$$\{u^p:\ |p|=r\}$$
of $\mathcal P_r(E^c\times\mathbb R^\ell,\mathbb R)$, where
$$u^p=u_1^{p_1}\cdots u_q^{p_q},\qquad p=(p_1,\dots,p_q)\in\mathbb N^q,\qquad|p|=p_1+\cdots+p_q=r.$$
Let $d_s:=\dim E^s$ and let $\{v_1,\dots,v_{d_s}\}$ be a basis of $E^s$. Then
$$\{u^pv_j:\ |p|=r,\ j=1,\dots,d_s\}$$
is a basis of $\mathcal P_r(E^c\times\mathbb R^\ell,E^s)$.

By Lemma \ref{lem2},
$$h_r(t,\cdot)\in\mathcal P_r(E^c\times\mathbb R^\ell,E^s),\qquad H_r(t,\cdot)\in\mathcal P_r(E^c\times\mathbb R^\ell,E^s)$$
can be written uniquely in the form
\begin{eqnarray}\label{c7}
h_r(t,u)&=&\sum_{j=1}^{d_s}\sum_{|p|=r}u^p v_j\eta_{j,p}(t),\\
H_r(t,u)&=&\sum_{j=1}^{d_s}\sum_{|p|=r}u^p v_j\zeta_{j,p}(t),\nonumber
\end{eqnarray}
where every $\eta_{j,p}(t)$ and $\zeta_{j,p}(t)$ is continuous. After substituting \eqref{c7} into \eqref{c4}, since $\{u^pv_j\}$ is a basis in $\mathcal P_r(E^c\times\mathbb R^\ell,E^s)$, the corresponding coefficients are equal. Collecting the coefficients $\eta_{j,p}(t)$ and $\zeta_{j,p}(t)$ into
column vectors
$$c(t),d(t)\in\mathbb R^{M_r},\qquad M_r=d_s\nu_r.$$
After fixing the above bases, we use the same symbol $S$ for the restriction of
the stable linear part to $E^s$ and for its representing matrix. 

We compute the three terms on the left-hand side of \eqref{c4}. First,
$$\frac{\partial h_r}{\partial t}(t,u)=\sum_{j=1}^{d_s}\sum_{|p|=r}u^pv_j\frac{\dd\eta_{j,p}(t)}{\dd t},$$
which corresponds to $\dd c(t)/\dd t$. Next,
$$Sh_r(t,u)=\sum_{j=1}^{d_s}\sum_{|p|=r}u^p\,S(v_j)\eta_{j,p}(t),$$
so on the coefficient space this term is represented by
$$(I_{\nu_r}\otimes S)c(t),$$
where $\otimes$ stands for the Kronecker product. Finally, define the linear operator
$$\mathcal N_r:\mathcal P_r(E^c\times\mathbb R^\ell,\mathbb R)\to\mathcal P_r(E^c\times\mathbb R^\ell,\mathbb R),\qquad\mathcal N_r\varphi:=D_u\varphi(u)\,\Omega_{\rm ext}u.$$
Let $\mathbf N_r$ be the matrix representing $\mathcal N_r$ in the monomial
basis $\{u^p:|p|=r\}$. Then
$$D_uh_r(t,u)(\Omega_{\rm ext}u)=\sum_{j=1}^{d_s}\sum_{|p|=r}\mathcal N_r(u^p)v_j\eta_{j,p}(t),$$
which means that on the coefficient space this term is represented by
$$(\mathbf N_r\otimes I_{d_s})c(t).$$

Therefore the equation \eqref{c4} is equivalent to
$$\frac{\dd c(t)}{\dd t}+k(t)(\mathbf N_r\otimes I_{d_s})c(t)-k(t)(I_{\nu_r}\otimes S)c(t)=d(t),$$
that is,
\begin{equation}\label{c8}
	\frac{\dd c(t)}{\dd t}=k(t)B_r\,c(t)+d(t),
\end{equation}
where
$$B_r:=(I_{\nu_r}\otimes S)-(\mathbf N_r\otimes I_{d_s}).$$
The scalar coefficient $k(t)$ is inherited directly from the linear part of the original
equation. Moreover, the term $d(t)$ in equation \eqref{c8} is obtained from
the $r$-th homogeneous term $H_r$, whose coefficients are determined recursively by the
Taylor coefficients of the original nonlinearities $P$ and $Q$. The coefficients of the
original equation are of class $C^1$ with respect to $t$. Then $k$ is continuous, and all
Taylor coefficients of $P$ and $Q$ are continuous functions of $t$. By induction, the same
is therefore true for the coefficients of $H_r$, and hence for the vector $d(t)$.

Consequently, $k(t)B_r$ and $d(t)$ are continuous on $\mathbb R$, and the
coefficient comparison proves the asserted equivalence.
\end{proof}

After fixing the monomial basis
$\{u^p:|p|=r\}$ and the basis $\{v_1,\ldots,v_{d_s}\}$ of $E^s$,
we use the same symbol $S$ for the matrix representing the restriction of
the stable linear part to $E^s$.

\begin{lemma}\label{lem:homological-matrix-hurwitz}
	 The matrix
	$$B_r=I_{\nu_r}\otimes S-\mathbf N_r\otimes I_{d_s}$$
	is Hurwitz when $\Omega$ is nilpotent or skew-symmetric.
\end{lemma}
\begin{proof}
	
	First, we record a spectral observation which is independent of these two
	cases. Over $\mathbb C$, choose bases in which $S$ and $\mathbf N_r$ are in
	Schur triangular form. Thus
	$$P^{-1}SP=T_S,\qquad Q^{-1}\mathbf N_rQ=T_N,$$
	where $T_S$ and $T_N$ are upper triangular matrices whose diagonal entries are
	the eigenvalues of $S$ and $\mathbf N_r$, respectively, counted with algebraic
	multiplicity. We have
	$$(Q\otimes P)^{-1}B_r(Q\otimes P)=I_{\nu_r}\otimes T_S-T_N\otimes I_{d_s}.$$
	The two matrices $I_{\nu_r}\otimes T_S$ and $T_N\otimes I_{d_s}$ are upper
	triangular, and hence so is their difference. Its diagonal entries are precisely
	$$\lambda-\mu,\qquad\lambda\in\sigma(S),\quad\mu\in\sigma(\mathbf N_r),$$
	where $\sigma(A)$ represents the spectrum of the matrix $A$. Therefore
	$$\sigma(B_r)=\{\lambda-\mu:\lambda\in\sigma(S),\ \mu\in\sigma(\mathbf N_r)\}.$$
	
	We now determine the possible real parts of the eigenvalues of $\mathbf N_r$.
	Recall that for a homogeneous polynomial $\varphi$ in $u$, the operator
	$\mathcal N_r$ is defined by
	$$\mathcal N_r\varphi(u)=D_u\varphi(u)\Omega_{\rm ext}u.$$
	Define the pull-back action
	$$(\mathcal T_t\varphi)(u):=\varphi(e^{t\Omega_{\rm ext}}u),$$
	then $\mathcal T_t$ is a one-parameter group on
	$\mathcal P_r(E^c\times\mathbb R^\ell,\mathbb R)$. Indeed,
	$\mathcal T_{t+s}\varphi(u)=\varphi(e^{(t+s)\Omega_{\rm ext}}u)=\varphi(e^{t\Omega_{\rm ext}}e^{s\Omega_{\rm ext}}u)=\mathcal T_s(\mathcal T_t\varphi)(u)$, and hence
	$\mathcal T_{t+s}=\mathcal T_s\mathcal T_t$. Differentiating this action at $t=0$ gives
	$$\left.\frac{\dd}{\dd t}\right|_{t=0}\mathcal T_t\varphi=D_u\varphi(u)\Omega_{\rm ext}u=\mathcal N_r\varphi.$$
	Thus $\mathcal N_r$ is the infinitesimal generator of the one-parameter
	linear group $\mathcal T_t$. After passing to the monomial basis, this
	generator is represented by the matrix $\mathbf N_r$, and the action
	$\mathcal T_t$ is represented by the matrix $T_t$. Since the polynomial space
	is finite dimensional, the matrix-valued linear initial value problem
	$$\frac{\dd T_t}{\dd t}=\mathbf N_rT_t,\qquad T_0=I,$$
	has the unique solution
	$T_t=e^{t\mathbf N_r}$.
	
	If $\Omega$ is nilpotent, then $\Omega_{\rm ext}$ is nilpotent and
	$e^{t\Omega_{\rm ext}}$ is a polynomial in $t$. Since
	$\varphi$ is a polynomial, $\mathcal T_t\varphi=\varphi(e^{t\Omega_{\rm ext}}u)$ is also a
	polynomial in $t$. Therefore the matrix $e^{t\mathbf N_r}$ is polynomial in
	$t$, which implies that $\mathbf N_r$ is nilpotent. In particular,
	$\sigma(\mathbf N_r)=\{0\}$.
	
	If $\Omega$ is skew-symmetric, then $\Omega_{\rm ext}$ is skew-symmetric on
	$\mathcal U$ and $e^{t\Omega_{\rm ext}}$ is orthogonal for all
	$t\in\mathbb R$. Hence the induced operators $\mathcal T_t$ are bounded for
	all $t\in\mathbb R$, and so are their representing matrices $T_t$. Since
	$T_t=e^{t\mathbf N_r}$, no eigenvalue of $\mathbf N_r$ can have nonzero real
	part. Otherwise $e^{t\mathbf N_r}$ would grow exponentially either as
	$t\to+\infty$ or as $t\to-\infty$. Thus
	$\sigma(\mathbf N_r)\subset i\mathbb R$.
	
	Finally, since $S$ is Hurwitz, every $\lambda\in\sigma(S)$ satisfies
	$\operatorname{Re}\lambda<0$.
	In the nilpotent case, every $\mu\in\sigma(\mathbf N_r)$ equals $0$; in the
	skew-symmetric case, every $\mu\in\sigma(\mathbf N_r)$ has zero real part.
	Therefore, in both cases,
	$\operatorname{Re}(\lambda-\mu)<0$.
	Thus every eigenvalue of $B_r$ has negative real part. In particular,
	$0\notin\sigma(B_r)$, and $B_r$ is invertible.
\end{proof}

For any finite-dimensional normed space $Y$, we use the notation
$$C_b(\mathbb R,Y):=\left\{f\in C(\mathbb R,Y):\|f\|_\infty:=\sup_{t\in\mathbb R}\|f(t)\|<\infty\right\},$$
and
$$C_b^1(\mathbb R,Y):=\left\{f\in C^1(\mathbb R,Y):\|f\|_\infty+\|\dot f\|_\infty<\infty\right\}.$$
Here
$\|\dot f\|_\infty:=\sup_{t\in\mathbb R}\|\dot f(t)\|$. We equip
$C_b^1(\mathbb R,Y)$ with the norm
$\|f\|_{C_b^1}:=\|f\|_\infty+\|\dot f\|_\infty$.
When $Y=\mathbb R$, we simply write $C_b(\mathbb R)$ and
$C_b^1(\mathbb R)$.  

\begin{lemma}\label{lem:bounded-complete-coefficient}
	Assume $d\in C_b(\mathbb R,\mathbb R^{M_r})$.  Then equation
	\eqref{c8} has a unique bounded solution  on  $\mathbb R$
	\begin{equation}\label{eq:bounded-complete-coefficient}
	c_r^*(t)=\int_{-\infty}^{t}\exp\!\left(B_r\int_s^tk(\tau)\,\dd\tau\right)d(s)\,\dd s,
\end{equation}
	which belongs to $C_b^1(\mathbb R,\mathbb R^{M_r})$.  Consequently,
	the degree-$r$ equation \eqref{c4} has a unique bounded
 homogeneous polynomial solution on $\mathbb{R}$.
\end{lemma}
\begin{proof}
	By Lemma~\ref{lem:homological-matrix-hurwitz}, $B_r$ is Hurwitz.  Hence
	there are constants $M_r\ge1$ and $\omega_r>0$ such that
	$$\|e^{\xi B_r}\|\le M_re^{-\omega_r\xi},\qquad\xi\ge0.$$
	For $t\ge s$, condition \eqref{c1} gives
	$K(t,s):=\int_s^tk(\tau)\,\dd\tau\ge\underline k(t-s)$, and therefore
	$$\left\|e^{K(t,s)B_r}\right\|\le M_re^{-\omega_r\underline k(t-s)}.$$
	Thus the integral in \eqref{eq:bounded-complete-coefficient} converges
	absolutely for every $t\in\mathbb R$, and its norm is bounded uniformly in
	$t$, with
	$$\|c_r^*\|_\infty\le\frac{M_r}{\omega_r\underline k}\,\|d\|_\infty.$$
	Differentiation
	under the integral sign gives
	$\dot c_r^*=k(t)B_rc_r^*+d(t)$. The right hand side is bounded and
	continuous, so $c_r^*\in C_b^1(\mathbb R,\mathbb R^{M_r})$.  Hence
	$c_r^*$ is a bounded solution of \eqref{c8}.

	For uniqueness, let $c_1,c_2$ be any two bounded complete solutions of
	\eqref{c8}. Set $z:=c_1-c_2$.
	Then $z$ is defined and bounded on $\mathbb R$ and satisfies the
	homogeneous equation
	$$\dot z(t)=k(t)B_rz(t).$$
	Fix $t_0\in\mathbb R$.  For every $T<t_0$, the evolution formula gives
	$$z(t_0)=e^{K(t_0,T)B_r}z(T).$$
	Writing $L:=\sup_{t\in\mathbb R}\|z(t)\|<\infty$, the preceding
	exponential estimate yields
	$$\|z(t_0)\|\le M_re^{-\omega_r\underline k(t_0-T)}\|z(T)\|\le M_rLe^{-\omega_r\underline k(t_0-T)}.$$
	Letting $T\to-\infty$ makes the right-hand side tend to zero, so
	$z(t_0)=0$.  Since $t_0$ is chosen arbitrarily, $z\equiv0$ on $\mathbb R$, and
	therefore $c_1=c_2$. This
	proves uniqueness within the bounded complete class.
\end{proof}

\begin{corollary}\label{cor:bounded-jet-recursion}
	For every $2\le r\le\rho$, the coefficient $h_r(t,u)$ can be solved
	recursively as the unique bounded complete degree-$r$ polynomial.  The
	coefficients through order $r$ depend only on the Taylor coefficients of
	$P,Q$ through order $r$.
\end{corollary}
\begin{proof}
	For $r=2$, uniform boundedness of the derivatives of $P,Q$ makes the
	forcing vector $d_2$ bounded and continuous.  Lemma
	\ref{lem:bounded-complete-coefficient} gives the unique bounded $h_2$.
	Inductively, $H_r$ contains only the already determined
	$h_2,\ldots,h_{r-1}$ and the Taylor coefficients of $P,Q$ through order
	$r$; hence its coefficient vector $d_r$ is again bounded and continuous.
	The same lemma completes the induction.  Terms of order $r+1$ and higher
	do not occur in $H_r$.
\end{proof}

By Corollary~\ref{cor:bounded-jet-recursion}, for each $2\le r\le\rho$, equation (3.12) admits a unique bounded homogeneous polynomial solution $h_r(t,u)$, which can be determined recursively. By Proposition 3.6 of \cite{PotzscheRasmussen2006}, all integral manifolds possess the same Taylor expansion with respect to the state-space variable $u$ up to order $\rho$.
Thus, although the graph function $h$ representing the center manifold may not be unique, its Taylor expansion up to order $\rho$ is uniquely determined. In the remainder of the paper, we use $h_r$ to denote the uniquely determined homogeneous term of degree $r$ in this Taylor expansion.

\begin{theorem}
	\label{thm:reduced-vector-field-jet} By a suitable choice of
	the higher-order terms $P(t,x_c,\varepsilon,v)$ and $Q(t,x_c,\varepsilon,v)$, the Taylor coefficients of degrees $2$ through $\rho$ with respect to $(x_c,\varepsilon)$ of the parameterized reduced vector field $R(t,x_c,\varepsilon)$ can be prescribed arbitrarily within the class satisfying $R(t,0,\varepsilon)=0$, provided their coefficient functions belong to $C_b^1(\mathbb R)$.
\end{theorem}

\begin{proof}
		Equivalently, given any $E^c$-valued polynomial vector field
	$$\widehat R(t,x_c,\varepsilon)=k(t)\Omega x_c+\sum_{r=2}^{\rho}\widehat a_r(t,x_c,\varepsilon),$$
	where each $\widehat a_r(t,\cdot)$ is homogeneous of degree $r$ in $(x_c,\varepsilon)$ with coefficient functions in $C_b^1(\mathbb R)$, and satisfies $\widehat a_r(t,0,\varepsilon)=0$, one can choose
	the terms of order $2$ through $\rho$ of the original nonlinearities so that the parameterized reduced vector field has the prescribed Taylor expansion through order $\rho$.
	
	We start from the identity
	$$R(t,u)=k(t)\Omega_cu+P\bigl(t,u,h(t,u)\bigr).$$
	In particular,
	$$R(t,0)=0,\qquad D_uR(t,0)=k(t)\Omega_c.$$
	Hence the constant and linear terms of $R$ are fixed, and only the terms of order at
	least two are available for adjustment. Write
	$$P(t,u,v)=\sum_{r=2}^\rho P_r(t,u,v)+\mathcal{O}(|(u,v)|^{\rho+1}),\qquad Q(t,u,v)=\sum_{r=2}^\rho Q_r(t,u,v)+\mathcal{O}(|(u,v)|^{\rho+1}),$$
	where $P_r$ and $Q_r$ are homogeneous of degree $r$ in $(u,v)$. Likewise, write
	$$h(t,u)=\sum_{r=2}^\rho h_r(t,u)+o(|u|^{\rho}),\qquad R(t,u)=k(t)\Omega_cu+\sum_{r=2}^\rho a_r(t,u)+o(|u|^{\rho}),$$
	where $h_r$ and $a_r$ are homogeneous of degree $r$ in $u$.
	
	Notice that the bounded complete Taylor coefficients of the graph
	function $h$ are obtained uniquely and recursively from the invariance
	equation by Corollary~\ref{cor:bounded-jet-recursion}. For fixed linear part
	$(k(t)\Omega,k(t)S)$ and for each integer $j\ge2$, the coefficients
	$h_2,\dots,h_j$
	are uniquely determined within the bounded complete class by the coefficients of
	$$P_2,\dots,P_j,\qquad Q_2,\dots,Q_j.$$
	Therefore, if the vector field is modified only in order $j+1$,  the Taylor expansion of
	$h$ through order $j$ remains unchanged.
	
	We proceed by induction. For $r=2$, since $h(t,u)=O(|u|^2)$, the quadratic part of $R$ comes only from the quadratic part of $P$ evaluated at $v=0$. Thus
	$$a_2(t,u)=P_2(t,u,0).$$
	Hence the second-order coefficient of the reduced vector field can be prescribed
	arbitrarily by choosing the pure center quadratic term of the original system. We set
	$$P_2(t,u,0):=\widehat a_2(t,u).$$
	Then $a_2(t,u)=\widehat a_2(t,u)$.
	
	Assume now that, for some integer $p$ with $2\le p<\rho$, the terms of the original
	system have already been chosen in such a way that
	$$a_r(t,u)=\widehat a_r(t,u),\qquad r=2,\dots,p.$$
	We show that the terms of order $p+1$ can be chosen so as to make $a_{p+1}(t,u)=\widehat a_{p+1}(t,u)$.
	Modify the original system by adding homogeneous terms of degree $p+1$:
	$$\check P(t,u,v)=P(t,u,v)+\Delta P_{p+1}(t,u,v),$$
	$$\check Q(t,u,v)=Q(t,u,v)+\Delta Q_{p+1}(t,u,v),$$
	where $\Delta P_{p+1}$ and $\Delta Q_{p+1}$ are homogeneous of degree $p+1$ in
	$(u,v)$. Let
	$$\check h(t,u)=\sum_{r=2}^\rho\check h_r(t,u)+o(|u|^{\rho}),\qquad\check R(t,u)=k(t)\Omega_cu+\sum_{r=2}^\rho\check a_r(t,u)+o(|u|^{\rho})$$
	denote the corresponding center manifold and reduced vector field. Therefore
	$$\check h_r=h_r,\qquad r=2,\dots,p,$$
	and consequently
	$$\check a_r=a_r,\qquad r=2,\dots,p.$$
	Thus the only coefficient of the reduced vector field that can change at this stage is the
	one of degree $p+1$. Since
	$$\check R(t,u)=k(t)\Omega_cu+\check P\bigl(t,u,\check h(t,u)\bigr),$$
	its homogeneous part of degree $p+1$ is given by
	$$\check a_{p+1}(t,u)=\bigl[\check P\bigl(t,u,\check h(t,u)\bigr)\bigr]_{p+1}.$$
	Because $\check h(t,u)=O(|u|^2)$,
	every monomial in $\Delta P_{p+1}(t,u,v)$ that contains at least one factor of $v$
	becomes, after the substitution $v=\check h(t,u)$, a term of degree at least $p+2$
	in $u$. Hence these terms do not contribute to the coefficient of degree
	$p+1$ in the reduced vector field. The only part of $\Delta P_{p+1}$ that contributes in
	degree $p+1$ is the part
	$$\Delta P_{p+1}^c(t,u):=\Delta P_{p+1}(t,u,0).$$
	It follows that $\check a_{p+1}(t,u)=a_{p+1}(t,u)+\Delta P_{p+1}^c(t,u)$.
	Therefore the $(p+1)$-st coefficient of the reduced
	vector field can take any prescribed value by choosing
	$$\Delta P_{p+1}^c(t,u):=\widehat a_{p+1}(t,u)-a_{p+1}(t,u).$$
	With this choice we obtain
	$\check a_{p+1}(t,u)=\widehat a_{p+1}(t,u)$,
	while the coefficients of orders $2,\dots,p$ remain unchanged.
	
	Thus we extend the agreement from order $p$ to order $p+1$. By induction,
	the coefficients of the reduced vector field can be chosen so that
	$$a_r(t,u)=\widehat a_r(t,u),\qquad r=2,\dots,\rho.$$

    This completes the proof by induction.
\end{proof}

\begin{proof}[Proof of Theorem~\ref{thm:structural-reduction}]
	Write
	$$K_A=\frac{A-A^T}{2},\qquad A=F_A+K_A.$$
	The matrix $K_A$ is skew-symmetric. Hence, if $\widehat B$ is a skew-symmetric matrix constructed for the symmetric part $F_A$, then the choice
	$B=\widehat B-K_A$
	is also skew-symmetric and satisfies $A+B=F_A+\widehat B$. 

	Let $\lambda_1\ge\cdots\ge\lambda_n$ be the eigenvalues of $F_A$. Proposition~\ref{prop:center-stable-decomposition} therefore provides a skew-symmetric $\widehat B$ and an invariant splitting
	$$\mathbb R^n=E^c\oplus E^s,\qquad\dim E^c=m,$$
	such that $\Omega=(F_A+\widehat B)|_{E^c}$ is nilpotent and $S=(F_A+\widehat B)|_{E^s}$ is Hurwitz. In the case where $F_A$ has precisely an $m$-dimensional zero eigenspace, Corollary~\ref{cor:skew-hurwitz-decomposition-zero-eigenvalue} instead permits $\Omega$ to be chosen skew-symmetric while the complementary block remains Hurwitz.

	With the preceding choice of $B$, the variational equation is $\dd x/\dd t=k(t)(A+B)x$. The bounds in \eqref{c1}, together with the nilpotent or skew-symmetric estimate on $E^c$ and the Hurwitz estimate on $E^s$, give the exponential trichotomy, with center space $E^c$, stable space $E^s$, and $E^u=\{0\}$. Consequently, a local parameterized center manifold of dimension $m$ is constructed. Its graph may depend on $t$, and its reduced vector field satisfies
	$$R(t,0,\varepsilon)=0,\qquad D_{x_c}R(t,0,0)=k(t)\Omega.$$

	Finally, fix any finite $\rho\ge2$. For $2\le r\le\rho$, let $\widehat a_r(t,x_c,\varepsilon)$ be arbitrary $E^c$-valued homogeneous polynomials of degree $r$, with coefficient functions in $C_b^1(\mathbb R)$ and subject to $\widehat a_r(t,0,\varepsilon)=0$. Theorem~\ref{thm:reduced-vector-field-jet} allows the higher-order terms of $P(t,x_c,\varepsilon,v)$ and $Q(t,x_c,\varepsilon,v)$ to be chosen recursively so that the degree-$r$ Taylor coefficient of the reduced vector field is exactly $\widehat a_r(t,x_c,\varepsilon)$ for every $2\le r\le\rho$. Lemma~\ref{lem:realizability-projected-jets} realizes the corresponding compatible projected finite jet by $f_s(t,u,\varepsilon)$ in the original variables. This completes the proof.
\end{proof}

\subsection{Nilpotent Shilnikov realization}

 Under the factorization \eqref{force-factorization}, write $u=(x_c,\varepsilon)\in E^c\times\mathcal U_\varepsilon$ and $v\in E^s$. In the corresponding $(u,v)$ coordinates the nonlinear terms have the same scalar time factor,
 $$P(t,u,v)=k(t)P_0(u,v),\qquad Q(t,u,v)=k(t)Q_0(u,v).$$
 We now look for a center-manifold graph whose Taylor coefficients are independent of $t$.
 Suppose inductively that $h_2,\ldots,h_{r-1}$ have already been chosen independent of $t$.
 Then the $r$-th homogeneous term in the right-hand side of \eqref{c4} has the form
 $$H_r(t,u)=k(t)H_{r,0}(u),$$
 because both $P$ and $Q$ contain the common factor $k(t)$.
 If $h_r(t,u)=h_r(u)$, then the coefficient vector $c$ of $h_r$ satisfies
 $$0=k(t)B_rc+k(t)d_0,$$
 where $d_0$ is the coefficient vector of $H_{r,0}$.
 Since $k(t)>0$ and $B_r$ is invertible by the invertibility argument above, this equation has the unique solution
 $$c=-B_r^{-1}d_0.$$
 This is also the value obtained from
 \eqref{eq:bounded-complete-coefficient}. Hence it is the unique bounded
 coefficient, independent of $t$, and the construction proceeds
 recursively.
 Apply the recursive construction in Theorem~\ref{thm:reduced-vector-field-jet} within the factorized class, using $\widehat a_r(t,x_c,\varepsilon)=k(t)\widehat a_{r,0}(x_c,\varepsilon)$ and corrections of the form $\Delta P_r(t,u,v)=k(t)\Delta P_{r,0}(u,v)$. The preceding recursion then gives a local center manifold with graph function independent of $t$, $v=h(x_c,\varepsilon)$, where $h(0,\varepsilon)=0$ and $D_{x_c}h(0,0)=0$, and on it
$$R(t,x_c,\varepsilon)=k(t)(\Omega x_c+R_0(x_c,\varepsilon)),\qquad R_0(0,\varepsilon)=0,\qquad R_0(x_c,\varepsilon)=O(\|(x_c,\varepsilon)\|^2).$$

For each fixed $\varepsilon$, before showing chaotic behavior in the reduced vector field $R(t,x_c,\varepsilon)$, we apply the global time reparametrization
$$s=s(t):=\int_{0}^tk(\tau)\,\dd\tau.$$
By \eqref{c1}, $s(t)$ is a strictly monotone $C^1$ diffeomorphism from
$\mathbb R$ onto $\mathbb R$. Consider the nonautonomous system
$$\frac{\dd x}{\dd t}=k(t)R_1(x)$$
and the reparametrized autonomous system
$$\frac{\dd x}{\dd s}=R_1(x).$$
If $x(t)$ is a solution of the former and we set $\widetilde x(s)=x(t(s))$, then
$$\frac{\dd\widetilde x}{\dd s}=R_1(\widetilde x).$$
Conversely, every solution of the autonomous system gives rise, through the inverse change
$t=t(s)$, to a solution of the original system. Hence both systems have exactly the same
phase curves in the state space and only the speed along these curves is changed.

As a consequence, all geometric orbit structures are preserved under the time
reparametrization, in particular equilibria, homoclinic and heteroclinic connections. Homoclinic and heteroclinic connections, as well as global bifurcations, are mechanisms for generating rich and complex dynamical behavior \cite{DiazGelfertRams2022,DiazRochaViana1996,GorodetskiKaloshin2009,
	MatheusMoreiraPalis2018,Melbourne1989,MoraViana1993}. In particular, horseshoes, symbolic dynamics, and the associated chaotic invariant sets produced by a Shilnikov homoclinic configuration are preserved. Thus, once chaos is established for the reduced autonomous system in the $s$-time scale, the same horseshoe-type chaotic dynamics is present in the original system.

After the time reparametrization, the dynamics on the center manifold is denoted by
\begin{equation}\label{eq:Y}
	Y(x_c,\varepsilon)=\Omega x_c+R_0(x_c,\varepsilon).
\end{equation}

Recall briefly the Shilnikov mechanism. Let $p$ be a saddle-focus equilibrium of a smooth vector field
$$\frac{\dd x}{\dd s}=X(x),\qquad x\in\mathbb{R}^3.$$
Assume that the eigenvalues of $DX(p)$ are
$\lambda$ and $-\rho\pm i\omega$,
where $\lambda>0$, $\rho>0$, and $\omega\ne0$. Then $\dim W^u(p)=1$ and $\dim W^s(p)=2$. A Shilnikov homoclinic orbit is a nontrivial orbit, denoted by $\Gamma_h$, lying in the intersection of the stable and unstable manifolds of the saddle-focus equilibrium.  It is known that every sufficiently small neighborhood of $\Gamma_h$ contains infinitely many saddle periodic orbits if $\lambda>\rho$ and chaotic behavior occurs in this case. Here, chaotic behavior in the sense of Robert L. Devaney is characterized by sensitive dependence on initial conditions, topological transitivity, and a dense set of periodic points.

Consider a sufficiently smooth three-parameter unfolding of a vector field on $\mathbb R^3$ of the nilpotent singularity, whose nilpotent linear part is conjugate to
$$J=\begin{pmatrix}
	0&1&0\\
	0&0&1\\
	0&0&0
\end{pmatrix}.$$
As stated in \cite{IbanezRodriguez2005}, after a coordinate change and the introduction of rescaled parameters $(c,n,u)$, the normal form can be reduced to
\begin{equation}\label{nilpotent}
	y\frac{\partial}{\partial x}+z\frac{\partial}{\partial y}+(c^2-y+nz-\frac{x^2}{2}+\mathcal R_N(u))\frac{\partial}{\partial z},
\end{equation}
where
$|\mathcal R_N(u)|\le C_N|u|$
for a constant $C_N>0$.

By analyzing the intersection of the two-dimensional stable manifold
and the one-dimensional unstable manifold, \cite{IbanezRodriguez2005}
proved that the heteroclinic connection breaks and a homoclinic loop
is formed.

\begin{theorem}[\cite{IbanezRodriguez2005}]
	\label{thm:nilpotent-shilnikov}
	There exist parameter values arbitrarily close to $(c,n,u)=(c^*,0,0)$, such that the unfolded vector field \eqref{nilpotent} has a Shilnikov homoclinic orbit to a saddle-focus equilibrium. 
\end{theorem}

\begin{proof}[Proof of Corollary~\ref{cor:three-dimensional-shilnikov}]
Consider \eqref{eq:Y}. Since $\dim E^c=3$ and $\operatorname{rank}\Omega=2$, the nilpotent center block $\Omega$ is conjugate to the single Jordan block $J$. By Theorem~\ref{thm:reduced-vector-field-jet}, the nonlinear terms and three effective parameter directions can be chosen so that this autonomous reduced family is a generic three-parameter unfolding of the codimension-three nilpotent singularity. Following \cite{NijholtPereiraQueirozTuraev2023}, suitable coordinate and parameter changes, together with a family blow-up, bring \eqref{eq:Y} into the form \eqref{nilpotent}.  Theorem~\ref{thm:nilpotent-shilnikov} then gives parameter values for which the autonomous reduced system has a Shilnikov homoclinic orbit.
	
Finally, \eqref{c1} implies that $$s(t)=\int_0^tk(\tau)\,\dd\tau$$ is an orientation-preserving diffeomorphism of $\mathbb R$.  Hence the inverse time change preserves both the forward and backward branches of the homoclinic orbit and changes only its time parametrization. The flow of $R(t,x_c,\varepsilon)$ therefore has the same Shilnikov homoclinic loop and the associated chaotic dynamics.
\end{proof}

\subsection{Skew-symmetric case}

Consider \eqref{eq:Y}. In this subsection we assume that at the zero parameter value, its linear part is
orthogonally conjugate to
$$B_p=\begin{pmatrix}
	0&p&0\\
	-p&0&0\\
	0&0&0
\end{pmatrix},\qquad p>0,$$
and hence has spectrum $\{0,\pm ip\}$. 

We briefly review the relevant results from \cite{BroerVegter1984}.
For a $C^\infty$ two-parameter Hopf-zero unfolding, after a change of coordinates and a rescaling of the two
unfolding parameters, the authors studied equations of the following form
\begin{samepage}
\begin{equation}\label{smooth-hopf-zero-normal-form-relaxed}
	\begin{aligned}
	\frac{\dd\varphi}{\dd s}&=1,\\
	\frac{\dd r}{\dd s}&=r(\nu-a_1z)+r\mathcal H_r(r^2,z,\mu,\nu),\\
	\frac{\dd z}{\dd s}&=-\mu+b_1r^2+b_2z^2+b_3\mu z+b_4\nu z+b_5\mu\nu+\mathcal H_z(r^2,z,\mu,\nu),
\end{aligned}
\end{equation}
where $\mathcal H_r$ and $\mathcal H_z$ denote the higher-order terms in the $r$- and $z$-equations, respectively. 
\end{samepage}
Here $(\mu,\nu)$ are the two unfolding parameters,
$(r,\varphi,z)$ are cylindrical phase variables with $r\geq0$, and
$a_1,b_1,\ldots,b_5$ are real normal-form coefficients.
Under the assumption
\begin{equation}\label{broer-vegter-open-conditions-relaxed}
	0<a_1<2b_2,\qquad b_1>0,\qquad b_4>0,
\end{equation}
we recall the following consequence of Lemma~A in
\cite{BroerVegter1984}.  For sufficiently small $\mu^2+\nu^2$ and
$\mu>0$, \eqref{smooth-hopf-zero-normal-form-relaxed} has two
hyperbolic saddle points $p_\pm(\mu,\nu)$ on the $z$-axis, with $z$-coordinate
$$z_\pm(\mu,\nu)=\pm\sqrt{\frac{\mu}{b_2}}+O\!\left(\sqrt{\mu^2+\nu^2}\right).$$
The interval between $p_+$ and $p_-$ on the $z$-axis is a
one-dimensional heteroclinic connection.  Moreover, there exists a
smooth curve
$\Gamma$
such that, for $(\mu,\nu)\in\Gamma$, the two-dimensional invariant
manifolds of $p_+$ and $p_-$ coincide. 

In what follows $W^u_j(p)$ and $W^s_j(p)$ denote, respectively, a $j$-dimensional
unstable or stable manifold of the equilibrium $p$.
Geometrically, these coincident two-dimensional invariant manifolds
form an invariant globe with $p_+$ and $p_-$ as its two poles, while
the one-dimensional connection lies on the rotation axis inside the
globe.  More precisely,
$$W^u_1(p_+)=W^s_1(p_-),\qquad W^u_2(p_-)=W^s_2(p_+).$$

On $\Gamma$, two coincident invariant manifolds appear. A suitable $C^\infty$-flat perturbation breaks the axial and globe connections while leaving every finite Taylor coefficient unchanged. The authors showed that the flat terms can be chosen so that the perturbed family has an immersed curve of Shilnikov bifurcation parameters accumulating at the Hopf-zero point.

\begin{theorem}[Broer--Vegter \cite{BroerVegter1984}]
	\label{thm:smooth-hopf-zero-shilnikov-relaxed}
Let $\overline X_{\mu,\nu}$ denote the full $C^\infty$
rotationally symmetric family given by
\eqref{smooth-hopf-zero-normal-form-relaxed}.  Assume that
its second-order Taylor coefficients satisfy the open
conditions \eqref{broer-vegter-open-conditions-relaxed}, and let
$\Gamma$ be the parameter curve obtained above.
	
Then there exists a $C^\infty$ family of perturbations $p_{\mu,\nu}$, flat at the Hopf--zero point, such that $X_{\mu,\nu}=\overline X_{\mu,\nu}+p_{\mu,\nu}$ has a sequence of Shilnikov bifurcations at parameter values $(\mu_j,\nu_j)\in\Gamma$, which accumulate at $(\mu,\nu)=(0,0)$.
	
More precisely, there exists an immersed $C^\infty$ curve $C$ of Shilnikov bifurcation parameter values contained in a wedge-shaped neighborhood $W$ of $\Gamma\setminus\{(0,0)\}$ whose width is flat as $\mu\downarrow0$. Moreover, $\partial C\cap W=\varnothing$, and $C$ intersects $\Gamma$ transversally at the parameter values $(\mu_j,\nu_j)$. A transverse crossing of $C$ produces a Shilnikov bifurcation.
\end{theorem}

\begin{remark}
	The preceding Broer--Vegter proof is intrinsically $C^\infty$: it creates the Shilnikov loop by choosing a nonzero flat remainder. Therefore this argument does not apply directly to the analytic category. The analytic Hopf-zero problem was treated by  \cite{BaldomaIbanezSeara2020}. By analyzing the exponentially small splitting of invariant manifolds, the authors prove the existence of Shilnikov homoclinic orbits at parameter values accumulating at the Hopf-zero point under analytic genericity hypotheses.
\end{remark}

Consider vector field \eqref{eq:Y} satisfying $Y(0,\varepsilon)=0$. When $\Omega=B_p$, after cylindrical coordinate transformation and rescaling we obtain the system 
\begin{equation}\label{eq:persistent-hopf-zero-normal-form}
	\begin{aligned}
	\frac{\dd\varphi}{\dd s}&=1,\\
	\frac{\dd r}{\dd s}&=r\Bigl(a_{11}\varepsilon_1+a_{12}\varepsilon_2-a_1\zeta+\widetilde{\mathcal H}_r(r^2,\zeta,\varepsilon_1,\varepsilon_2)\Bigr),\\
	\frac{\dd\zeta}{\dd s}&=(a_{21}\varepsilon_1+a_{22}\varepsilon_2)\zeta+b_1r^2+b_2\zeta^2+\widetilde{\mathcal H}_z(r^2,\zeta,\varepsilon_1,\varepsilon_2),
\end{aligned}
\end{equation}
where $\widetilde{\mathcal H}_r$ and
$\widetilde{\mathcal H}_z$ are $C^\infty$ functions satisfying
$$\widetilde{\mathcal H}_r=O\!\left(\left\|(r,\zeta,\varepsilon_1,\varepsilon_2)\right\|^2\right),\qquad\widetilde{\mathcal H}_z=O\!\left(\left\|(r,\zeta,\varepsilon_1,\varepsilon_2)\right\|^3\right).$$
 We also impose
$$\widetilde{\mathcal H}_z(0,0,\varepsilon_1,\varepsilon_2)=0,$$
so that $r=\zeta=0$ remains an equilibrium for every parameter
value.

Set
$$A_p:=\begin{pmatrix}
	a_{11}&a_{12}\\
	a_{21}&a_{22}
\end{pmatrix}.$$

\begin{theorem}
	\label{thm:persistent-to-broer-vegter}
	Assume that
	$$\det A_p\neq0,\qquad0<a_1<2b_2,\qquad b_1>0.$$
	For every fixed $b_4>0$, the higher-order terms in
	\eqref{eq:persistent-hopf-zero-normal-form} can be chosen so that
	a parameter-dependent translation and a parameter change
	transform \eqref{eq:persistent-hopf-zero-normal-form} into the  form
	\eqref{smooth-hopf-zero-normal-form-relaxed}.
	Moreover, every Shilnikov parameter value $(\mu_j,\nu_j)$ in Theorem~\ref{thm:smooth-hopf-zero-shilnikov-relaxed}, with
	$\mu_j>0$, has two preimages
	$(\varepsilon_{1,j}^{\pm},\varepsilon_{2,j}^{\pm})$, which tend to $(0,0)$ as $j\to\infty$.
\end{theorem}

\begin{proof}
	First introduce the parameters $(\sigma,\tau)$ by
	\begin{equation}\label{eq:first-persistent-parameter-change}
	\begin{pmatrix}
	\sigma\\
	\tau
\end{pmatrix}=A_p\begin{pmatrix}
	\varepsilon_1\\
	\varepsilon_2
\end{pmatrix}.
\end{equation}
	The assumption $\det A_p\neq0$ makes this an invertible linear
	parameter change. Applying it to \eqref{eq:persistent-hopf-zero-normal-form}, we obtain
	$$\begin{aligned}
	\frac{\dd\varphi}{\dd s}&=1,\\
	\frac{\dd r}{\dd s}&=r\bigl(\sigma-a_1\zeta+\widetilde{\mathcal H}_r(r^2,\zeta,\varepsilon_1,\varepsilon_2)\bigr),\\
	\frac{\dd\zeta}{\dd s}&=\tau\zeta+b_1r^2+b_2\zeta^2+\widetilde{\mathcal H}_z(r^2,\zeta,\varepsilon_1,\varepsilon_2).
\end{aligned}$$
	Fix $b_4>0$ and set
	$D:=2b_2+a_1b_4$.
	Since $0<a_1<2b_2$, we have $D>0$. Define
	\begin{equation}\label{eq:second-persistent-parameter-change}
	h=\frac{\tau-b_4\sigma}{D},\qquad\nu=\sigma+a_1h=\frac{2b_2\sigma+a_1\tau}{D},\qquad\mu=\tau h-b_2h^2.
\end{equation}
	Make the parameter-dependent translation
	$z=\zeta+h$.
	Substitution \eqref{eq:second-persistent-parameter-change} into the radial equation gives
	$$\begin{aligned}
	\frac{\dd r}{\dd s}&=r\Bigl(\sigma+a_1h-a_1z+\widetilde{\mathcal H}_r(r^2,z-h,\sigma,\tau)\Bigr)\\
	&=r\bigl(\nu-a_1z+\widetilde{\mathcal H}_r(r^2,z-h,\sigma,\tau)\bigr),
\end{aligned}$$
	Similarly,
	$$\begin{aligned}
	\frac{\dd z}{\dd s}&=\tau(z-h)+b_1r^2+b_2(z-h)^2+\widetilde{\mathcal H}_z(r^2,z-h,\sigma,\tau)\\
	&=-\mu+b_1r^2+b_2z^2+b_4\nu z+\widetilde{\mathcal H}_z(r^2,z-h,\sigma,\tau),
\end{aligned}$$
	By \eqref{eq:second-persistent-parameter-change}, $\mu$ is a quadratic form in  $(\sigma,\tau)$ and $(h,\nu)$ is linear in $(\sigma,\tau)$. Replacing the parameters $(\sigma,\tau)$ in the higher-order terms by $(\mu,\nu)$ yields for some constants $a_2,b_3,b_5$
	\begin{align*}
	\widetilde{\mathcal H}_r(r^2,z-h,\sigma,\tau)&=a_2\mu+{\mathcal H}_r(r^2,z,\mu,\nu),\\
	\widetilde{\mathcal H}_z(r^2,z-h,\sigma,\tau)&=b_3\mu z+b_5\mu\nu+{\mathcal H}_z(r^2,z,\mu,\nu),
\end{align*}
	where ${\mathcal H}_r(r^2,z,\mu,\nu)=O(\|(r^2,z,\mu,\nu)\|^2)$ and ${\mathcal H}_z(r^2,z,\mu,\nu)=O(\|(r^2,z,\mu,\nu)\|^3)$. We choose the higher-order term such that $a_2=0$.
	Therefore \eqref{eq:persistent-hopf-zero-normal-form} becomes
	$$\begin{aligned}
	\frac{\dd\varphi}{\dd s}&=1,\\
	\frac{\dd r}{\dd s}&=r(\nu-a_1z+\mathcal H_r),\\
	\frac{\dd z}{\dd s}&=-\mu+b_1r^2+b_2z^2+b_3\mu z+b_4\nu z+b_5\mu\nu+\mathcal H_z,
\end{aligned}$$
	which is the form
	\eqref{smooth-hopf-zero-normal-form-relaxed}.

	It remains to determine the preimages of the Shilnikov parameters.
	Solving the relations defining $h$ and $\nu$ for
	$(\sigma,\tau)$ gives
	$$\tau=2b_2h+b_4\nu,\qquad\sigma=\nu-a_1h.$$
	Consequently,
	$$\mu=\tau h-b_2h^2=b_2h^2+b_4\nu h.$$
	Let $(\mu_j,\nu_j)$ be a Shilnikov parameter value supplied by
	Theorem~\ref{thm:smooth-hopf-zero-shilnikov-relaxed}. Since
	$\mu_j>0$, the equation
	$$b_2h^2+b_4\nu_jh-\mu_j=0$$
	has two distinct real solutions
	$$h_j^\pm=\frac{-b_4\nu_j\pm\sqrt{b_4^2\nu_j^2+4b_2\mu_j}}{2b_2}.$$
	Define
	$$\tau_j^\pm=2b_2h_j^\pm+b_4\nu_j=\pm\sqrt{b_4^2\nu_j^2+4b_2\mu_j}$$
	and
	$$\sigma_j^\pm=\nu_j-a_1h_j^\pm.$$
	Finally, use the inverse of
	\eqref{eq:first-persistent-parameter-change} and set
	$$\begin{pmatrix}
	\varepsilon_{1,j}^{\pm}\\
	\varepsilon_{2,j}^{\pm}
\end{pmatrix}=A_p^{-1}\begin{pmatrix}
	\sigma_j^\pm\\
	\tau_j^\pm
\end{pmatrix}.$$
	Thus
	$(\varepsilon_{1,j}^{\pm},\varepsilon_{2,j}^{\pm})$ are two preimages of
	$(\mu_j,\nu_j)$.
	
	Finally, as $j\to\infty$, $(\mu_j,\nu_j)\to(0,0)$ implies
	$$h_j^\pm\to0,\qquad\sigma_j^\pm\to0,\qquad\tau_j^\pm\to0.$$
	Since $A_p$ is invertible,
	$(\varepsilon_{1,j}^{\pm},\varepsilon_{2,j}^{\pm})\to(0,0)$.
\end{proof}

\begin{proof}[Proof of Corollary \ref{cor:persistent-skew-shilnikov}]
	Consider \eqref{eq:Y}. Since $\Omega$ is orthogonally conjugate to
	$B_p$, after a change of coordinates, we may assume that the dynamical system becomes \eqref{eq:persistent-hopf-zero-normal-form}. By
	Theorem~\ref{thm:reduced-vector-field-jet}, the Taylor terms
	of degree two and three can be chosen to satisfy the hypotheses
	of Theorem~\ref{thm:persistent-to-broer-vegter}. 
	Theorem~\ref{thm:smooth-hopf-zero-shilnikov-relaxed}, together
	with Theorem~\ref{thm:persistent-to-broer-vegter}, then gives
	parameter values
	$(\varepsilon_{1,j}^{\pm},\varepsilon_{2,j}^{\pm})$ tending to $(0,0)$ as $j\to\infty$
	for which \eqref{eq:Y} has a Shilnikov
	homoclinic loop.
	
	Finally, \eqref{c1} implies that
	$$s(t)=\int_0^tk(\tau)\,\dd\tau$$
	is an orientation-preserving diffeomorphism of $\mathbb R$.
	The inverse time change preserves the homoclinic phase curve
	and changes only its time parametrization. Hence the flow of  $R(t,x_c,\varepsilon)$ has the same Shilnikov loop and the associated chaotic
	dynamics.
\end{proof}

\begin{spacing}{1.2}
	\noindent\textbf{Declaration of competing interest}. 
	
	The authors declare no conflicts of interest.\\

	\noindent\textbf{Acknowledgements}
	
	The work of L.H. Chen is partially supported by Guangxi University International Cooperation Exchange and Study Program.
	The work of W.P. Yan is partially supported by NSFC. No 12671189 and the special foundation for Guangxi Ba Gui Scholars. 
	The work of Y. Li is partially supported by NSFC (Grant Nos. 12531009, 12471183, 12071175),
	National Basic Research Program of China (Grant No. 2013CB834100), and the Natural Science Foundation
	of Jilin Province (Grant No. 20200201253JC).\\
	
	\noindent\textbf{Data availability}
	
	Data sharing not applicable to this article as no data sets were generated or analyzed during the current
	study.
\end{spacing}

\end{document}